\documentclass{article}

     \PassOptionsToPackage{numbers,sort, compress}{natbib}

    \usepackage[preprint]{neurips_2020}

\usepackage[utf8]{inputenc} 
\usepackage[T1]{fontenc}    
\usepackage{hyperref}       
\usepackage{url}            
\usepackage{booktabs}       
\usepackage{amsfonts}       
\usepackage{nicefrac}       
\usepackage{microtype}      

\usepackage{wrapfig}

\usepackage[font=small]{caption}

\usepackage{subcaption}
\usepackage{graphicx}
\usepackage{caption}
\usepackage{amsfonts}
\usepackage{hyperref}
\hypersetup{colorlinks=true,citecolor=blue}
\usepackage{algorithm}
\usepackage{algpseudocode}
\usepackage{framed}
\usepackage{footnote}
\usepackage{amsmath}  
\usepackage{amssymb}   
\usepackage{amsthm}  
\usepackage{comment}  
\usepackage{enumitem}
\usepackage[colorinlistoftodos,prependcaption,textsize=tiny]{todonotes}
\numberwithin{equation}{section}
\usepackage{xcolor}

\usepackage{soul}

\usepackage{array}
\usepackage{arydshln}
\newtheorem{theorem}{Theorem}[section]

\newtheorem{lemma}[theorem]{Lemma}
\newtheorem{remark}[theorem]{Remark}
\newtheorem{assumption}[theorem]{Assumption}

\newtheorem{definition}[theorem]{Definition}

\newtheorem{innercustomthm}{Theorem}
\newenvironment{customthm}[1]
  {\renewcommand\theinnercustomthm{#1}\innercustomthm}
  {\endinnercustomthm}

\newtheorem{innercustomlemma}{Lemma}
\newenvironment{customlemma}[1]
  {\renewcommand\theinnercustomlemma{#1}\innercustomlemma}
  {\endinnercustomlemma}

\newtheorem{innercustomassum}{Assumption}
\newenvironment{customassum}[1]
  {\renewcommand\theinnercustomassum{#1}\innercustomassum}
  {\endinnercustomassum}  
\usepackage{bbm}
\newcommand{\V}{\tilde{V}}
\newcommand{\R}{\mathbb{R}}
\newcommand{\A}{\mathcal{A}}
\newcommand{\Ap}{A^{\perp}}
\newcommand{\gp}{g^{\perp}}

\newcommand{\eqdef}{:=}
\newcommand{\SONIA}{\texttt{SONIA}}

\title{SONIA: A Symmetric Blockwise Truncated Optimization Algorithm}

\author{%
 Majid Jahani \\
  Lehigh University\\
  Bethlehem, PA 18015 \\
  \url{maj316@lehigh.edu} \\
   \And
  Mohammadreza Nazari \\
  SAS Institute\\
  Cary, NC 27513 \\
  \url{mrza.nazari@gmail.com} \\
   \AND
    Rachael Tappenden \\
  University of Canterbury\\
  Christchurch 8041, New Zealand \\
  \url{rachael.tappenden@canterbury.ac.nz} \\
  \And
  Albert S. Berahas \\
  Lehigh University\\
  Bethlehem, PA 18015 \\
  \url{albertberahas@gmail.com} \\
  \And
  Martin Tak\' a\v c \\
  Lehigh University\\
  Bethlehem, PA 18015 \\
  \url{Takac.MT@gmail.com} \\
}

\begin{document}

\setlength{\abovedisplayskip}{3pt}
\setlength{\belowdisplayskip}{3pt}

\maketitle
\begin{abstract}
This work presents a new algorithm for empirical risk minimization. The algorithm bridges the gap between first- and second-order methods by computing a search direction that uses a second-order-type update in one subspace, coupled with a scaled steepest descent step in the orthogonal complement. To this end, partial curvature information is incorporated to help with ill-conditioning, while simultaneously allowing the algorithm to scale to the large problem dimensions often encountered in machine learning applications. Theoretical results are presented to confirm that the algorithm converges to a stationary point in both the strongly convex and nonconvex cases. A stochastic variant of the algorithm is also presented, along with corresponding theoretical guarantees. Numerical results confirm the strengths of the new approach on standard machine learning problems. 
\end{abstract}




\vspace{-5pt}
\section{Introduction}
\vspace{-5pt}

This paper presents a novel optimization algorithm for empirical risk minimization:
\begin{align}   \label{eq:prob}
    {\min}_{w \in \mathbb{R}^d} F(w) := \tfrac{1}{n} {\textstyle\sum}_{i=1}^n f(w;x^i,y^i) = \tfrac{1}{n} {\textstyle\sum}_{i=1}^n f_i(w),
\end{align}
where $\{(x^i,y^i)\}_{i=1}^n$ are training examples (observations), and $f_i: \mathbb{R}^d \rightarrow \mathbb{R}$ is the composition of a prediction function (parameterized by $w \in \mathbb{R}^d$) and a loss function associated with the $i$th training observation (sample). Problems of the form \eqref{eq:prob} arise in a wide variety of machine learning applications \cite{bishop2006pattern,chang2011libsvm,friedman2001elements,LeCun2015}.
The main challenge of solving these problems stems from the fact that they are often high-dimensional and nonlinear, and may be nonconvex.


For many machine learning applications, a common approach is to employ first-order methods such as the Stochastic Gradient method (SGD). 
SGD and its variance-reduced, adaptive and distributed variants \cite{robbins1951stochastic, Nguyen2017, johnson2013accelerating,Harikandeh2015,Reddi2016,schmidt2017minimizing,duchi2011adaptive,kingma2014adam,recht2011hogwild} are popular because they are simple to implement and have low per-iteration cost.
However, these methods often require significant tuning efforts (for each problem) to ensure practical performance, and they struggle on ill-conditioned problems.

One avenue for mitigating the aforementioned issues is the use of second-order and quasi-Newton methods \cite{nocedal_book,Fletcher1987,dennis1977quasi}. These methods, in the deterministic setting, are relatively insensitive to their associated hyper-parameters and are able to alleviate the effects of ill-conditioning. Unfortunately, a drawback of these methods is that they often do not scale sufficiently well with the high dimensionality (both $n$ and $d$) typical in machine learning and deep learning problems. Thus, the computational burden of using deterministic higher-order methods is often deemed to be too high. 

Recently attention has shifted towards stochastic second-order \cite{byrd2011use,martens2010deep,bollapragada2016exact,Roosta-Khorasani2018} and quasi-Newton \cite{schraudolph2007stochastic,byrd2016stochastic,curtis2016self,mokhtari2015global,berahas2020robust,berahas2016multi} methods. These methods attempt to combine the speed
of Newton’s method and the scalability of first-order methods by incorporating curvature information
in a judicious manner, and have proven to work well for several machine learning tasks \cite{berahas2017investigation,xu2020second}. However, the question of how to balance the accuracy in the gradient and Hessian approximation is yet unresolved, and as such these methods often perform on par with their first-order variants.

Several other attempts have been made to balance the first- versus second-order trade-off, in order to find ways of incorporating partial curvature information to help with ill-conditioning at an acceptable cost. For example, variants of coordinate descent methods that perform second-order-type updates restricted to a low dimensional subspace, are prototypical methods in this niche \cite{richtarik2014iteration,richtarik2016parallel,Fountoulakis2018,Tappenden2016}. However, while progress has been made, the gap between first- and second-order methods remains.


In this paper, we propose the  \textbf{S}ymmetric bl\textbf{O}ckwise tru\textbf{N}cated optim\textbf{I}zation \textbf{A}lgorithm (\SONIA{}). \SONIA{} aims to bridge the gap between first- and second-order methods, but is different in nature to coordinate descent methods because at every iteration a step in the \emph{full dimensional space} is generated. The search direction consists of two components. The first component lies in an $m$-dimensional subspace (where $m\ll d$ is referred to as the `memory' and is user defined), and is generated using a second-order approach. The second component of the update lies in the orthogonal complement, and is an inexpensive scaled steepest descent update. The combination of the two components allows for the overall search direction to explore the full-dimensional space at every iteration.

\vspace{-5pt}
\paragraph{Contributions} Our contributions can be summarized as follows:
\begin{itemize}[noitemsep,nolistsep,topsep=0pt,leftmargin=15pt]
\item \textit{Novel Optimization Algorithm}. We propose \SONIA{} for solving empirical risk minimization problems that attempts to bridge the gap between first- and second-order methods. The algorithm judiciously incorporates curvature information in one subspace (whose dimension is determined by the user) and takes a gradient descent step in the complement of that subspace. As such, at every iteration, \SONIA{} takes a step in full dimensional space while retaining a low per-iteration cost and storage, similar to that of limited memory quasi-Newton methods.
\item \textit{Theoretical Analysis}. We derive convergence guarantees for \SONIA{}, both in deterministic and stochastic regimes, for strongly convex and nonconvex optimization problems. These guarantees match those of popular quasi-Newton methods such as L-BFGS.
\item \textit{Stochastic Variant of \SONIA{}}. We develop and analyze a stochastic variant of \SONIA{} that uses stochastic gradient and Hessian approximations in lieu of the true gradient and Hessian.
\item \textit{Competitive Numerical Results}. We investigate the empirical performance of the deterministic and stochastic variants of \SONIA{} on strongly convex (logistic regression) and nonconvex (nonlinear least squares) problems that arise in machine learning. Our proposed methods are competitive with the algorithms of choice in both the deterministic and stochastic settings.
\end{itemize}

\vspace{-5pt}
\paragraph{Organization} Related works are described in Section \ref{sec:rel_work}.  Section \ref{sec:SONIA} presents our proposed algorithm, \SONIA{}, and its stochastic variant. We show the theoretical properties of our proposed method in Section \ref{sec:theory}. Numerical results on deterministic and stochastic problems are reported in Section \ref{sec:num_res}. Finally, in Section \ref{sec:fin_rem} we provide some final remarks and discuss possible avenues for future research.


\vspace{-5pt}
\section{Related Work}\label{sec:rel_work}
\vspace{-5pt}
As in this work, the following works employ iterate updates of the form
\begin{eqnarray}\label{wupdate}
    w_{k+1} = w_k + \alpha_k p_k,
\end{eqnarray}
where $p_k \in \mathbb{R}^d$ is the search direction and $\alpha_k>0$ is the step length or learning rate. 

The work in \cite{paternain2019newton} proposes a Newton-type algorithm for nonconvex optimization problems. At each iteration the construction and  eigenvalue decomposition (full dimensional) of the Hessian is required, small (in modulus) eigenvalues are truncated, and a Newton-like search direction is generated using the truncated inverse Hessian instead of 
the true inverse Hessian. The method works well in practice and is guaranteed to converge to local minima, but is expensive.

Quasi-Newton methods--methods that compute search directions using (inverse) Hessian approximations that are constructed using past iterate and gradient information--represent some of the most effective algorithms for minimizing nonlinear objective functions. This class of nonlinear optimization algorithms includes BFGS, DFP and SR1; see \cite{Fletcher1987,nocedal_book,dennis1977quasi} and the references therein.

The Symmetric Rank One (SR1) update is a special case of a rank one quasi-Newton method \cite{dennis1977quasi}.  It is the unique symmetric rank-1 update that satisfies the secant condition $B_{k+1} s_{k} = y_{k}$, where $s_k = w_k - w_{k-1}$ and $ y_k = \nabla F(w_k)- \nabla F(w_{k-1})$ are the curvature pairs and the Hessian approximation is updated at every iteration via:  
$ 
B_{k+1} = B_k + \tfrac{(y_k-B_ks_k)(y_k-B_ks_k)^T}{s_k^T(y_k-B_ks_k)}.
$ 
Hessian approximations using SR1 updates are not guaranteed to be positive definite, and while this was originally seen as a drawback, it is arguably viewed as an advantage in the context of nonconvex optimization.
Limited memory variants exist where there is a fixed memory size $m$, and only the last $m$ curvature pairs are kept and used to construct the Hessian approximation. Let $S_k=[s_{k-m+1},\dots,s_k] \in \mathbb{R}^{d \times m}$ and $Y_k=[y_{k-m+1},\dots,y_k] \in \mathbb{R}^{d \times m}$ denote matrices consisiting of the $m$ most recent curvature pairs.  As studied in \cite{byrd1994}, the compact form of L-SR1 is as follows:
\begin{equation}\label{SR1compact}
B_k = B_0 +(Y_k - B_0S_k) (L_k+D_k+L_k^T - S_k^TB_0S_k)^{-1} (Y_k - B_0S_k)^T,
\end{equation}
where $S_k^TY_k = L_k + D_k + U_k$, $L_k$ denotes the strictly lower triangular part, $D_k$ is the diagonal and $U_k$ denotes the strictly upper triangular part of $S_k^TY_k$, respectively, and $B_0$ is an initial approximation (usually set as $B_0 = \eta I$, $\eta>0$).
A key observation is that while $B_{k+1}$ is a full dimensional $d \times d$ matrix, the inverse in \eqref{SR1compact} is a small $m\times m$ matrix. Recent works that employ the compact L-SR1 update include \cite{erway2019trust,Brust2017,berahas2019quasi}, where \eqref{SR1compact} defines the quadratic model within a trust region algorithm.

Rather than maintaining a history of the $m$ most recent curvature pairs, the work \cite{berahas2019quasi} proposes a \emph{sampled} variant of the L-SR1 update. In that work, at each iteration $k\geq 0$, $m$ directions $\{s_1,\dots,s_m\}$ are sampled around the current iterate $w_k$ and stored as $S_k=[s_1,\dots,s_m]\in \R^{d\times m}$. Next, the gradient displacement vectors are computed via 
\begin{equation}\label{curvepairs}
    Y_k = \nabla^2 F(w_k) S_k,
\end{equation}
and the matrix $B_0$ is set to zero.
In this way, previous curvature information is `forgotten', and local curvature information is `sampled' around the current iterate. Moreover, depending on the way the vectors $\{s_1,\dots,s_m\}$ are sampled, one can view \eqref{curvepairs} as a sketch of the Hessian \cite{woodruff2014sketching}.

The approach proposed here combines quasi-Newton updates for indefinite Hessians \cite{erway2019trust} with sampled curvature pairs \cite{berahas2019quasi}. Moreover, a truncation step (as in \cite{paternain2019newton}) allows us to avoid checking conditions on the curvature pairs, and ensures that the Hessian approximations constructed are positive definite. The subspace generation is based on the user-defined hyper-parameter $m$ (as such the user has full control over the computational cost of each step), and an eigenvalue decomposition step (which is performed in reduced dimension as so is cheap).

\vspace{-5pt}
\section{\textbf{S}ymmetric bl\textbf{O}ckwise tru\textbf{N}cated optim\textbf{I}zation \textbf{A}lgorithm (\SONIA{})}\label{sec:SONIA}
\vspace{-5pt}
In this section, we present our proposed algorithm. We begin by motivating and describing the deterministic variant of the method and then discuss its stochastic counterpart. We end this section by discussing the per iteration complexity of \SONIA{}.


\vspace{-5pt}
\subsection{Deterministic \SONIA{}}\label{sec:det_SONIA}
\vspace{-5pt}
The \SONIA{} algorithm generates iterates according to the update \eqref{wupdate}. The search direction $p_k$ consists of two components; the first component lies in one subspace and is a second-order based update, while the second component lies in the orthogonal complement and is a scaled steepest descent direction. We now describe how the subspaces are built at each iteration as well as how to compute the second-order component of the search direction. 

The algorithm is initialized with a user defined memory size $m\ll d$. At each iteration $k\geq 0$ of  \SONIA{}, $m$ directions $\{s_1,\dots,s_m\}$ are randomly sampled, and curvature pair matrices $S_k$ and $Y_k$ are constructed via \eqref{curvepairs}. 
Setting $B_0 = 0$, and substituting into \eqref{SR1compact} gives the compact form of the Hessian approximation used in this work\footnote{If $S_k^TY_k$ has full rank, then the pseudo-inverse in \eqref{SLSR1} is simply the inverse.}:
\begin{eqnarray}\label{SLSR1}
B_k = Y_k (Y_k^T S_k)^{\dagger} Y_k^T.
\end{eqnarray}
Similar to the strategy in \cite{erway2019trust}, using the `thin' $QR$ factorization of $Y_k =Q_kR_k$, where $Q_k \in \mathbb{R}^{d \times m}$ has orthonormal columns and $R_k\in \mathbb{R}^{m \times m}$ is an upper triangular matrix, \eqref{SLSR1} gives
\begin{equation}\label{eq:2}    
B_k = Q_kR_k (Y_k^T S_k)^{\dagger} R_k^TQ_k^T.
\end{equation}

Note that $B_k$ is symmetric because, by \eqref{curvepairs}, $ (Y_k^T S_k)^{\dagger} = (S_k^T\nabla^2 F(w_k) S_k)^{\dagger}\in \R^{m\times m}$ is symmetric. Thus, by the spectral decomposition, $R_k(Y_k^T S_k)^{\dagger}R_k^T = V_k \Lambda_k V_k^T$, where the columns of $V_k \in \R^m$ form an orthonormal basis (of eigenvectors), and $\Lambda_k \in \mathbb{R}^{m \times m}$ is a diagonal matrix containing the corresponding eigenvalues. Substituting this into \eqref{eq:2}
gives $B_k = Q_kV_k \Lambda_k V_k^TQ_k^T.$
Since $Q_k$ has orthonormal columns and $V_k$ is an orthogonal matrix, it is clear that 
\begin{equation}\label{tildeV}
    \V_k \eqdef Q_kV_k \in \mathbb{R}^{d \times m}
\end{equation}
has orthonormal columns. Finally, the low rank decomposition of the Hessian approximation $B_k$ is
\begin{equation}\label{eq:4}
B_k = \V_k \Lambda_k \V_k^T.
\end{equation}
The following definition is motivated by \cite[Definition~2.1]{paternain2019newton}.
\begin{definition}\label{def:Ak}
Let $B_k$, $\V_k$ and $\Lambda_k$ be the matrices in \eqref{eq:4}, and let $\epsilon > 0$. The truncated inverse Hessian approximation of $B_k$ is $A_k\eqdef \V_k |\Lambda_k|^{-1}_{\epsilon} \V_k^T$,
where
$(|\Lambda_k|_{\epsilon})_{ii} = 
\max\{
|\Lambda_k|_{ii}, \epsilon\}$.
\end{definition}

Definition~\ref{def:Ak} explains that any eigenvalues in $\Lambda_k$ below the threshold $\epsilon$ are truncated and set to $\epsilon$. This is useful for several reasons. Firstly, it ensures that the search direction consists only of directions with non-negligible curvature. Moreover, unlike classical quasi-Newton methods that enforce conditions on the curvature pairs to guarantee that the (inverse) Hessian approximations are well-defined and that the updates are stable, \SONIA{} utilizes a truncation step (as described in Definition~\ref{def:Ak}) and as such avoids the need for any such safe-guards. The reason for this is that even if $Y_k^T S_k$ is rank deficient, the truncation step ensures that $A_k$ has full rank. This is especially important with SR1-type methods that require matrix-vector products in the checks.

Before we proceed, we make a few more comments about our algorithmic choice of constructing the gradient differencing curvature pairs via \eqref{curvepairs}. As mentioned above, this ensures that $Y_k^T S_k$ is symmetric which is a fundamental component of our approach for three main reasons. Firstly, it ensures that the Hessian approximations constructed are symmetric. Secondly, it allows us to utilize the spectral decomposition. And, thirdly, unlike the classical SR1 method that utilizes only the lower triangular part of the $Y_k^T S_k$ matrix to construct Hessian approximations (see \eqref{SR1compact}) and as such throws away possibly useful curvature information, our approach allows us to use all curvature information collected at every iteration. Moreover, one can show that constructing curvature pairs in this manner guarantees that the secant equations hold for all curvature pairs, and that the Hessian approximations are scale invariant.

Next we discuss the subspace decomposition. The gradient is orthogonally decomposed as:
\begin{equation}\label{gradF}
\nabla F(w_k) = g_k + g_k^{\perp},
\quad \mbox{where}\quad 
g_k = \V_k\V_k^T \nabla F(w_k) \quad 
\text{and} \quad g_k^{\perp} = (I-\V_k\V_k^T) \nabla F(w_k).
\end{equation}
Clearly, $g_k \in {\rm range}(\V_k\V_k^T)$ and $\gp_k \in {\rm ker}(\V_k\V_k^T) \equiv {\rm range}(I - \V_k\V_k^T)$. Vectors $g_k$ and $\gp_k$ are orthogonal because the subspaces ${\rm range}(\V_k\V_k^T)$ and ${\rm range}(I - \V_k\V_k^T)$ are orthogonal complements (i.e., $g_k^T\gp_k =  \nabla F(w_k)^T\V_k\V_k^T(I-\V_k\V_k^T) \nabla F(w_k) = 0$).



The \SONIA{} search direction is 
%
\begin{align}\label{p}
    p_k &\eqdef -\V_k |\Lambda_k|_{\epsilon}^{-1} \V_k^T \nabla F(w_k) - \rho_k (I-\V_k\V_k^T) \nabla F(w_k)  = -\V_k |\Lambda_k|_{\epsilon}^{-1} \V_k^T g_k - \rho_k \gp_k,
\end{align}
where for all $k \geq 0$, 
\begin{equation}\label{rhok}
    \rho_k \in (0,\lambda_k^{\min}], \qquad \text{and} \qquad \lambda_k^{\min} \eqdef \min_{i} \{[|\Lambda_k|^{-1}_\epsilon]_{ii}\}.
\end{equation}
The first component of the search direction lies in the subspace ${\rm range }(\V_k\V_k^T)$, while the second component lies in the orthogonal complement.

\begin{lemma}\label{pequiv}
The search direction $p_k$ in \eqref{p} is equivalent to $p_k = -\A_k \nabla F(w_k)$, where
\begin{align}\label{eq:A_k_matrix}
\A_k \eqdef \V_k |\Lambda_k|_{\epsilon}^{-1} \V_k^T + \rho_k (I-\V_k\V_k^T).
\end{align}
\end{lemma}

The search direction $p_k$ in \eqref{p} can be interpreted as follows. If the memory is chosen as $m =0$, then $p_k$ is simply a scaled steepest descent direction (in this setting, $\rho_k$ can be any positive number). On the other hand, if $m = d$, then $p_k$ incorporates curvature information in the full dimensional space.
If $0 < m < d$, then the algorithm is a hybrid of a second-order method in ${\rm range} (\V_k\V_k^T)$ and steepest descent in the orthogonal complement ${\rm null}(\V_k\V_k^T)$. Thus, this algorithm bridges the gap between first- and second-order methods. 

\begin{remark}
The following remarks are made regarding the search direction $p_k$.
\begin{itemize}[noitemsep,nolistsep,topsep=0pt,leftmargin=15pt]
    \item The first component of the search direction vanishes only if  (i) the memory size is $m = 0$, or if (ii) $\nabla F(w) \in {\rm null}(\V\V^T)$.
    \item The second component of the search direction vanishes only if (i) the memory size is $m = d$, or if (ii) ${\rm range}(\V\V^T) \equiv \R^d$.  
\end{itemize}
\end{remark}

The \SONIA{} algorithm is presented in Algorithm~\ref{alg:fastsr1}.

\begin{algorithm}
 {\small 
\caption{\SONIA{}}
  \label{alg:fastsr1}
 {\bf Input:} $w_{0}$ (initial iterate), $m$ (memory), $\epsilon$ (truncation parameter). 
 
  \begin{algorithmic}[1]
  \For {$k=0,1,2,...$}
    \State Compute the gradient $\nabla F(w_k)$
    \State Construct a random matrix $S_k \in \mathbb{R}^{d \times m}$ and set $Y_k$ via \eqref{curvepairs} 
    \State Compute the $QR$ factorization of $Y_k (=Q_kR_k)$
    \State Compute the spectral decomposition of $R_k(Y_k^T S_k)^{\dagger}R_k^T (= V_k \Lambda_k V_k^T)$
    \State Construct $\V_k (=Q_kV_k)$ via \eqref{tildeV}
    \State Truncate the eigenvalues of $\Lambda_k$ to form $|\Lambda_k|_\epsilon$ and set $\rho_k$ via \eqref{rhok}
    \State Decompose the gradient $\nabla F(w_k) (= g_k + g_k^{\perp})$ via \eqref{gradF}
    \State Compute the search direction
    $p_k$ via \eqref{p}
    \State Select the steplength $\alpha_{k} >0$, and set $w_{k+1}=w_k+\alpha_{k}p_{k}$
\EndFor
  \end{algorithmic}
  }
\end{algorithm}

\vspace{-5pt}
\subsection{Stochastic \SONIA{}}
\vspace{-5pt}
The \SONIA{} algorithm presented in Section \ref{sec:det_SONIA}, requires a gradient evaluation and a Hessian-matrix product (to construct $Y_k$) at every iteration. For many machine learning applications $n$ and $d$ are large, and thus the required computations can be prohibitively expensive. To overcome these difficulties, we present a stochastic variant of the \SONIA{} algorithm that employs a mini-batch approach.

Stochastic \SONIA{} chooses a set $\mathcal{I}_k \subset [n]$, and the new iterate is computed as follows:
\begin{equation}\label{eq:stochUpd}
    w_{k+1}=w_k  - \alpha_{k}\mathcal{A}_k \nabla F_{\mathcal{I}_k}(w_k), \quad \text{where }\;  \nabla F_{\mathcal{I}_k}(w_k) = \tfrac{1}{|\mathcal{I}_k|}\textstyle{\sum}_{i \in \mathcal{I}_k} \nabla F_i (w_k)
\end{equation}
and $\mathcal{A}_k$ is the stochastic inverse truncated Hessian approximation.
Stochastic \SONIA{} uses stochastic Hessian-matrix products to construct $Y_k$, i.e., 
$
    Y_k = \nabla^2 F_{\mathcal{J}_k}(w_k)S_k,
$ 
where $\mathcal{J}_k \subset [n]$. It is important to note that for the theory (see Section \ref{sec:theory}) the sample sets $\mathcal{I}_k$ and $\mathcal{J}_k$ need to be chosen independently.

\vspace{-5pt}
\subsection{Discussion about Complexity of \SONIA{}}
\vspace{-5pt}
\begin{wraptable}{r}{7cm}
\vspace{-0.9cm}
\caption{Summary of Computational Cost and Storage (per iteration).}
\vspace{-0.1cm}
\label{tbl:cost_SONIA}
\centering
\begin{small}
\begin{tabular}{lcc}
\toprule
\textbf{method} & \textbf{computational cost} & \textbf{storage}  \\  \midrule
\textbf{NCN \cite{paternain2019newton}} & $\mathcal{O}(nd^2 + d^3)$ & $\mathcal{O}(d^2)$   \\ \hdashline
\textbf{LBFGS \cite{liu1989limited}} & $\mathcal{O}(nd)$ & $\mathcal{O}(d)$   \\ \hdashline
\textbf{LSR1 \cite{lu1996study}} & $\mathcal{O}(nd)$ & $\mathcal{O}(d)$   \\ \hdashline
\textbf{\SONIA{}} & $\mathcal{O}(nd)$ & $\mathcal{O}(d)$  \\
 \midrule
\end{tabular}
\end{small}
\vspace{-0.5cm}
\end{wraptable}
The per iteration complexity of \SONIA{} consists of: $(1)$ a Hessian-matrix product ($\mathcal{O}(mnd)$); $(2)$ a $QR$ factorization of an $d \times m$ matrix ($\mathcal{O}(dm^2)$); $(3)$ a pseudo-inverse of an $m \times m$ matrix ($\mathcal{O}(m^3)$); and, $(4)$ a spectral decomposition of an $m \times m$ matrix ($\mathcal{O}(m^3)$). The computational cost and storage requirement for the \SONIA{} algorithm are presented in Table \ref{tbl:cost_SONIA}\footnote{Note, these computations are on top of the function/gradient evaluations that are common to all methods.}, where we compare the cost and storage to popular limited-memory quasi-Newton methods and the NCN method \cite{paternain2019newton}. Note that the \SONIA{} algorithm was developed for the regime where $m \ll d,n$. As is clear from Table \ref{tbl:cost_SONIA}, \SONIA{} has similar cost and storage to LBFGS and LSR1, and is significantly more efficient, in both regards, to the NCN method. We should note that the computational cost and storage requirements for stochastic \SONIA{} are $\mathcal{O}(d)$.

\vspace{-5pt}
\section{Theoretical Analysis}\label{sec:theory}
\vspace{-5pt}
Here, theoretical results for \SONIA{} are presented, in the deterministic and stochastic settings, for both strongly convex and nonconvex objective functions. Before we present the main theorems, we state two preliminary Lemmas that are used throughout this section. Proofs can be found in Appendix~\ref{sec:app_theory}.


\begin{assumption} \label{assum:diff} The function $F$ is twice continuously differentiable.
\end{assumption}

\begin{lemma}\label{lem.pos_def}
The matrix $\A_k$ in \eqref{eq:A_k_matrix} is positive definite for all $k\geq 0$.
\end{lemma}

\begin{lemma} \label{lem:1} If Assumption \ref{assum:diff} holds, there exist constants $0 < \mu_1 \leq \mu_2$ such that the inverse truncated Hessian approximations $\{ \mathcal{A}_k\}$ generated by Algorithm \ref{alg:fastsr1} satisfy,
\begin{align}   \label{eq:bnd_Hess1}
    \mu_1 I \preceq \mathcal{A}_k \preceq \mu_2 I,\qquad \text{for all } k=0,1,\dots.
\end{align}
\end{lemma}

\vspace{-5pt}
\subsection{Deterministic Setting}
\vspace{-5pt}
\paragraph{Strongly Convex Functions} The following assumption is standard for strongly convex functions.

\begin{assumption} \label{assum:strong_conv} There exist positive constants $\mu$ and $L$ such that $\mu I \preceq \nabla^2F(w) \preceq L I,$ $\forall w \in \mathbb{R}^d.$
\end{assumption}

\begin{theorem} \label{thm:detSC}
Suppose that Assumptions \ref{assum:diff} and \ref{assum:strong_conv} hold, and let $F^{\star} = F(w^{\star})$, where $w^{\star}$ is the minimizer of $F$. Let $\{w_k\}$ be the iterates generated by Algorithm~\ref{alg:fastsr1}, where $0 <  \alpha_k = \alpha  \leq \frac{\mu_1}{\mu_2^2 L}$, and $w_0$ is the starting point. Then, for all $k\geq 0$, 
$
 F(w_k) - F^{\star}  \leq  ( 1-\alpha \mu \mu_1  )^k  [ F(w_0) - F^{\star}  ].
$ 
\end{theorem}
Theorem \ref{thm:detSC} shows that \SONIA{} converges at a linear rate to the optimal solution of \eqref{eq:prob}. The step length range prescribed by \SONIA{} depends on $\mu_1$ and $\mu_2$, as does the rate. This is typical for limited memory quasi-Newton methods \cite{liu1989limited,berahas2019quasi}. In the worst-case, the matrix $\A_k$ can make the limit in Theorem \ref{thm:detSC} significantly worse than that of the first-order variant if the update has been unfortunate and generates ill-conditioned matrices. However, this is rarely observed in practice. 
\vspace{-5pt}
\paragraph{Nonconvex Functions} The following assumptions are needed for the nonconvex case.

\begin{assumption} \label{assum:boundedf} The function $ F$ is bounded below by a scalar $\widehat F$.
\end{assumption}

\begin{assumption} \label{assum:lip} The gradients of $F$ are $L$-Lipschitz continuous for all $w \in \mathbb{R}^d$.
\end{assumption}

\begin{theorem} \label{thm_non_MB}
Suppose that Assumptions \ref{assum:diff},  \ref{assum:boundedf} and \ref{assum:lip} hold. Let $\{w_k\}$ be the iterates generated by Algorithm~\ref{alg:fastsr1}, where
$0 <  \alpha_k = \alpha  \leq \frac{\mu_1}{\mu_2^2  L}$,
and $w_0$ is the starting point. Then, for any $T > 1$,
$\tfrac{1}{T}\textstyle{\sum}_{k=0}^{T-1} \| \nabla F(w_k) \|^2   \leq \tfrac{2[ F(w_0) - \hat{F}]}{\alpha \mu_1 T } \xrightarrow[]{T \rightarrow \infty}0.$
\end{theorem}
Theorem~\ref{thm_non_MB} bounds the average norm squared of the gradient of $F$, and shows that the iterates spend increasingly more time in regions where the objective function has small gradient. From this result, one can show that the iterates, in the limit, converge to a stationary point of $F$.


\vspace{-5pt}
\subsection{Stochastic Setting}
\vspace{-5pt}
Here, we present theoretical convergence results for the stochastic variant of \SONIA{}. 
Note that, in this section $\mathbb{E}_{\mathcal{I}_k}[\cdot]$ denotes the conditional expectation given $w_k$, whereas $\mathbb{E}[\cdot]$ denotes the total expectation over the full history. In this setting me make the following standard assumptions.

\begin{assumption} \label{assum:boundedStochGrad} There exist a  constant $\gamma$ such that $\mathbb{E}_{\mathcal{I}}[\|\nabla F_{\mathcal{I}}(w) - \nabla F(w)\|^2] \leq \gamma^2$.
\end{assumption}
\begin{assumption} \label{assum:unbiasedEstimGrad} $\nabla F_{\mathcal{I}}(w)$ is an unbiased estimator of the gradient, i.e., $\mathbb{E}_{\mathcal{I}}[\nabla F_{\mathcal{I}}(w)] = \nabla F(w)$, where the samples $\mathcal{I}$ are drawn independently.
\end{assumption}


\vspace{-5pt}
\paragraph{Strongly Convex Functions}
\vspace{-5pt}


\begin{theorem} \label{thm:stochSC}
Suppose that Assumptions \ref{assum:diff}, \ref{assum:strong_conv}, \ref{assum:boundedStochGrad} and \ref{assum:unbiasedEstimGrad} hold, and let $F^{\star} = F(w^{\star})$, where $w^{\star}$ is the minimizer of $F$. Let $\{w_k\}$ be the iterates generated by Algorithm~\ref{alg:fastsr1}, where $0 <  \alpha_k = \alpha  \leq \frac{\mu_1}{\mu_2^2 L}$, and $w_0$ is the starting point. Then,  for all $ k\geq 0$, 

\vskip-5pt 
$\mathbb{E}[F(w_{k}) - F^\star]  \leq   ( 1 - \alpha \mu_1\lambda  )^k  ( F(w_0) - F^\star - \tfrac{\alpha \mu_2^2 \gamma^2 L}{2\mu_1 \lambda} ) +  \tfrac{\alpha \mu_2^2 \gamma^2 L}{2\mu_1 \lambda}.$
\end{theorem}

The bound in Theorem \ref{thm:stochSC} has two components: $(1)$ a term decaying linearly to zero, and $(2)$ a term identifying the neighborhood of convergence. Notice that a larger step length yields a more favorable constant in the linearly decaying term, at the cost of an increase in the size of the neighborhood of convergence. As in the deterministic case, the step length range prescribed by \SONIA{} depends on $\mu_1$ and $\mu_2$, as does the rate. Thus, this result is weaker than that of its first-order variant if the update has been unfortunate and generates ill-conditioned matrices. This is seldom observed in practice. 


One can establish convergence of \SONIA{} to the optimal solution $w^\star$ by employing a sequence of step lengths that converge to zero (see \cite{robbins1951stochastic}), but at the slower, sub-linear rate. Another way to achieve exact convergence is to employ variance reduced gradient approximations \cite{johnson2013accelerating,schmidt2017minimizing}, and achieve linear convergence, at the cost of computing the full gradient every so often, or increased storage.

\vspace{-2mm}
\paragraph{Non-convex Functions}
\begin{theorem} \label{thm:stoch_non}
Suppose that Assumptions \ref{assum:diff}, \ref{assum:boundedf}, \ref{assum:lip}, \ref{assum:boundedStochGrad} and \ref{assum:unbiasedEstimGrad} hold. Let $\{w_k\}$ be the iterates generated by Algorithm~\ref{alg:fastsr1}, where $0 <  \alpha_k = \alpha  \leq \frac{\mu_1}{\mu_2^2 L}$, and $w_0$ is the starting point.

\vskip-5pt 
Then, for all $k\geq 0$,
  $\mathbb{E} [\tfrac{1}{T}\textstyle{\sum}_{k=0}^{T-1} \| \nabla F(w_k) \|^2  ]   \leq \tfrac{2[ F(w_0) - \widehat{F}]}{\alpha \mu_1 T} + \tfrac{\alpha \mu_2^2 \gamma^2 L }{\mu_1} \xrightarrow[]{T \rightarrow \infty} \tfrac{\alpha \mu_2^2 \gamma^2 L }{\mu_1}.$
\end{theorem}
Theorem \ref{thm:stoch_non} bounds the average norm squared of the gradient of $F$, in expectation, and shows that, in expectation, the iterates spend increasingly more time in regions where the objective function has small gradient. The difference with the deterministic setting is that one cannot show convergence to a stationary point; this is due to the variance in the gradient approximation employed. One can establish such convergence under an appropriate step length schedule (diminishing step lengths).

\section{Numerical Experiments}\label{sec:num_res}
\vspace{-5pt}
In this section, we present numerical experiments on several standard machine learning problems, and compare the empirical performance of \SONIA{} with that of state-of-the-art first- and second-order methods\footnote{See Section \ref{sec:tbl_algs} for details about all algorithms considered in this section.}, in both the stochastic and deterministic settings\footnote{All the codes to reproduce the experimental results will be released upon publication.}. We considered 4 different classes of problems: $(1)$ deterministic and stochastic logistic regression (stronlgy convex); and, $(2)$ deterministic and stochastic nonlinear least squares (nonconvex), and report results on 2 standard machine learning datasets\footnote{\texttt{a1a} and \texttt{gistte}. Available at: \url{https://www.csie.ntu.edu.tw/~cjlin/libsvmtools/datasets/}.}. For brevity we report only a subset of the results here and defer the rest to Appendix \ref{apndx:add}.

We compared the performance of \SONIA{} to algorithms with computational cost and storage requirements linear in both $n$ and $d$. As such, we did not compare against NCN \cite{paternain2019newton} and full-memory quasi-Newton methods. Our metric for comparison was the number of effective passes (or epochs), which we calculated as the number of function, gradient and Hessian-vector (or matrix) evaluations; see Appendix \ref{sec:Hess_mat} for more details. We tuned the hyper-parameters of each method individually for every instance; see Appendix \ref{sec:impDetndTuning} for a complete description of the tuning efforts. Where applicable, the regularization parameter was chosen from the set $\lambda \in \{10^{-3},10^{-4},10^{-5},10^{-6}\}$. The memory size was set to $\min\{d,64 \}$. Finally, the truncation parameter was set to $\epsilon =  10^{-5}$ and $\rho_k = \max_{i} \{[|\Lambda_k|^{-1}_\epsilon]_{ii}\}$; we found that these choices gave the best performance. We also performed sensitivity analysis for \SONIA; see Appendices \ref{sec:sensMMR} and \ref{sec:sensEPS}.
\vspace{-8pt}
\subsection{Deterministic Setting}
\vspace{-5pt}
In the deterministic setting, we compared the performance of \SONIA{} to that of Gradient Descent, L-BFGS \cite{liu1989limited}, L-SR1 \cite{lu1996study}, NEST+ \cite{Nesterov04} and Newton CG \cite{nocedal_book}. We implemented the algorithms  with adaptive procedures for selecting the steplength (e.g., Armijo backtracking procedure \cite{nocedal_book}) and/or computing the step (e.g., trust-region subroutine \cite{nocedal_book}). Note, Newton CG was implemented with a line search for strongly convex problems and with a trust region for nonconvex problems.\vspace{-5pt}

\begin{figure}[h!]
\centering
{\includegraphics[trim=10 110 10 110,clip, width=0.8\textwidth]{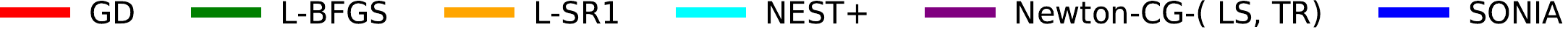}
}
\begin{subfigure}{.46\textwidth}
  	\centering
  \includegraphics[width=0.48\textwidth]{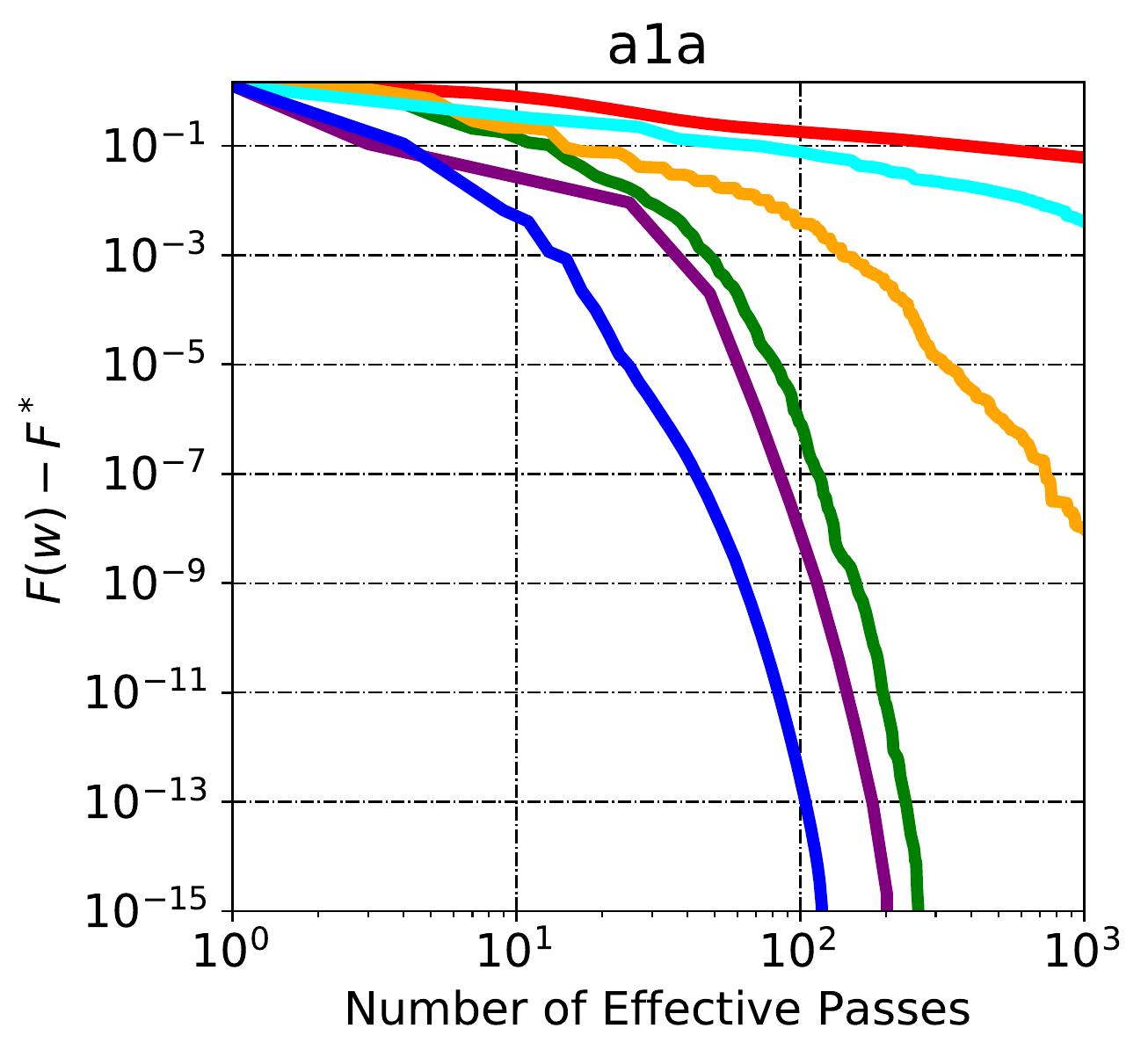}
	\includegraphics[width=0.48\textwidth]{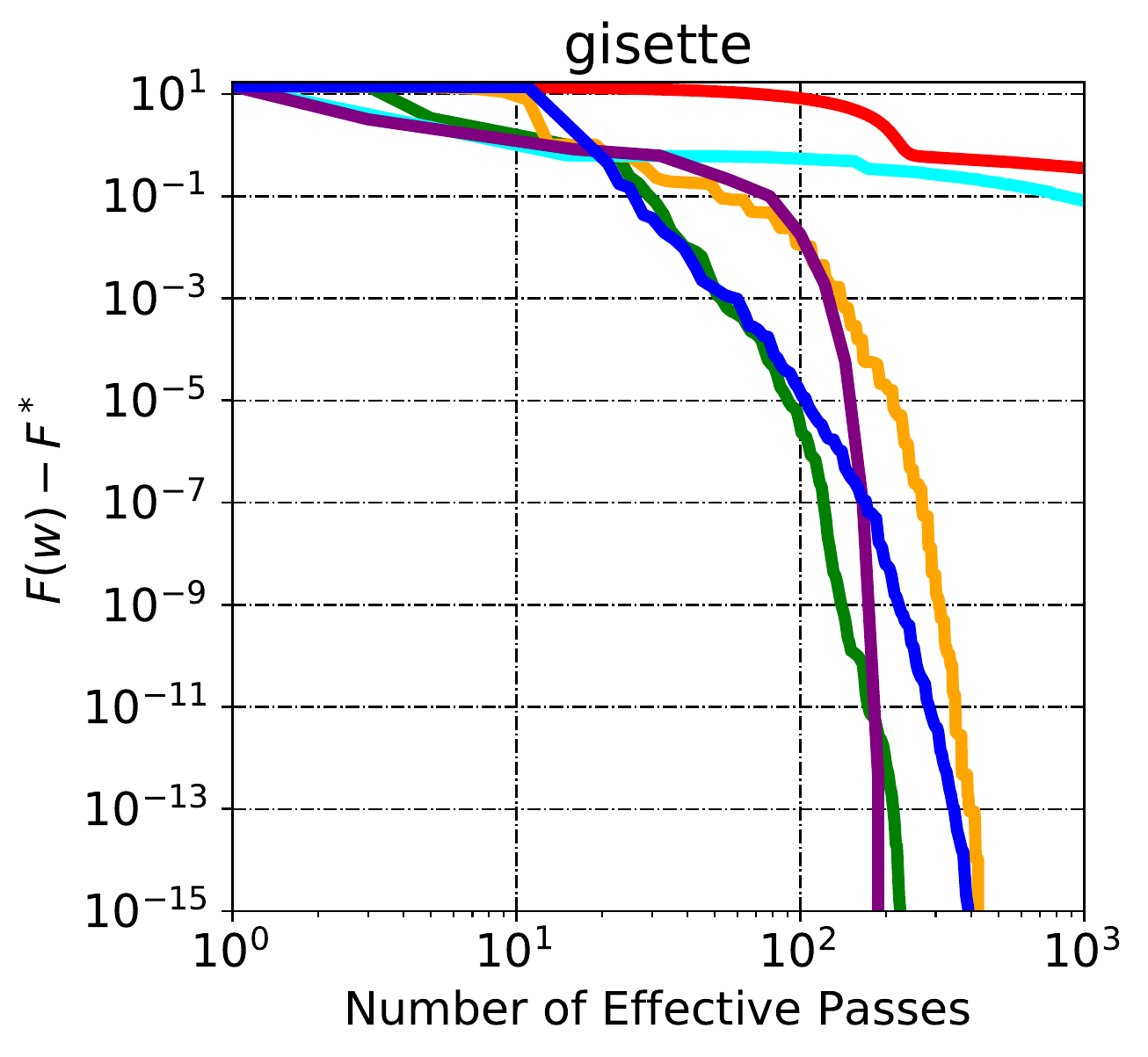}
	
	\includegraphics[width=0.48\textwidth]{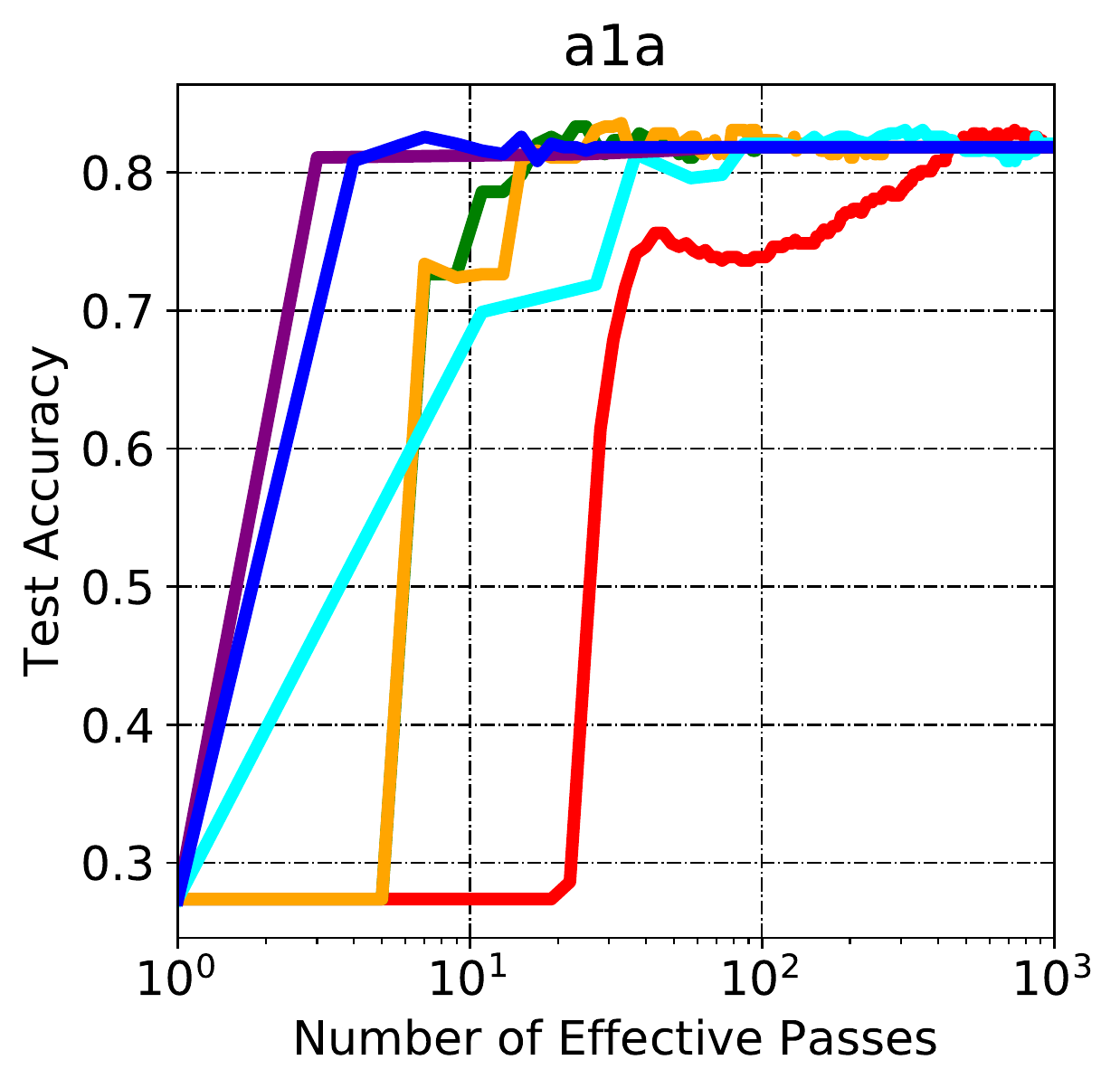}
	\includegraphics[width=0.48\textwidth]{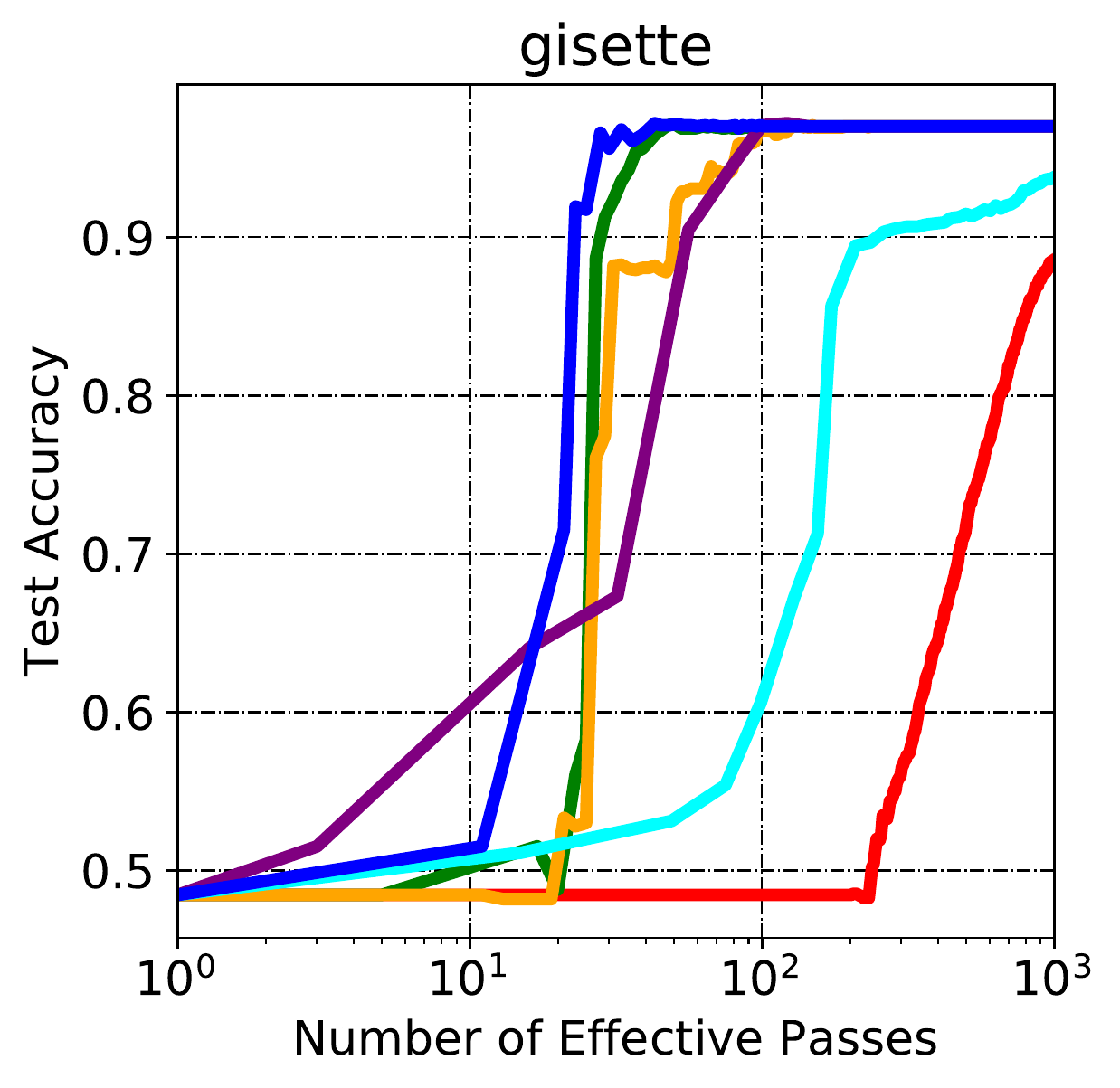}
  \caption{Comparison of optimality gap ($F(w) - F^\star$) and Test  Accuracy for different algorithms on Logistic Regression Problems.}
  \label{fig:StConv_LIBSVM}
\end{subfigure}
\hspace{0.2cm}
\begin{subfigure}{.46\textwidth}
	\centering
  \includegraphics[width=0.48\textwidth]{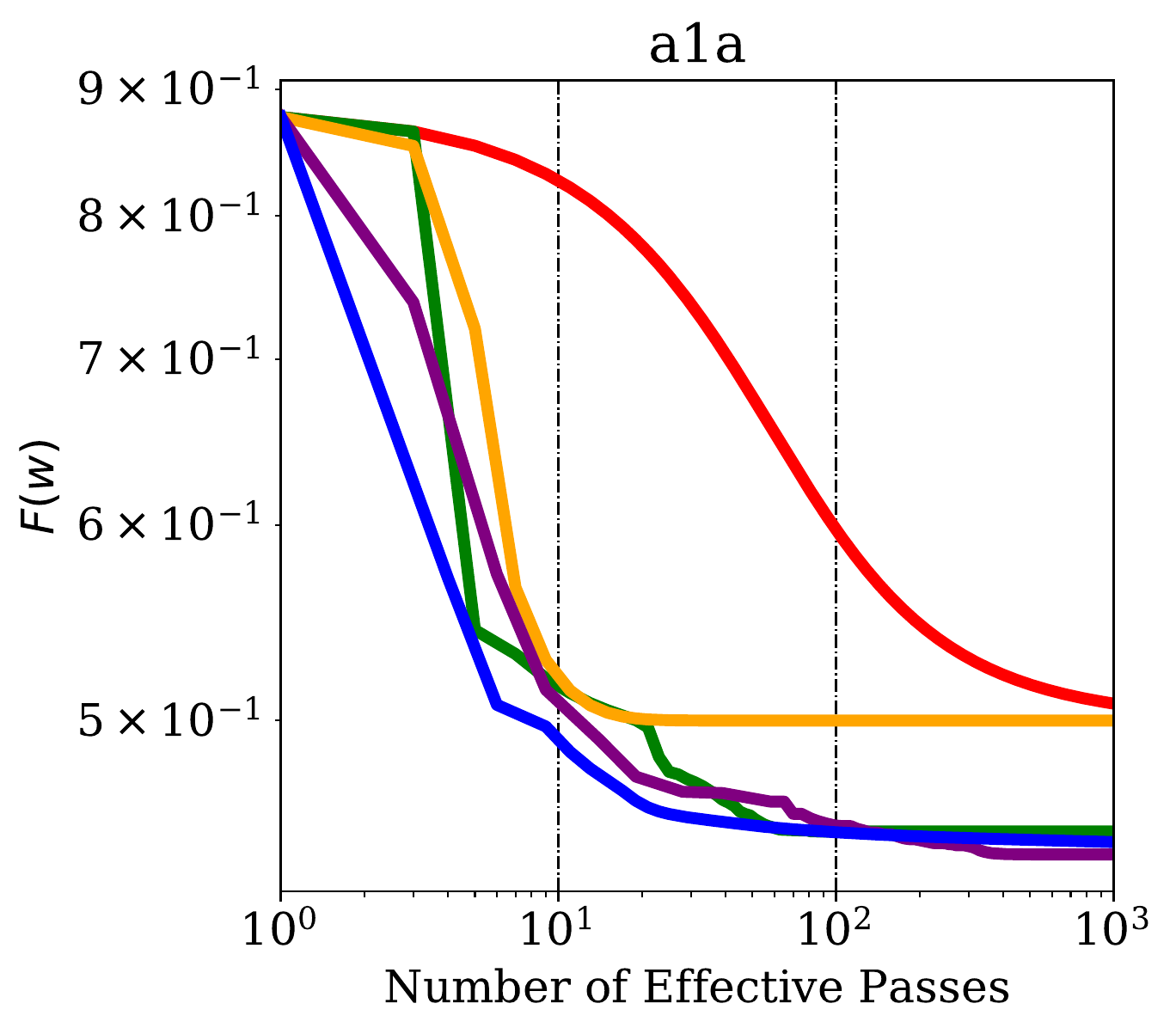}
	\includegraphics[width=0.48\textwidth]{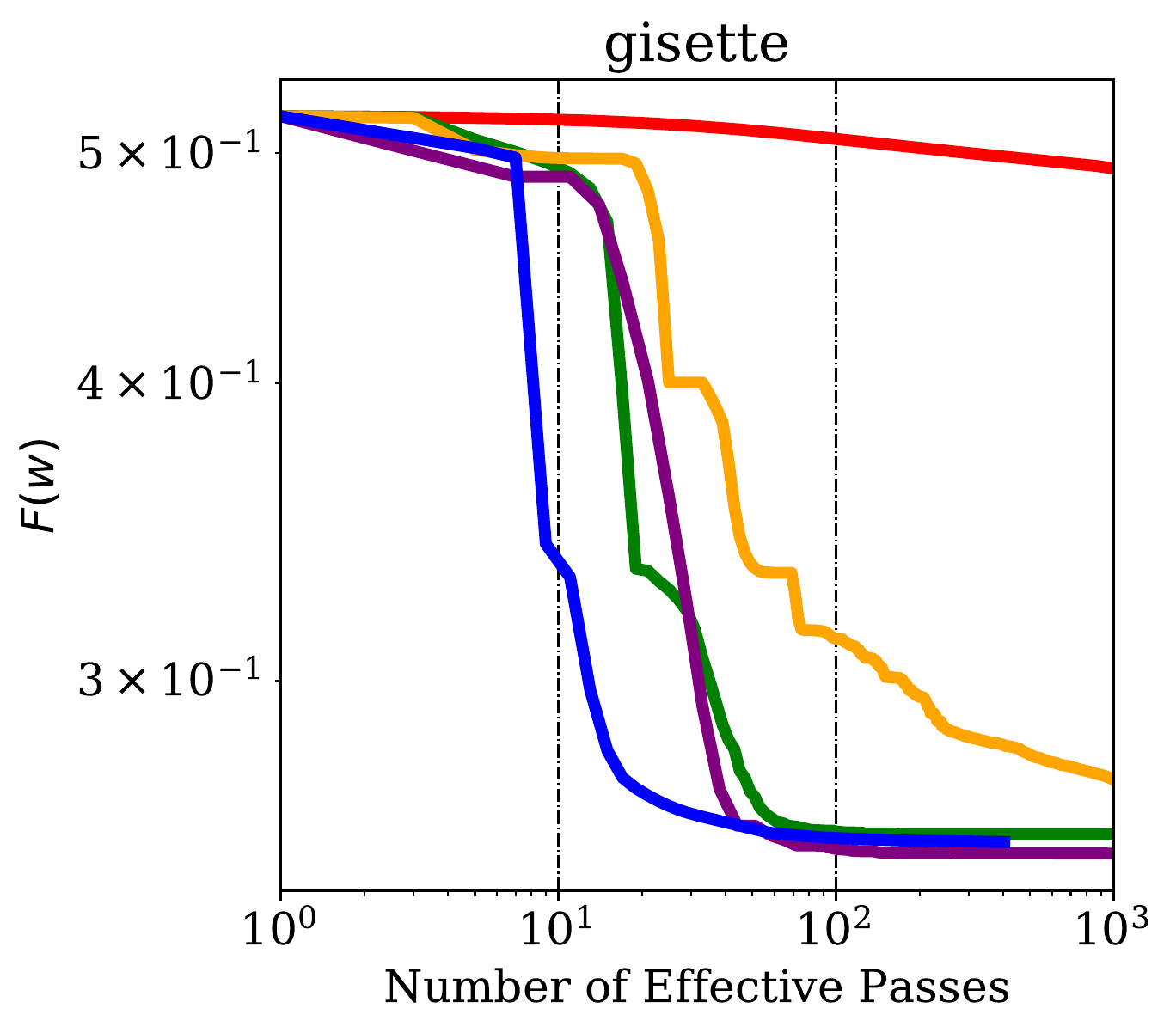}
	
	\includegraphics[width=0.48\textwidth]{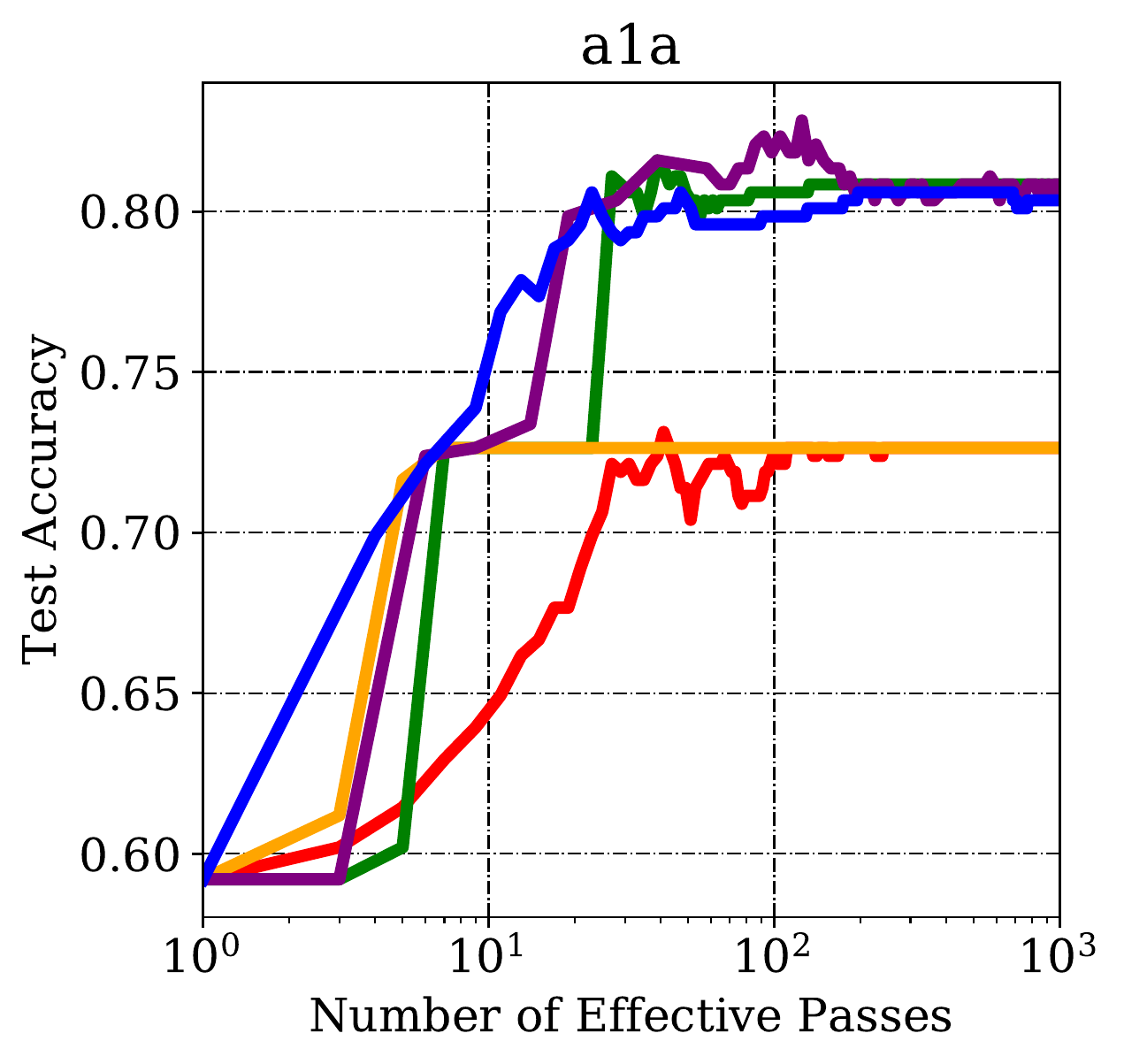}
	\includegraphics[width=0.48\textwidth]{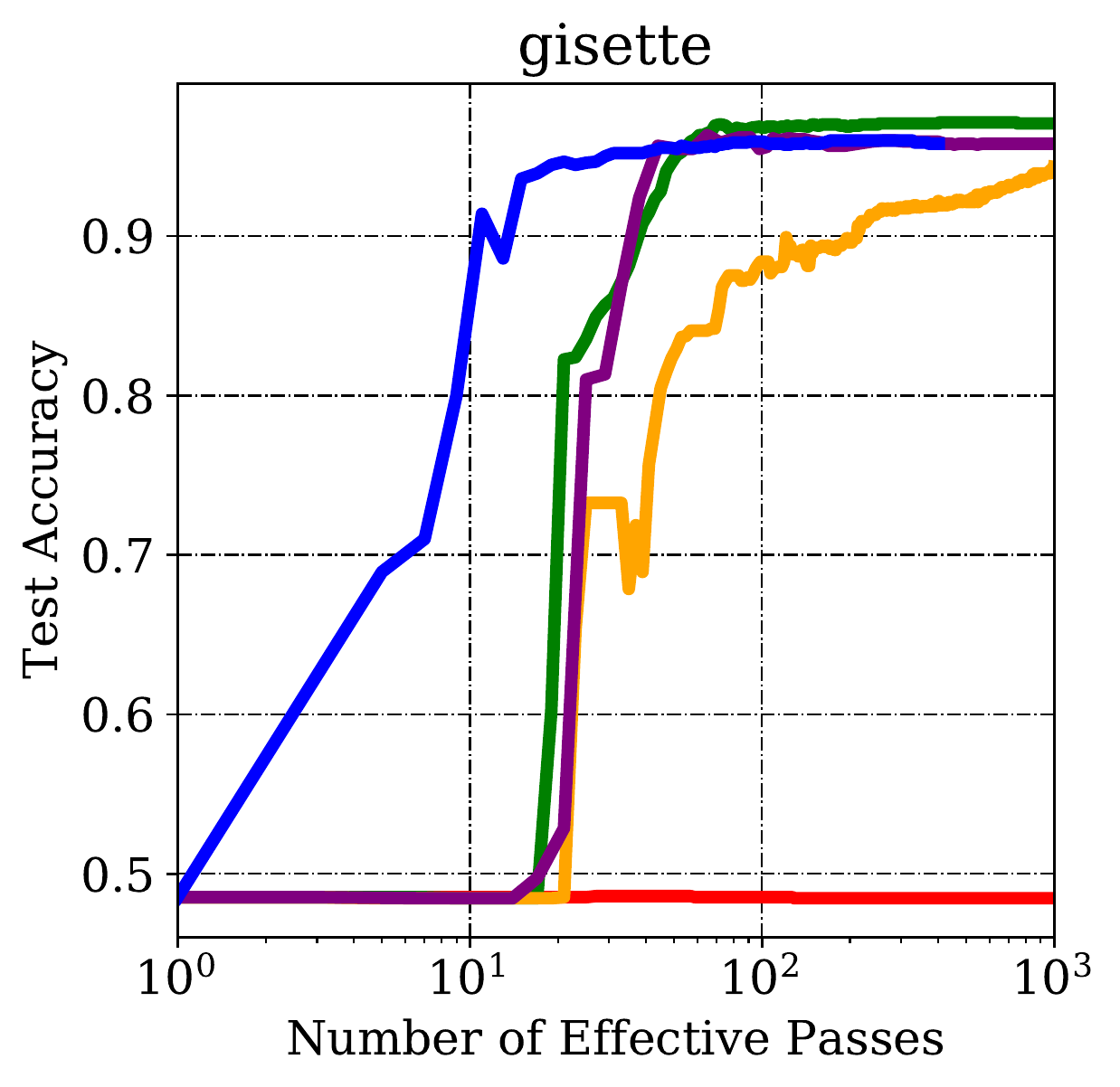}
  \caption{Comparison of objective function ($F(w)$) and Test  Accuracy for different algorithms on Non-Linear Least Squares Problems.}
  \label{fig:DetNonConv_LIBSVM}
\end{subfigure}\vspace{-4pt}
\caption{Deterministic Problems. Datasets: \texttt{a1a}, and \texttt{gistte}.}
\label{fig:fig1}
\end{figure}
\vspace{-10pt}
\paragraph{Deterministic Logistic Regression}
We considered $\ell_2$ regularized logistic regression problems, $F(w) = \tfrac{1}{n}\textstyle{\sum}_{i=1}^n \log  (1+e^{-y_ix_i^Tw} ) + \tfrac{\lambda}{2}\|w\|^2$. Figure \ref{fig:StConv_LIBSVM} shows the performance of the methods in terms of optimality gap ($F(w) - F^\star$\footnote{To find $w^\star$ we ran the ASUESA algorithm \cite{ma2017underestimate}; see Section \ref{sec:tbl_algs}.}) and testing accuracy versus number of effective passes. As is clear, the performance of \SONIA{} is on par or better than that of the other methods. Similar behavior was observed on other datasets; see Appendix \ref{sec:moreStConvResults}.
\vspace{-5pt}
\paragraph{Deterministic Non-linear Least Square}
We considered non-linear least squares problems, $ F(w) = \tfrac{1}{n}\textstyle{\sum}_{i=1}^n   (y_i - \tfrac{1}{1+e^{-x_i^Tw}}  )^2$, described in \cite{xu2020second}. Figure \ref{fig:DetNonConv_LIBSVM} shows the performance of the methods in terms of objective function and testing accuracy versus number of effective passes. As is clear, the performance of \SONIA{} is always better than the other methods in the initial stages of training, and the final objective and testing accuracy is comparable to the best method for each problem. We should note that in Figure \ref{fig:DetNonConv_LIBSVM} we report results for a single starting point (as is done in \cite{xu2020second}), but report that the performance of \SONIA{} was stable with respect to the starting point.

\vspace{-5pt}
\subsection{Stochastic Setting}
\vspace{-5pt}
In the stochastic setting, we compared the performance of \SONIA{} to that of SGD \cite{bottou2018optimization}, SARAH \cite{Nguyen2017} and SQN \cite{byrd2016stochastic}. We implemented the algorithms with fixed steplength rules, and tuned this parameter as well as the batch size for every problem; see Section \ref{sec:impDetndTuning} for more details.

\begin{figure}[ht]
\centering
{	\includegraphics[trim=10 100 10 110,clip, width=0.95\textwidth]{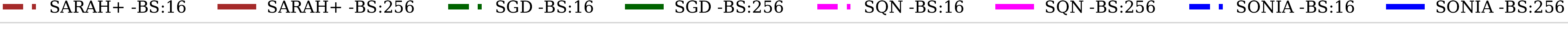}
}
\begin{subfigure}{.46\textwidth}
\centering
    \includegraphics[width=0.48\textwidth]{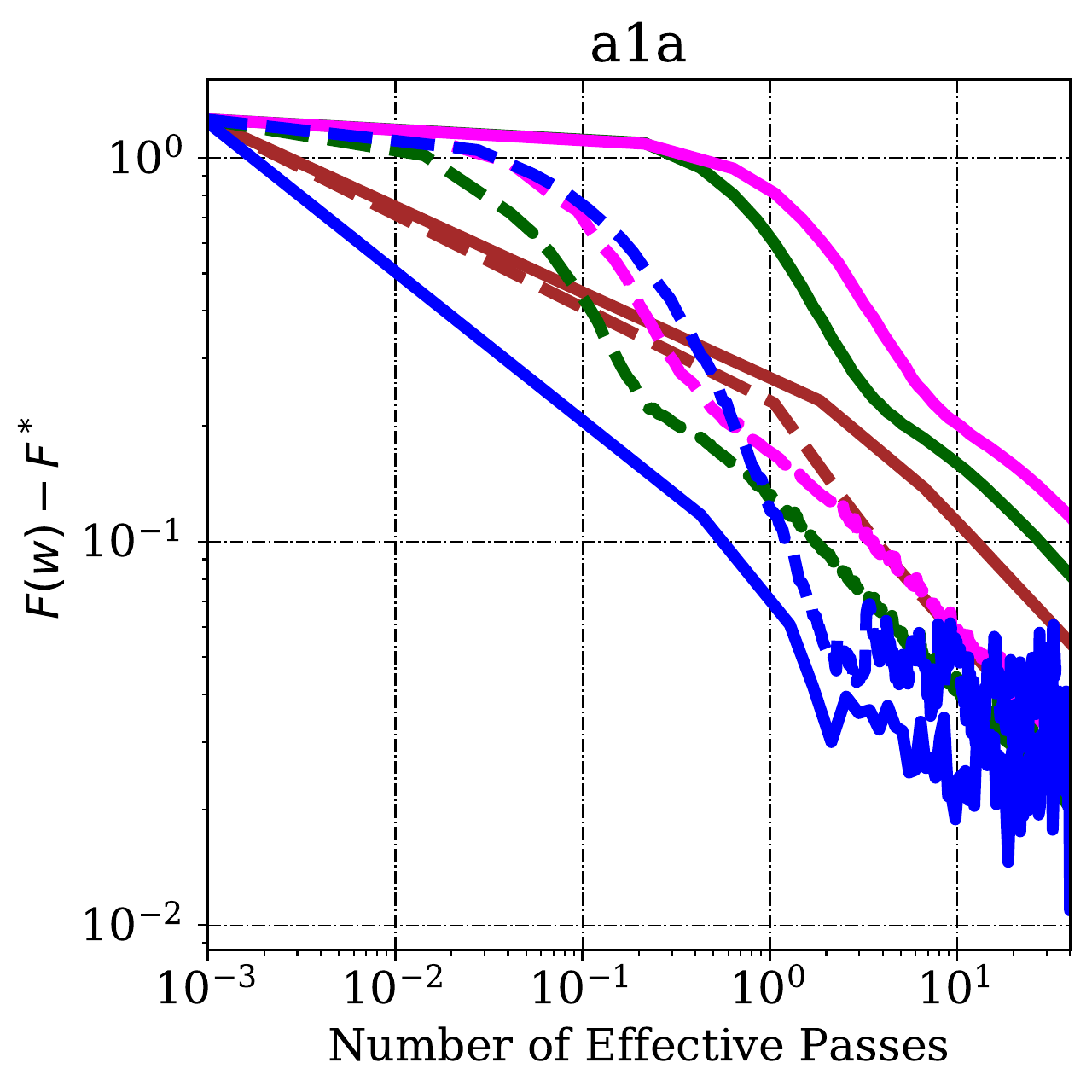}
	\includegraphics[width=0.48\textwidth]{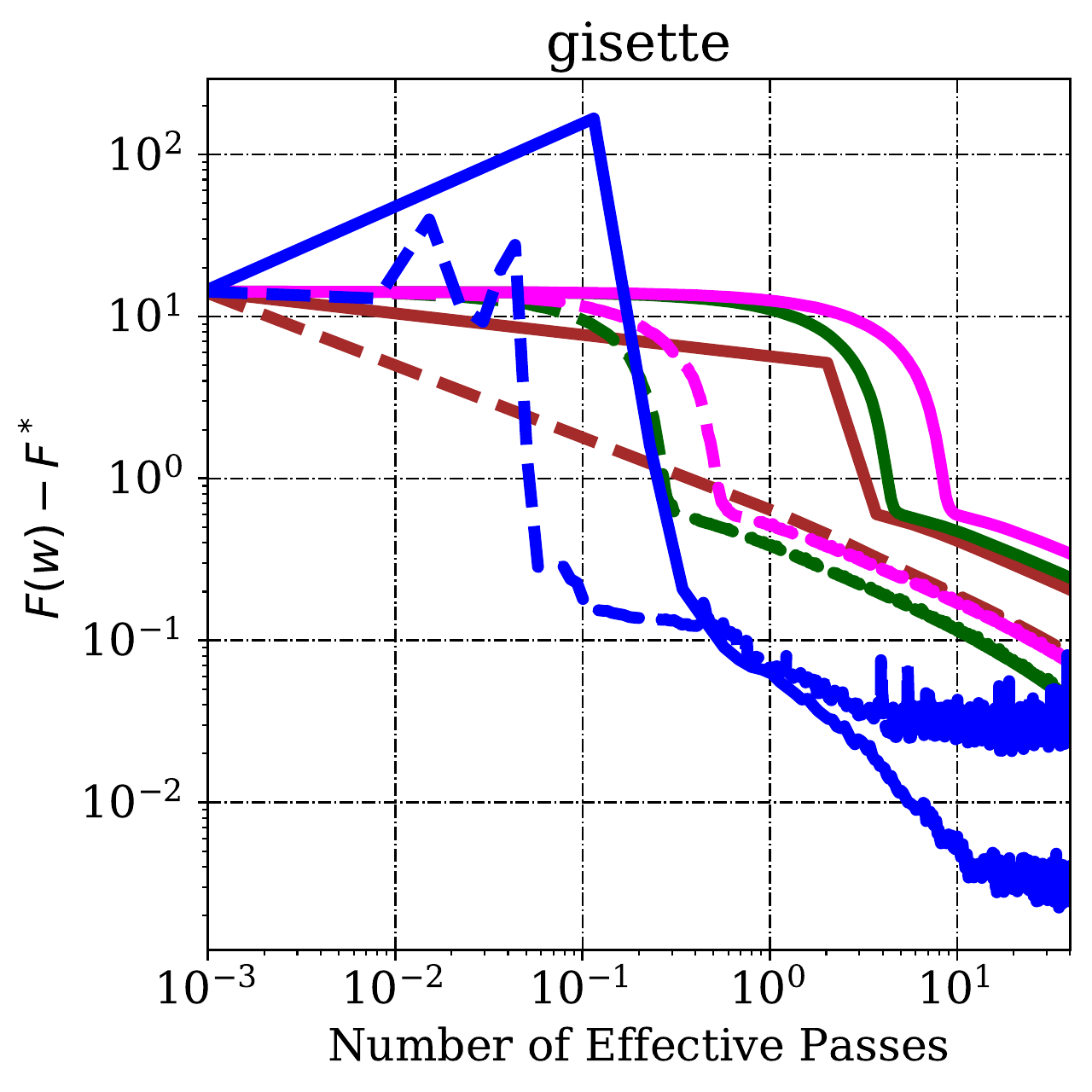}
	
	\includegraphics[width=0.48\textwidth]{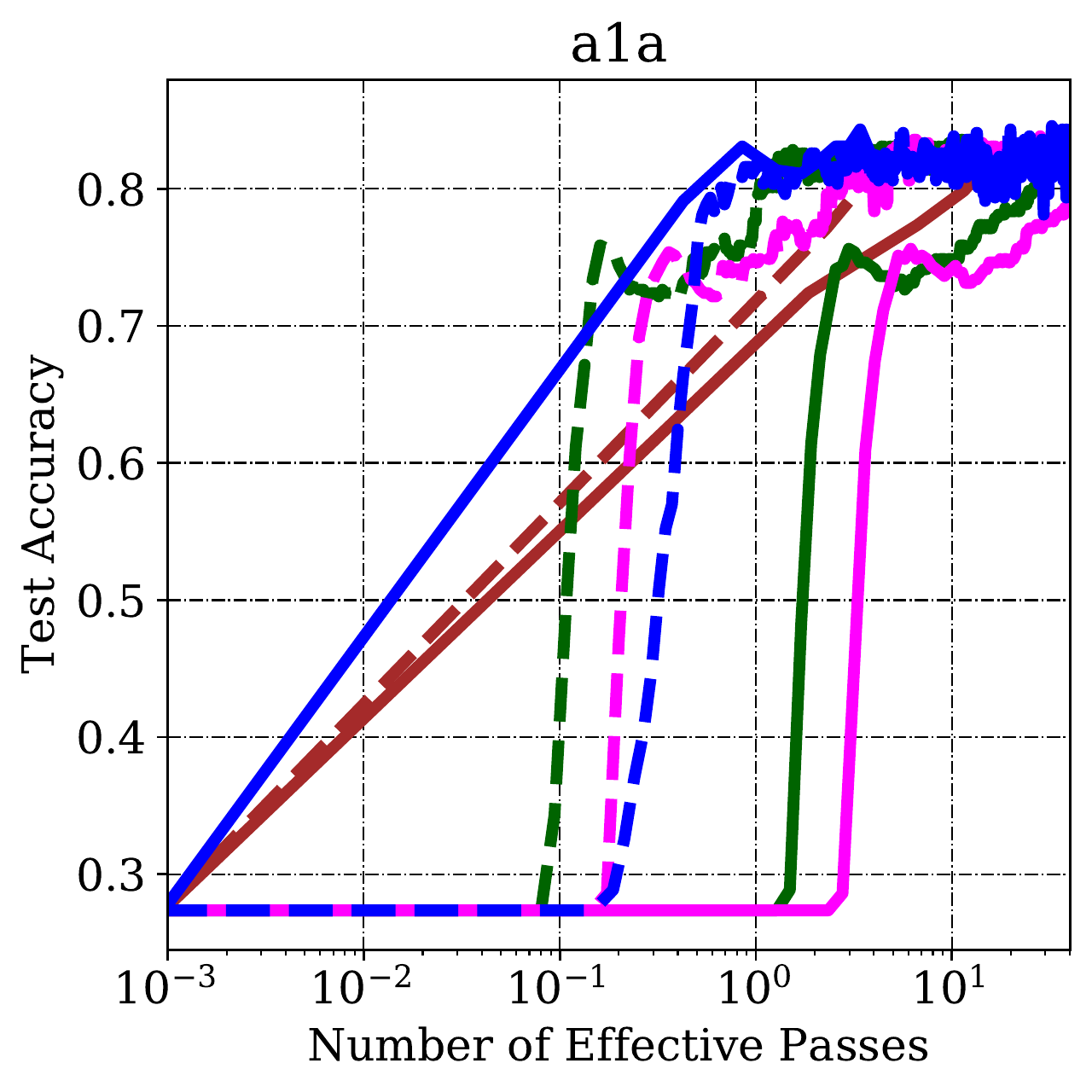}
	\includegraphics[width=0.48\textwidth]{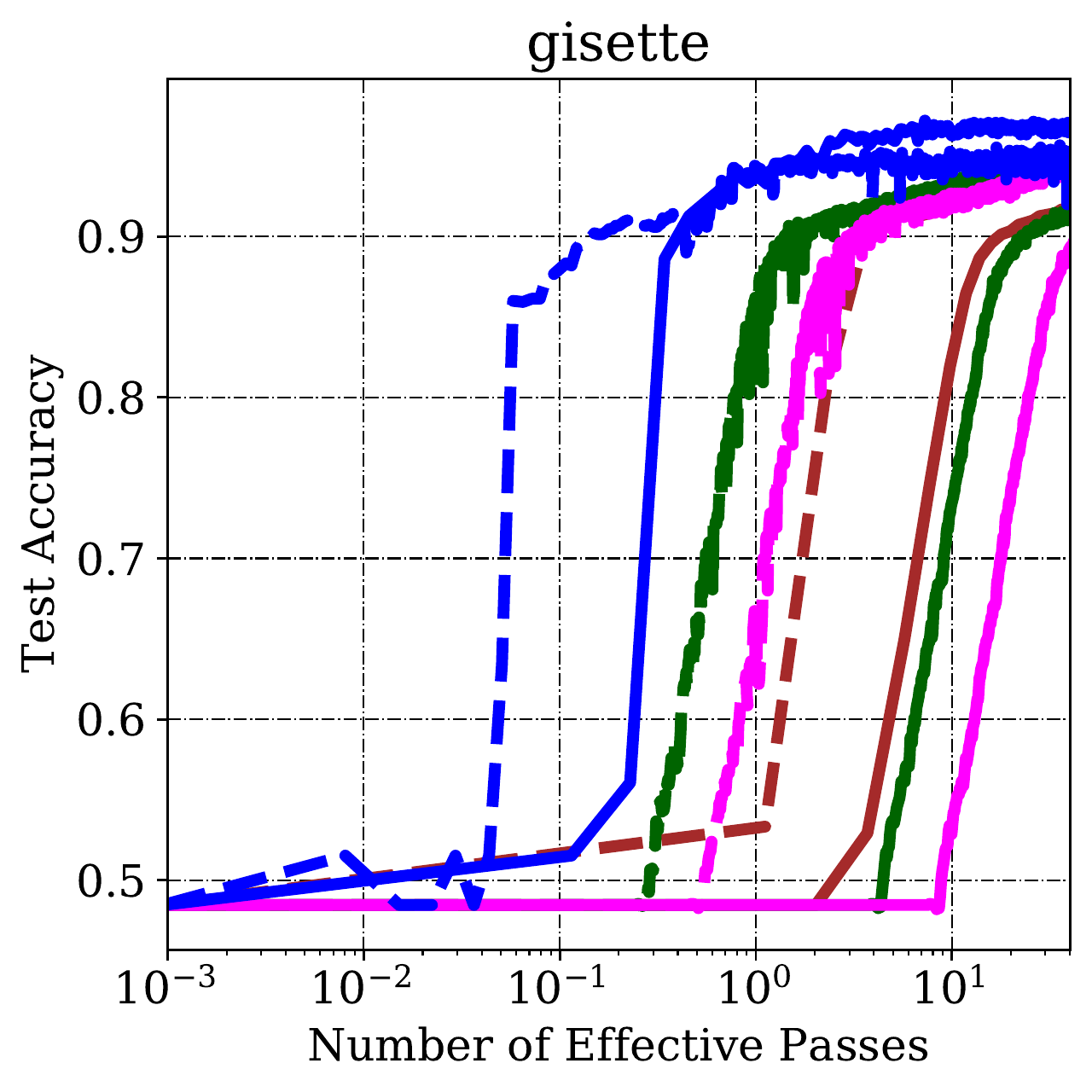}
  \caption{Comparison of optimality gap ($F(w) - F^\star$) and Test  Accuracy for different algorithms on Logistic Regression Problems.}
  \label{fig:StochSConv_LIBSVM}
\end{subfigure}
\hspace{0.2cm}
\begin{subfigure}{.46\textwidth}
	\centering
    \includegraphics[width=0.48\textwidth]{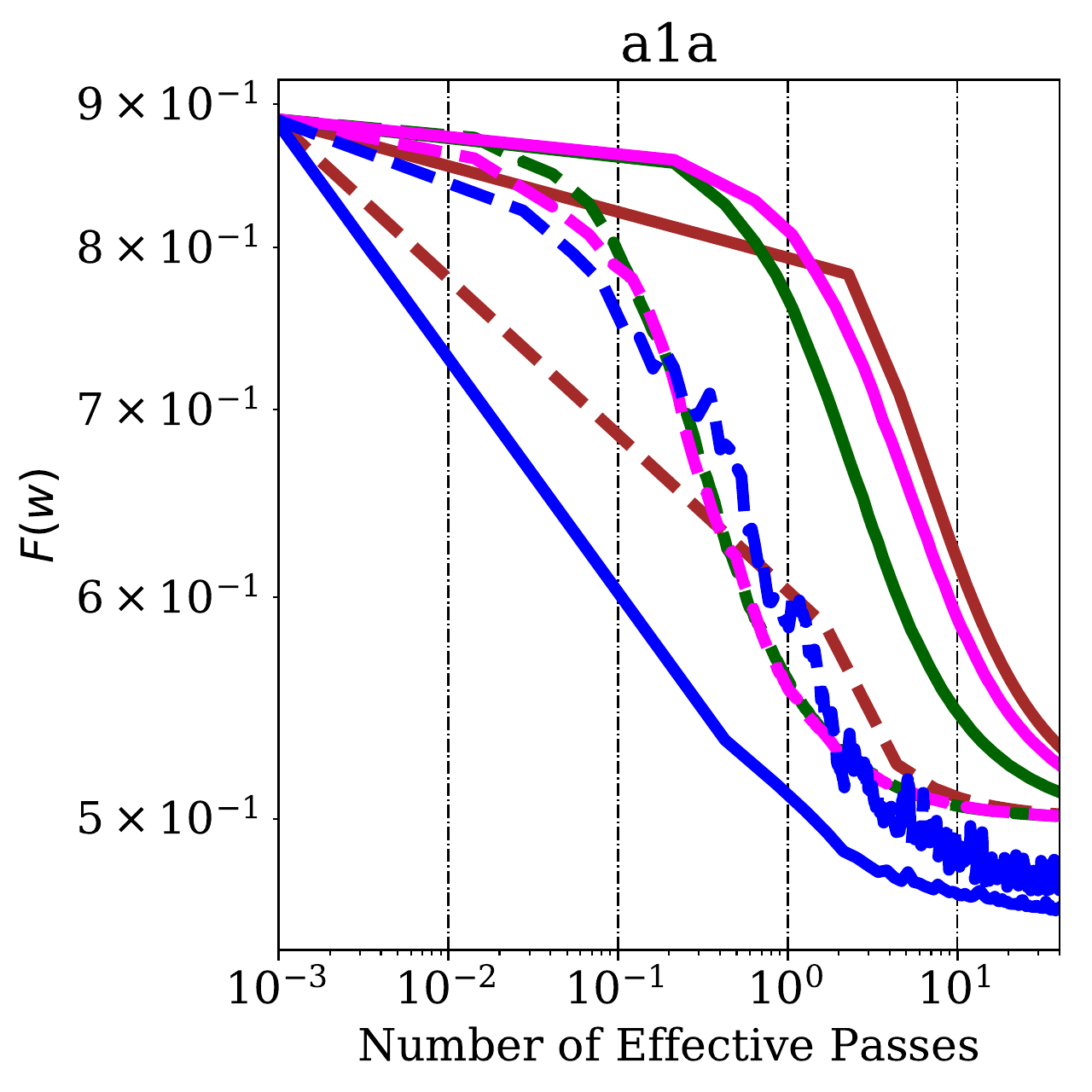}
	\includegraphics[width=0.48\textwidth]{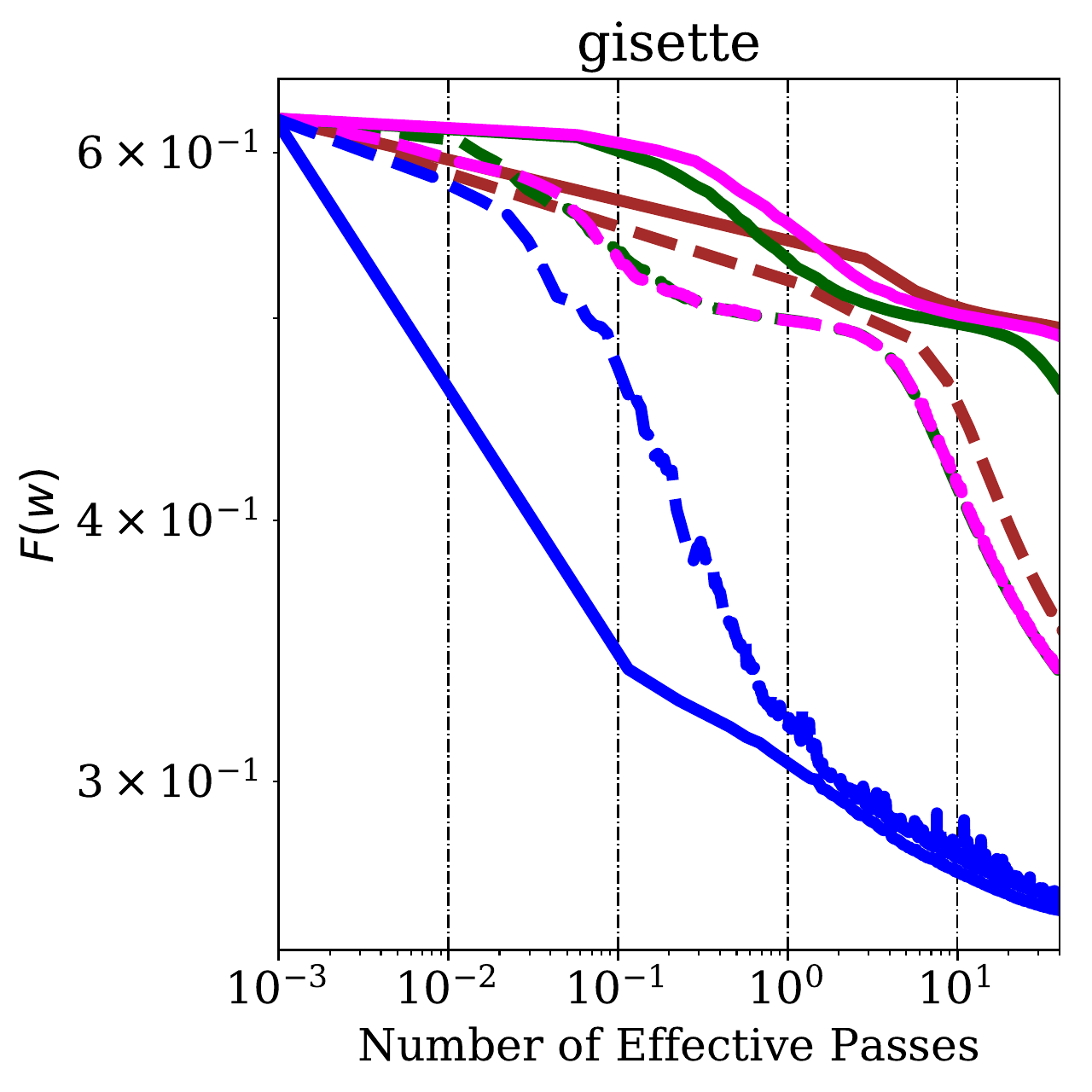}
	
	\includegraphics[width=0.48\textwidth]{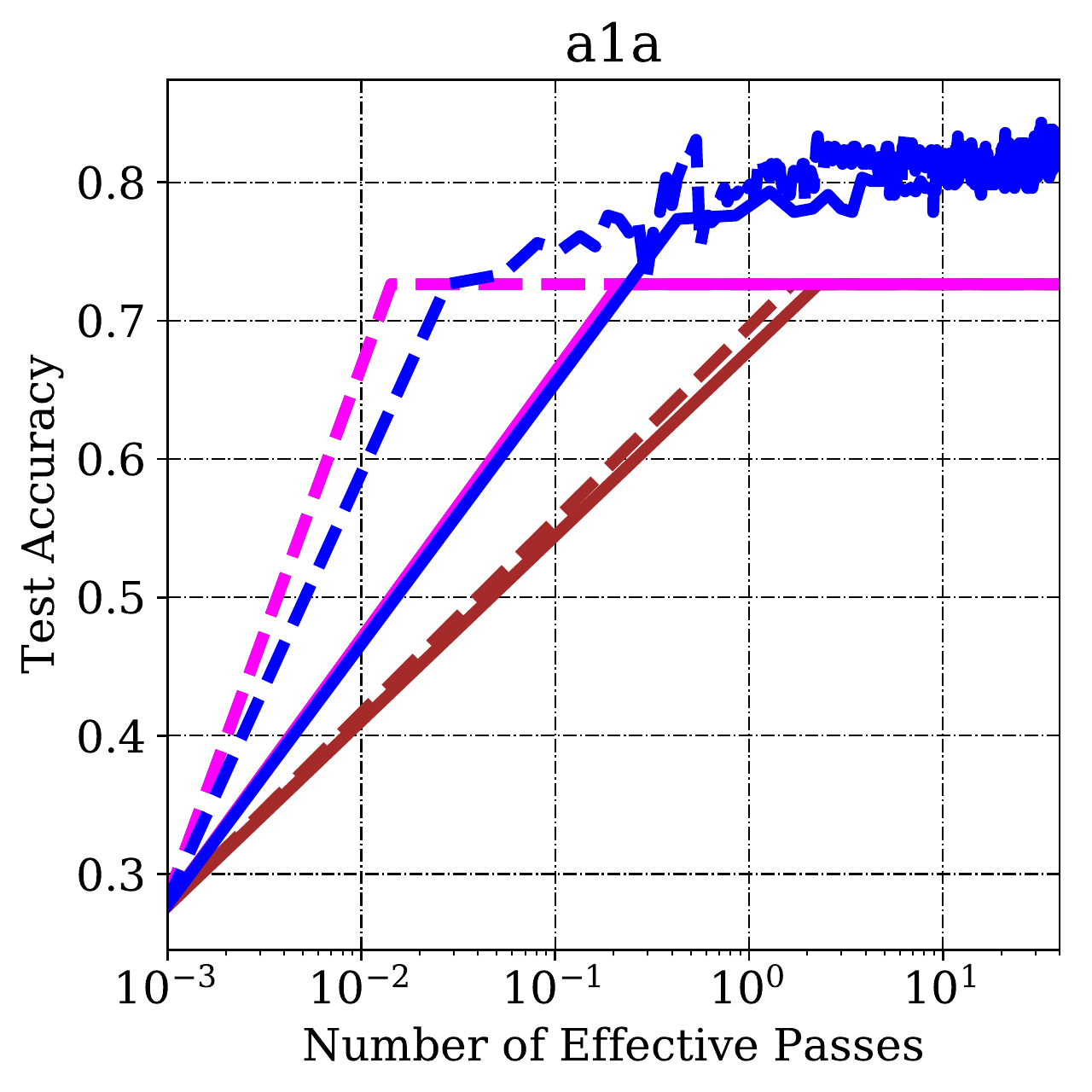}
	\includegraphics[width=0.48\textwidth]{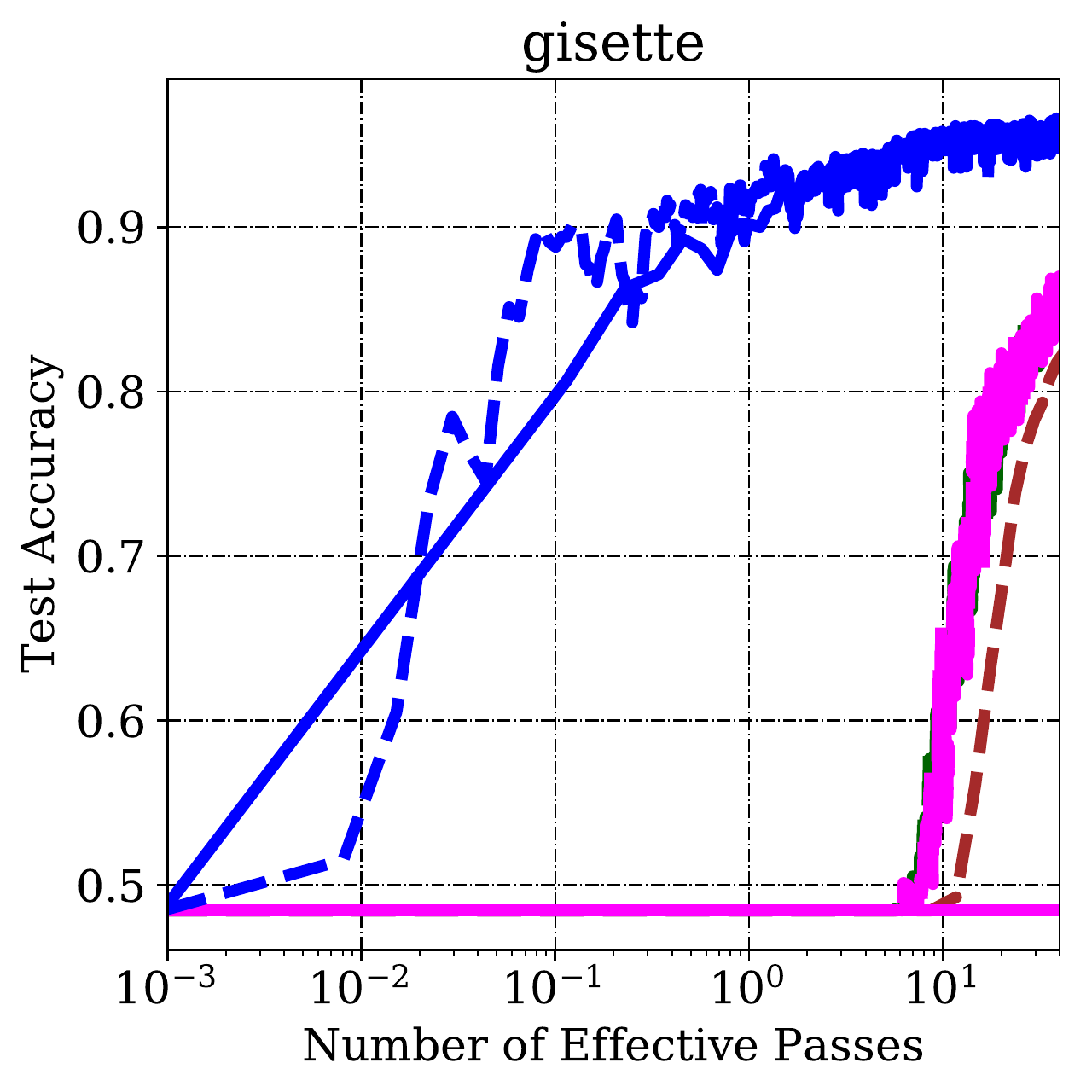}
  \caption{Comparison of objective function ($F(w)$) and Test  Accuracy for different algorithms on Non-Linear Least Squares Problems.}
  \label{fig:StochNonConv_LIBSVM}
\end{subfigure}
\caption{Stochastic Problems. Datasets: \texttt{a1a} and \texttt{gistte}.}
\label{fig:fig}
\end{figure}

\vspace{-5pt}
\paragraph{Stochastic Logistic Regression}
Figure \ref{fig:StochSConv_LIBSVM} shows the performance of the stochastic methods on logistic regression problems. We show results for every method in the small batch regime ($16$) and in the large batch regime ($256$). As is clear, the stochastic variant of \SONIA{} is competitive with the other methods. We should also mention that as predicted by the theory, using a larger batch size (lower variance in the stochastic gradient approximation) allows for \SONIA{} to convergence to a smaller neighborhood around the optimal solution. For more results see Section \ref{sec:moreStochStCOnvResults}.
\vspace{-11pt}
\paragraph{Stochastic Non-Linear Least Square}
Figure \ref{fig:StochNonConv_LIBSVM} shows the performance of the stochastic methods on (nonconvex) nonlinear least squares problems. As is clear, the stochastic variant of \SONIA{} outperforms the other methods for all problems reported. Within a very small number of epochs, \SONIA{} is able to achieve high testing accuracy. We attribute the success of \SONIA{} in the stochastic nonconvex regime to the fact that useful curvature information is incorporated in the search direction. For more results see Section \ref{sec:moreStochNonCOnvResults}.
\vspace{-5pt}
\section{Final Remarks and Future Works}\label{sec:fin_rem}
\vspace{-5pt}
This paper describes a deterministic and stochastic variant of a novel optimization method, \SONIA{}, for empirical risk minimization. The method attempts to bridge the gap between first- and second-order methods by computing a search direction that uses a second-order-type update in one subspace, coupled with a scaled steepest descent step in the orthogonal complement. Numerical results show that the method is efficient in both the deterministic and stochastic settings, and theoretical guarantees confirm that \SONIA{} converges to a stationary point for both strongly convex and nonconvex functions.


Future research directions include: $(1)$ developing adaptive memory variants of \SONIA{}, $(2)$ exploring stochastic \SONIA{} variants that use an adaptive number of samples for gradient/Hessian approximations (following the ideas from \cite{byrd2012sample,mokhtari2016adaptive,jahani2018efficient,bollapragada2018adaptive,friedlander2012hybrid}), or that employ variance reduced gradients, and $(3)$  a thorough numerical investigation for deep learning problems to test the limits of the methods.

\subsubsection*{Acknowledgements} This work was partially supported by the U.S. National Science Foundation, under award numbers NSF:CCF:1618717, NSF:CMMI:1663256 and NSF:CCF:1740796.

\bibliographystyle{spmpsci}
\bibliography{references}

\newpage

\appendix

\section{Theoretical Results and Proofs}\label{sec:app_theory}

\subsection{Assumptions}

\begin{customassum}{\ref{assum:diff}} $F$ is twice continuously differentiable.
\end{customassum}

\begin{customassum}{\ref{assum:strong_conv}} There exist positive constants $\mu$ and $L$ such that
\begin{align*}
    \mu I \preceq \nabla^2F(w) \preceq L I, \quad \text{for all } w \in \mathbb{R}^d.
\end{align*}
\end{customassum}

\begin{customassum}{\ref{assum:boundedf}} The function $ F(w)$ is bounded below by a scalar $\widehat F$.
\end{customassum}

\begin{customassum}{\ref{assum:lip}} The gradients of $F$ are $L$-Lipschitz continuous for all $w \in \mathbb{R}^d$.
\end{customassum}

\begin{customassum}{\ref{assum:boundedStochGrad}} There exist a  constant $\gamma$ such that $\mathbb{E}_{\mathcal{I}}[\|\nabla F_{\mathcal{I}}(w) - \nabla F(w)\|^2] \leq \gamma^2$.
\end{customassum}
\begin{customassum}{\ref{assum:unbiasedEstimGrad}} $\nabla F_{\mathcal{I}}(w)$ is an unbiased estimator of the gradient, i.e., $\mathbb{E}_{\mathcal{I}}[\nabla F_{\mathcal{I}}(w)] = \nabla F(w)$, where the samples $\mathcal{I}$ are drawn independently.
\end{customassum}

\subsection{Proof of Lemma~\ref{pequiv}}

\begin{customlemma}{\ref{pequiv}}
Let $\A_k \eqdef \V_k |\Lambda_k|_{\epsilon}^{-1} \V_k^T + \rho_k (I-\V_k\V_k^T)$. Then the search direction $p_k$ in \eqref{p} is equivalent to $p_k = -\A_k \nabla F(w_k)$.
\end{customlemma}

\begin{proof}
Define
\begin{equation}\label{A1A2}
   A_k \eqdef \V_k |\Lambda_k|_{\epsilon}^{-1} \V_k^T,\quad \text{ and } \quad \Ap_k \eqdef \rho_k (I-\V_k\V_k^T),
\end{equation}
so that by \eqref{eq:A_k_matrix} and \eqref{A1A2},
\begin{equation}\label{A}
    \A_k = A_k + \Ap_k.
\end{equation}
By \eqref{gradF}, $\gp_k = (I-\V\V^T)\gp_k$. Then
\begin{eqnarray}
    p_k &\overset{\eqref{p}+\eqref{A1A2}}{=}& -A_k g_k - \Ap_k\gp_k \nonumber \\
    &\overset{\eqref{gradF}}{=}& -A_k g_k - \Ap_k(\nabla F(w_k) - g_k)\nonumber \\
    &\overset{\eqref{gradF}}{=}& -A_k\V_k\V_k^T\nabla F(w_k) - \Ap_k(I-\V_k\V_k^T)\nabla F(w_k) \nonumber \\
    &\overset{\eqref{A1A2}}{=}& -\V_k |\Lambda_k|_{\epsilon}^{-1} \V_k^T\V_k\V_k^T\nabla F(w_k) - \rho_k(I-\V_k\V_k^T)(I-\V_k\V_k^T)\nabla F(w_k) \nonumber \\
    &=& -\V_k |\Lambda_k|_{\epsilon}^{-1} \V_k^T\nabla F(w_k) - \rho_k(I-\V_k\V_k^T)\nabla F(w_k) \nonumber \\
    &\overset{\eqref{A1A2}+\eqref{A}}{=}& -\mathcal{A}_k \nabla F(w_k).
\end{eqnarray}
\end{proof}

\subsection{Proof of Lemma \ref{lem.pos_def}}

\begin{customlemma}{\ref{lem.pos_def}} The matrix $\A_k$ in \eqref{A} is positive definite.
\end{customlemma}

\begin{proof} Let $y \in \mathbb{R}^d$ be any nonzero vector. One can then write $y = y_1 + y_2$, where $y_1 \in {\rm range}\{\tilde{V}\}$ and $y_2 \in  {\rm range}\{\tilde{V}\}^\perp$. Note that this implies that $\V^Ty_2 = 0$ and $(I - \V\V^T)y_1 = 0$. Also notice that $I - \V\V^T$ is a projection matrix, thus $(I - \V\V^T)^2 = I - \V\V^T$.

Let $D_1 = |\Lambda_k|_{\epsilon}^{-1}$ and $D_2 = \rho I_d$, where $\rho>0$. Thus, we have
\begin{align*}
    y^T \A_k y & = (y_1+y_2)^T \A_k(y_1+y_2) \\
    & = (y_1+y_2)^T \left(A_k+ \Ap_k\right)(y_1+y_2)
    \\
    & = (y_1+y_2)^T \left(\V_k D_1 \V_k^T+D_2(I - \V\V^T) \right)(y_1+y_2)\\
    & = (y_1+y_2)^T \left(\V_k D_1 \V_k^T+(I - \V\V^T)D_2(I - \V\V^T) \right)(y_1+y_2)\\
    & = y_1^T \V_k D_1 \V_k^Ty_1+y_2^T(I - \V\V^T)D_2(I - \V\V^T) y_2\\
    & > 0,
\end{align*}
The final strict inequality is due to the fact that since $y$ is a nonzero vector, this implies that at least one of the two vectors $y_1$ or $y_2$ are also nonzero, the specific decompositions of $y$ and the fact that both matrices $D_1$ and $D_2$ are positive definite.
\end{proof}

\subsection{Proof of Lemma \ref{lem:1}}


\begin{customlemma}{\ref{lem:1}} Let Assumption \ref{assum:diff} hold, and let $\rho_k$ satisfy \eqref{rhok} for all $k\geq 0$. Then there exist constants $0 < \mu_1 \leq \mu_2$ such that the sequence of matrices $\{ \mathcal{A}_k\}_{k\geq 0}$ generated by Algorithm \ref{alg:fastsr1} satisfies,
\begin{align*}   
\mu_1 I \preceq \mathcal{A}_k \preceq \mu_2 I,\qquad \text{for } k=0,1,2,\dots.
\end{align*}
\end{customlemma}







\begin{proof}
By \eqref{A}, $\mathcal{A}_k = A_k+\Ap_k$, so we can write
\begin{eqnarray*}
\mathcal{A}_k &=& A_k+\Ap_k\\
&=& \V_k |\Lambda_k|_{\epsilon}^{-1} \V_k^T +  \rho_k (I-\V_k\V_k^T)\\
&=& \V_k (|\Lambda_k|_{\epsilon}^{-1} -\rho_k I)\V_k^T +  \rho_k I\\
&\overset{\eqref{rhok}}{\succeq}& \rho_k I.
\end{eqnarray*}

Furthermore, $|\Lambda_k|_{\epsilon}^{-1} \preceq \tfrac1{\epsilon} I$ so that $\mathcal{A}_k  \preceq \tfrac1{\epsilon}I + \rho_k I \preceq \tfrac2{\epsilon} $. Defining $\mu_1 \eqdef \min_{0\leq i\leq k} \rho_i$ and $\mu_2 \eqdef \tfrac2{\epsilon}$, and noting that $0< \mu_1 \leq \mu_2$, gives the result.


\end{proof}

\subsection{Proof of Theorem \ref{thm:detSC}}

\begin{customthm}{\ref{thm:detSC}}
Suppose that Assumptions \ref{assum:diff} and \ref{assum:strong_conv} hold, and let $F^{\star} = F(w^{\star})$, where $w^{\star}$ is the minimizer of $F$. Let $\{w_k\}$ be the iterates generated by Algorithm~\ref{alg:fastsr1}, where $0 <  \alpha_k = \alpha  \leq \frac{\mu_1}{\mu_2^2 L}$, and $w_0$ is the starting point. Then, for all $k\geq 0$, 
\begin{align*}   
    F(w_k) - F^{\star}& \leq \big( 1-\alpha \mu \mu_1 \big)^k \left[ F(w_0) - F^{\star} \right].
\end{align*}
\end{customthm}

\begin{proof} We have that
\begin{align} 
    F(w_{k+1}) & = F(w_k -\alpha \mathcal{A}_k \nabla F(w_k))  \nonumber\\
    & \leq F(w_k) + \nabla F(w_k)^T (-\alpha \mathcal{A}_k \nabla F(w_k)) + \frac{L}{2}\| \alpha \mathcal{A}_k \nabla F(w_k)\|^2  \nonumber\\
    & \leq F(w_k) - \alpha \mu_1 \| \nabla F(w_k) \|^2  + \frac{\alpha^2 \mu_2^2 L}{2} \| \nabla F(w_k)\|^2 \nonumber\\
    & = F(w_k) - \alpha \left(\mu_1 - \alpha \frac{\mu_2^2 L}{2} \right)\| \nabla F(w_k) \|^2  \nonumber\\
    & \leq F(w_k) - \alpha \frac{\mu_1}{2} \| \nabla F(w_k) \|^2,  \label{eq:proof_step_app}
\end{align}
where the first inequality is due to Assumption \ref{assum:strong_conv}, the second inequality arises as a consequence of Lemma \ref{lem:1} and the last inequality is due to the choice of the steplength. By strong convexity \cite{nesterov2013introductory}, we have $2\mu(F(w)-F^\star) \leq \| \nabla F(w)\|^2$, and thus
\begin{align*}
    F(w_{k+1}) \leq F(w_k) - \alpha \mu \mu_1 (F(w_k)-F^\star).
\end{align*}
Subtracting $F^\star$ from both sides, 
\begin{align*}
    F(w_{k+1}) - F^\star \leq (1 - \alpha \mu \mu_1) (F(w_k)-F^\star).
\end{align*}
Recursive application of the above inequality yields the desired result.
\end{proof}

\subsection{Proof of Theorem \ref{thm_non_MB}}

\begin{customthm}{\ref{thm_non_MB}}
Suppose that Assumptions \ref{assum:diff},  \ref{assum:boundedf} and \ref{assum:lip} hold. Let $\{w_k\}$ be the iterates generated by Algorithm~\ref{alg:fastsr1}, where
$0 <  \alpha_k = \alpha  \leq \frac{\mu_1}{\mu_2^2  L}$,
and $w_0$ is the starting point. Then, for any $T > 1$,
\begin{align*}	
\frac{1}{T}\sum_{k=0}^{T-1} \| \nabla F(w_k) \|^2  &  \leq \frac{2[ F(w_0) - \hat{F}]}{\alpha \mu_1 T } \xrightarrow[]{T \rightarrow \infty}0.
\end{align*}
\end{customthm}

\begin{proof} We start with \eqref{eq:proof_step_app} \begin{align*} 
F(w_{k+1})  \leq F(w_k) - \alpha \frac{\mu_1}{2} \| \nabla F(w_k) \|^2.  
\end{align*}
Summing both sides of the above inequality from $k=0$  to $T -1$,
\begin{align*} 
	\sum_{k=0}^{T-1}(F(w_{k+1}) - F(w_k)) \leq  - \sum_{k=0}^{T-1}\alpha \frac{\mu_1}{2} \| \nabla F(w_k) \|^2. 
\end{align*}
The left-hand-side of the above inequality is a telescopic sum and thus,
\begin{align*}
	\sum_{k=0}^{T-1} \left[F(w_{k+1})-F(w_k) \right]  &= F(w_{T})-F(w_0)
	\geq \widehat{F} -F(w_0),
\end{align*}
where the inequality is due to $\hat{F} \leq F(w_T)$ (Assumption \ref{assum:boundedf}). Using the above, we have
\begin{align}   \label{eq:int_nonconvex}
    \sum_{k=0}^{T-1} \| \nabla F(w_k) \|^2  &  \leq \frac{2[ F(w_0) - \widehat{F}]}{\alpha \mu_1}.
\end{align}
Dividing \eqref{eq:int_nonconvex} by $T$ we conclude
\begin{align*}
    \frac{1}{T}\sum_{k=0}^{T-1} \| \nabla F(w_k) \|^2  &  \leq \frac{2[ F(w_0) - \widehat{F}]}{\alpha \mu_1 T}.
\end{align*}
\end{proof}

\subsection{Proof of Theorem \ref{thm:stochSC}}

\begin{customthm}{\ref{thm:stochSC}}
Suppose that Assumptions \ref{assum:diff}, \ref{assum:strong_conv}, \ref{assum:boundedStochGrad} and \ref{assum:unbiasedEstimGrad} hold, and let $F^{\star} = F(w^{\star})$, where $w^{\star}$ is the minimizer of $F$. Let $\{w_k\}$ be the iterates generated by Algorithm~\ref{alg:fastsr1}, where $0 <  \alpha_k = \alpha  \leq \frac{\mu_1}{\mu_2^2 L}$, and $w_0$ is the starting point. Then, for all $k\geq 0$, 
\begin{align*}   
  \mathbb{E}[F(w_{k}) - F^\star] & \leq  \left( 1 - \alpha \mu_1\lambda \right)^k \left( F(w_0) - F^\star - \frac{\alpha \mu_2^2 \gamma^2 L}{2\mu_1 \lambda}\right) +  \frac{\alpha \mu_2^2 \gamma^2 L}{2\mu_1 \lambda}.
\end{align*}
\end{customthm}

\begin{proof} We have that
\begin{align*} 
    F(w_{k+1}) & = F(w_k -\alpha \mathcal{A}_k \nabla F_{\mathcal{I}_k}(w_k))  \nonumber\\
    & \leq F(w_k) + \nabla F(w_k)^T (-\alpha \mathcal{A}_k \nabla F_{\mathcal{I}_k}(w_k)) + \frac{L}{2}\| \alpha \mathcal{A}_k \nabla F_{\mathcal{I}_k}(w_k)\|^2  \nonumber\\
    & \leq F(w_k) - \alpha\nabla F(w_k)^T  \mathcal{A}_k \nabla F_{\mathcal{I}_k}(w_k) + \frac{\alpha^2 \mu_2^2 L}{2}\|\nabla F_{\mathcal{I}_k}(w_k)\|^2
\end{align*}
where the first inequality is due to Assumption \ref{assum:strong_conv} and the second inequality is due to Lemma \ref{lem:1}. Taking the expectation over the sample $\mathcal{I}_k$, we have
\begin{align} 
    \mathbb{E}_{\mathcal{I}_k}[F(w_{k+1})] & \leq F(w_k) - \alpha \mathbb{E}_{\mathcal{I}_k}[ \nabla F(w_k)^T  \mathcal{A}_k \nabla F_{\mathcal{I}_k}(w_k)] + \frac{\alpha^2 \mu_2^2 L}{2}\mathbb{E}_{\mathcal{I}_k}[\|\nabla F_{\mathcal{I}_k}(w_k)\|^2] \nonumber\\
    & = F(w_k) - \alpha  \nabla F(w_k)^T  \mathcal{A}_k \nabla F(w_k) + \frac{\alpha^2 \mu_2^2 L}{2}\mathbb{E}_{\mathcal{I}_k}[\|\nabla F_{\mathcal{I}_k}(w_k)\|^2] \nonumber\\
    &\leq F(w_k) - \alpha \left( \mu_1 - \frac{\alpha \mu_2^2 L}{2}\right) \|\nabla F(w_k)\|^2 + \frac{\alpha^2 \mu_2^2 \gamma^2 L}{2} \nonumber\\
    &\leq F(w_k) - \frac{\alpha \mu_1}{2} \|\nabla F(w_k)\|^2 + \frac{\alpha^2 \mu_2^2 \gamma^2 L}{2}, \label{eq.key1}
\end{align}
where the second inequality is due to Lemma \ref{lem:1} and Assumption \ref{assum:boundedStochGrad} and the third inequality is due to the choice of the step length.
Since $F$ is strongly convex \cite{nesterov2013introductory}, we have
\begin{align*} 
    \mathbb{E}_{\mathcal{I}_k}[F(w_{k+1})] & \leq F(w_k) - \alpha \mu_1 \lambda (F(w_k) - F^\star) + \frac{\alpha^2 \mu_2^2 \gamma^2 L}{2}.
\end{align*}
Taking the total expectation over all batches $\mathcal{I}_0$, $\mathcal{I}_1$, $\mathcal{I}_2$,... and all history starting with $w_0$, and subtracting $F^\star$ from both sides, we have
\begin{align*}
    \mathbb{E}[F(w_{k+1}) - F^\star] & \leq \left( 1 - \alpha \mu_1\lambda \right) \mathbb{E}[F(w_k) - F^\star] + \frac{\alpha^2 \mu_2^2 \gamma^2 L}{2},
\end{align*}
where $0 \leq \left( 1 - \alpha \mu_1 \lambda\right) \leq 1$ by the step length choice. Subtracting $\frac{\alpha \mu_2^2 \gamma^2 L}{2\mu_1 \lambda}$ from both sides yields
\begin{align*}
    \mathbb{E}[F(w_{k+1}) - F^\star] - \frac{\alpha \mu_2^2 \gamma^2 L}{2\mu_1 \lambda} & \leq  \left( 1 - \alpha \mu_1\lambda \right) \left( \mathbb{E}[F(w_k) - F^\star] - \frac{\alpha \mu_2^2 \gamma^2 L}{2\mu_1 \lambda}\right).
\end{align*}
Recursive application of the above completes the proof.
\end{proof}

\subsection{Proof of Theorem \ref{thm:stoch_non}}

\begin{customthm}{\ref{thm:stoch_non}}
Suppose that Assumptions \ref{assum:diff}, \ref{assum:boundedf}, \ref{assum:lip}, \ref{assum:boundedStochGrad} and \ref{assum:unbiasedEstimGrad} hold, and let $F^{\star} = F(w^{\star})$, where $w^{\star}$ is the minimizer of $F$. Let $\{w_k\}$ be the iterates generated by Algorithm~\ref{alg:fastsr1}, where $0 <  \alpha_k = \alpha  \leq \frac{\mu_1}{\mu_2^2 L}$, and $w_0$ is the starting point. Then, for all $k\geq 0$, 
\begin{align*}   
  \mathbb{E}\left[\frac{1}{T}\sum_{k=0}^{T-1} \| \nabla F(w_k) \|^2 \right]  &  \leq \frac{2[ F(w_0) - \widehat{F}]}{\alpha \mu_1 T} + \frac{\alpha \mu_2^2 \gamma^2 L }{\mu_1} \xrightarrow[]{T \rightarrow \infty} \frac{\alpha \mu_2^2 \gamma^2 L }{\mu_1}.
\end{align*}
\end{customthm}

\begin{proof} Starting with \eqref{eq.key1} and taking the total expectation over all batches $\mathcal{I}_0$, $\mathcal{I}_1$, $\mathcal{I}_2$,... and all history starting with $w_0$ 
\begin{align*} 
    \mathbb{E}[F(w_{k+1}) - F(w_k)] 
    &\leq  - \frac{\alpha \mu_1}{2} \mathbb{E}[\|\nabla F(w_k)\|^2] + \frac{\alpha^2 \mu_2^2 \gamma^2 L}{2}.
\end{align*}
Summing both sides of the above inequality from $k=0$  to $T -1$,
\begin{align*} 
    \sum_{k=0}^{T-1}\mathbb{E}[F(w_{k+1}) - F(w_k)] 
    &\leq  - \frac{\alpha \mu_1}{2}\sum_{k=0}^{T-1} \mathbb{E}[\|\nabla F(w_k)\|^2] + \frac{\alpha^2 \mu_2^2 \gamma^2 L T}{2} \\
    & = - \frac{\alpha \mu_1}{2}\mathbb{E}\left[\sum_{k=0}^{T-1} \|\nabla F(w_k)\|^2\right] + \frac{\alpha^2 \mu_2^2 \gamma^2 L T}{2}.
\end{align*}
The left-hand-side of the above inequality is a telescopic sum and thus,
\begin{align*}
	\sum_{k=0}^{T-1} \mathbb{E}\left[F(w_{k+1})-F(w_k) \right]  &= \mathbb{E}[F(w_{T})]-F(w_0)
	\geq \widehat{F} -F(w_0),
\end{align*}
where the inequality is due to $\hat{F} \leq F(w_T)$ (Assumption \ref{assum:boundedf}). Using the above, we have
\begin{align*}  
    \mathbb{E}\left[\sum_{k=0}^{T-1} \| \nabla F(w_k) \|^2 \right]  &  \leq \frac{2[ F(w_0) - \widehat{F}]}{\alpha \mu_1} + \frac{\alpha \mu_2^2 \gamma^2 L T}{\mu_1}.
\end{align*}
Dividing by $T$ we conclude completes the proof.
\end{proof}

\clearpage
\section{Efficient Hessian-Matrix Computations}\label{sec:Hess_mat}

The \SONIA{} algorithm requires the computation of Hessian-maxtrix products for the construction of the curvature pairs. In this section, we describe how one can efficiently compute Hessian-matrix products for the problems studied in this paper. Moreover, we describe an efficient distributed algorithm for computing curvature pairs; see \cite{jahani2019scaling} for more details.

Assume that $\odot$ is the operator for component-wise product, $*$ is the standard multiplication of matrices, $\mathbbm{1}_n$ is the vector of ones with size $1\times n$, $X$ is the feature matrix and $Y$ is the label matrix. 

In the following, firstly, we present the efficient calculation of objective function, gradient and Hessian-matrix products for logistic regression problems. Next, we do the same for non-linear least square problems. Moreover, we describe a simple distributed methodology for computing Hessian-matrix products. Finally, further discussion is provided for the efficient computation of Hessian-matrix products.   

\subsection{Logistic Regression}
The objective function, gradient, Hessian and Hessian-matrix product for logistic regression problems are calculated efficiently as follows:
\begin{equation}
    F(w) = \dfrac{1}{n} \left( \mathbbm{1}_n * \log\left(1 + e^{-Y\odot X^Tw}\right)   \right) +\dfrac{\lambda}{2} \|w\|^2,
\end{equation}\label{eq:objlogReg}

\begin{equation}
   \nabla F(w) = \dfrac{1}{n} \left( X^T * \dfrac{-Y \odot e^{-Y\odot X^Tw}}{1 + e^{-Y\odot X^Tw}} \right) + \lambda w,
\end{equation}\label{eq:gradlogReg}

\begin{equation}
   \nabla^2 F(w) = \dfrac{1}{n} \left( X^T *\left[\dfrac{Y\odot Y \odot e^{-Y\odot X^Tw}}{\big( 1 + e^{-Y\odot X^Tw}\big)^2}\right]\odot X \right) + \lambda I_d,
\end{equation}\label{eq:hesslogReg}

\begin{equation}
   \nabla^2 F(w)*S = \dfrac{1}{n} \left( X^T *\left[\dfrac{Y\odot Y \odot e^{-Y\odot X^Tw}}{\big( 1 + e^{-Y\odot X^Tw}\big)^2}\right]\odot X*S \right) + \lambda S.
\end{equation}\label{eq:hess_mat_logReg}

\subsection{Non-Linear Least Squares}
The objective function, gradient, Hessian and Hessian-matrix product for non-linear least square are calculated efficiently as follows (where $\phi(z) = \dfrac{1}{1+e^{-z}}$):
\begin{equation}
   F(w) = \dfrac{1}{2n}\|Y - \phi(X^Tw)\|^2
\end{equation}\label{eq:objNLLS}

\begin{equation}
   \nabla F(w) = \dfrac{1}{n}\left( -X^T *[\phi(X^Tw) \odot(1- \phi(X^Tw))\odot(Y - \phi(X^Tw))] \right)
\end{equation}\label{eq:gradNLLS}

\begin{equation}
\label{eq:hessNLLS}
\begin{split}
   \nabla^2 F(w) =& \dfrac{1}{n} \left( -X^T *\left[\phi(X^Tw)\odot(1-\phi(X^Tw)) \odot(Y - 2(1+Y)\odot\phi(X^Tw) \right.\right.+\\ 
    &\left. \left. 3\phi(X^Tw)\odot\phi(X^Tw))\right]\odot X \right) 
\end{split}
\end{equation}

\begin{equation}
\begin{split}\label{eq:hess_mat_NLLS}
       \nabla^2 F(w)*S =& \dfrac{1}{n} \left( -X^T *\left[\phi(X^Tw)\odot(1-\phi(X^Tw)) \odot(Y - 2(1+Y)\odot\phi(X^Tw) \right.\right. +\\ 
       &\left. \left. 3\phi(X^Tw)\odot\phi(X^Tw))\right]\odot X * S \right) 
\end{split}
\end{equation}
Based on the above equations, one can note that the cost of objective function and gradient computations is $\mathcal{O}(nd)$ due to the calculation of $X^Tw$. The cost of Hessian-matrix products is $\mathcal{O}(mnd)$, and the cost is dominated by the calculation of $X*S$. Furthermore, by considering the fact that \SONIA{} was developed for the regime where $m \ll d,n$, the effective cost of Hessian-matrix product is $\mathcal{O}(nd)$. In summary, the cost of Hessian-matrix products is similar to the cost of objective function and gradient evaluations. 

\subsection{Distributed Algorithm}
For the cases where the Hessian has a compact representation (e.g., logistic regression and non-linear least square problems), we showed that the Hessian-matrix products can be efficiently calculated. However, this is not always the case. In the rest of this section, we justify that for general non-linear problems, such as deep learning, the Hessian-matrix product can be efficiently computed in a distributed environment. By following the study in \cite{jahani2019scaling}, one can note that in order to construct curvature information $(S_k, Y_k)$, a Hessian-matrix calculation is required in order to form $Y_k = \nabla^2 F(w_k)S_k$. The aforementioned Hessian-matrix products can be calculated efficiently in master-worker framework, summarized in Algorithm \ref{alg:calSYdist}. Each worker has a portion of the dataset, performs local computations, and then reduces the locally calculated information to the master node. This method is matrix-free (i.e., the Hessian approximation is never explicitly constructed).

\begin{algorithm}[]
{ 
\small
\caption{Construct new $(S_k,Y_k)$ curvature pairs}
  \label{alg:calSYdist}
 {\bf Input:} $w_k$ (iterate), $m$ (memory), $S_k = [\;]$, $Y_k =[\;]$ (curvature pair containers). 

\textbf{Master Node:} \hfill \textbf{Worker Nodes (\pmb{$i=1,2,\dots,\mathcal{K}$}):} 
  \begin{algorithmic}[1]
  \State {\color{green!40!black}{\it \textbf{Broadcast:}}} $S_k$ and $w_k$ \hfill  \textcolor{green!40!black}{$\pmb{\longrightarrow}$}  \hfill Compute   $Y_{k,i} = \nabla^2 F_{i}(w_k)S_k$
  \State {\color{orange!50!black}{\it \textbf{Reduce:}}} $Y_{k,i}$ to  $Y_k$ and $S_k^TY_{k,i}$ to  $Y_k^T S_k$ \hfill  \textcolor{orange!50!black}{$\pmb{\longleftarrow}$} \hfill Compute $S_k^TY_{k,i}$
  \end{algorithmic}
  {\bf Output:} $S^TY$, $Y_k$
  }
\end{algorithm}


\clearpage

\section{Method and Problem Details}\label{apndx:addNumExpDet}

\subsection{Table of Algorithms}\label{sec:tbl_algs}
In this section, we summarize the implemented algorithms in Section \ref{sec:num_res} in the Table \ref{tab:descAlgs}.
\begin{table}[H]\scriptsize
	\centering
	\begin{tabular}{ll} \toprule
		\textbf{Algorithm} & \textbf{Description and Reference}  \\ \midrule
		{\textcolor{cyan}{NEST+}}  & {Algorithm described in \cite[Chapter 2]{Nesterov04} with adaptive Lipschitz constant}   \\\hdashline
        {\textcolor{green}{L-BFGS}}  & {Limited memory BFGS} \cite{liu1989limited} \\\hdashline
		{\textcolor{orange}{L-SR1}}  & {Limited memory SR1} \cite{lu1996study}\\\hdashline
		{\textcolor{purple}{Newton-CG-TR}}   & {Newton method with conjugate gradient (CG) utilizing trust region (TR)\cite{nocedal_book}} \\\hdashline
		{\textcolor{purple}{Newton-CG-LS}}  & {Newton method with conjugate gradient (CG) utilizing line search (LS) \cite{nocedal_book}}  \\\hdashline
		{\textcolor{brown}{SARAH+}}  &  {Practical variant of SARAH \cite{Nguyen2017}}     \\\hdashline
		{\textcolor{magenta}{SQN}} &  {Stochastic Quasi-Newton \cite{byrd2016stochastic}} \\\hdashline
		{\textcolor{green!30!black}{SGD}} &  {Stochastic gradient method \cite{robbins1951stochastic}} \\\hdashline
		{\textcolor{red}{GD}} &  {Gradient descent} \\\hdashline
		{ASUESA}  & {Accelerated Smooth Underestimate Sequence Algorithm with adaptive Lipschitz constant \cite{ma2017underestimate}} \\\hdashline
		{\textcolor{blue}{SONIA}} &  {\textbf{S}ymmetric bl\textbf{O}ckwise tru\textbf{N}cated optim\textbf{I}nation \textbf{A}lgorithm} \\
		\midrule
	\end{tabular}
	\normalsize
	\caption{Description of implemented algorithms}
\label{tab:descAlgs}
\end{table}
In order to find $w^\star$ for the strongly convex problems, we used the ASUESA algorithm \cite{ma2017underestimate}. ASUESA constructs a sequence of lower bounds, and at each iteration of ASUESA the gap between the aforementioned lower bounds and objective function goes to zero at an optimal linear rate. One of the most important advantages of ASUESA is the natural stopping condition, which provides the user with a certificate of optimality. In other words, when the gap between objective function and the lower bounds is small enough, ASUESA is close enough to the optimal solution $w^\star$.

\subsection{Problem Details}\label{apndx:Prob_Details}

\begin{table}[htb]
	\centering
	\caption{Summary of two binary classification datasets and two multi-labels classification datasets}
	\begin{tabular}{lccc}
		\toprule
		\textbf{Dataset} & \textbf{\# of samples } & \textbf{\# of features}& \textbf{\# of categories} \\
		\midrule
		\texttt{rcv1}& 20,242& 47,326 & 2\\\hdashline
		\texttt{gisette}& 6000& 5,000 & 2\\\hdashline
		\texttt{a1a}& 1605& 119& 2\\\hdashline
		\texttt{ijcnn1}& 35000 & 22& 2\\
		\bottomrule
	\end{tabular}
	
	\label{tab:datasets}
\end{table}

\subsection{Implementation Details}\label{sec:impDetndTuning}
In the following sections, we describe the way that the algorithms (Table \ref{tab:descAlgs}) were implemented and tuned. In order to have fair comparisons, we considered the same number of hyper-parameter choices for a given problem for the algorithms that needed tuning. We consider the set of hyper-parameters for each single algorithm for the optimization problems discussed in Section \ref{sec:num_res}.

\subsubsection{Deterministic Strongly Convex Case}
\begin{itemize}
    \item \textbf{\textcolor{red}{GD}}: No tuning is needed. The learning rate is chosen by Armijo backtracking line search. 
    \item \textbf{\textcolor{green}{L-BFGS}}: Memory is chosen from the set $\{4,16,32,64\}$; $\epsilon_{\text{L-BFGS}} = 10^{-8}$ ($\epsilon$ for checking the curvature condition in L-BFGS method) and the learning rate is chosen by Armijo backtracking line search.
    \item \textbf{\textcolor{orange}{L-SR1}}: Memory is chosen from the set $\{4,16,32,64\}$ and $\epsilon_{\text{L-SR1}} = 10^{-8}$($\epsilon$ for checking the curvature condition in L-SR1 method). The search direction is calculated by trust region solver by the default setting reported in Algorithm 6.1 in  \cite{nocedal_book}. 
    \item \textbf{\textcolor{purple}{Newton-CG-LS}}: The learning rate is chosen by Armijo backtracking line search. The Newton system is solved according to Algorithm 7.1 in \cite{nocedal_book}.
    
    \item \textbf{\textcolor{cyan}{NEST+}}: For adaptive Lipschitz constant, we set the parameters $d_{\text{decrease}} \in \{1.1, 2\}$ and $u_{\text{increase}} \in \{1.1, 2\}$ according to the Algorithm 4.1 in \cite{nesterov2013gradient}. 
    \item \textbf{\textcolor{blue}{SONIA}}: Memory is chosen from the set $\{4,16,32,64\}$; $\epsilon_{\text{SONIA}} = 10^{-5}$ (truncated $\epsilon$) and the learning rate is chosen by Armijo backtracking line search.
\end{itemize}
\subsubsection{Deterministic Non-convex Case}
\begin{itemize}
    \item \textbf{\textcolor{red}{GD}}: No tuning is needed. The learning rate is chosen by Armijo backtracking line search. 
    \item \textbf{\textcolor{green}{L-BFGS}}: Memory is chosen from the set $\{4,16,32,64\}$; $\epsilon_{\text{L-BFGS}} = 10^{-8}$($\epsilon$ for checking the curvature condition in L-BFGS method) and the learning rate is chosen by Armijo backtracking line search.
    \item \textbf{\textcolor{orange}{L-SR1}}: Memory is chosen from the set $\{4,16,32,64\}$ and $\epsilon_{\text{L-SR1}} = 10^{-8}$($\epsilon$ for checking the curvature condition in L-SR1 method). The search direction is calculated by trust region solver by the default setting reported in Algorithm 6.1 in  \cite{nocedal_book}. 
    \item \textbf{\textcolor{purple}{Newton-CG-TR}}: The Newton system is solved according to CG-Steihaug method (Algorithm 7.2) in \cite{nocedal_book}.
    
    \item \textbf{\textcolor{blue}{SONIA}}: Memory is chosen from the set $\{4,16,32,64\}$; $\epsilon_{\text{SONIA}} = 10^{-5}$ (truncated $\epsilon$) and the learning rate is chosen by Armijo backtracking line search.
\end{itemize}


\subsubsection{Stochastic Strongly Convex Case}
\begin{itemize}
    \item \textbf{\textcolor{green!30!black}{SGD}}: The learning rate is chosen from the set $\{1, 0.5, 0.1, 0.05, 0.01, 0.005, 0.001\}$ and the batch size is from the set $\{16, 256\}$. 
    \item \textbf{\textcolor{magenta}{SQN}}: The learning rate is chosen from the set $\{1, 0.5, 0.1, 0.05, 0.01, 0.005, 0.001\}$ and the batch size is from the set $\{16, 256\}$. Moreover, we set $L_{\text{SQN}} = 1$, meaning that it checks to accept/reject the curvature information at every iteration. Also, we set $\epsilon_{\text{SQN}} = 10^{-8}$ ($\epsilon$ for checking the curvature condition in SQN method). By checking the sensitivity analysis of SQN w.r.t different memories, we notice SQN is not sensitive to the choice of memory, then we set memory $m=64$. 
    \item \textbf{\textcolor{brown}{SARAH+}}: The learning rate is chosen from the set $\{4, 2, 1, 0.1, 0.01, 0.001\}$\footnote{The reason for this choice is that the learning for SARAH is selected to be approximately by $\dfrac{1}{L}$ where $L$ is the Lipschitz constant.}. Also, we consider the batch sizes 16 and 256.
    \item \textbf{\textcolor{blue}{SONIA}}: The learning rate is chosen from the set $\{1, 0.5, 0.1, 0.05, 0.01, 0.005, 0.001\}$ and the batch size is from the set $\{16, 256\}$. Also, we set $\epsilon_{\text{SONIA}} = 10^{-5}$ (truncated $\epsilon$ ) and memory $m=64$. 
\end{itemize}
\subsubsection{Stochastic Non-convex Case}
\begin{itemize}
    \item \textbf{\textcolor{green!30!black}{SGD}}: The learning rate is chosen from the set $\{1, 0.5, 0.1, 0.05, 0.01, 0.005, 0.001\}$.  
    \item \textbf{\textcolor{magenta}{SQN}}: The tuning for this case is similar to the stochastic strongly convex case. The candidate learning rate set is $\{1, 0.5, 0.1, 0.05, 0.01, 0.005, 0.001\}$ and the batch size is from the set $\{16, 256\}$. Moreover, we set $L_{\text{SQN}} = 1$. Also, we set $\epsilon_{\text{SQN}} = 10^{-8}$ ($\epsilon$ for checking the curvature condition in SQN method) and memory $m=64$.   
    \item \textbf{\textcolor{brown}{SARAH+}}: The learning rate is chosen from the set $\{1, 0.5, 0.1, 0.05, 0.01, 0.005, 0.001\}$. Also, we consider the batch sizes 16 and 256.
    \item \textbf{\textcolor{blue}{SONIA}}: Similar to the previous case, the learning rate is chosen from the set $\{1, 0.5, 0.1, 0.05, 0.01, 0.005, 0.001\}$ and the batch size is from the set $\{16, 256\}$. Also, we set $\epsilon_{\text{SONIA}} = 10^{-5}$ (truncated $\epsilon$ ) and memory $m=64$.
\end{itemize}

\subsubsection{Required Hardware and Software}
All the algorithms are implemented in Python 3 and ran on Intel(R) Xeon(R) CPUs.

\newpage

\section{Additional Numerical Experiments}\label{apndx:add}
In this section, we present additional numerical results in order to compare the performance of \SONIA{} with the state-of-the-art first- and second-order methods described in Table \ref{tab:descAlgs} on the datasets reported in Table \ref{tab:datasets}. 
\begin{enumerate}
    \item Section \ref{sec:moreStConvResults}: deterministic strongly convex case ($\ell_2$ regularized logistic regression).
    \item Section \ref{sec:moreNonConvResults}: deterministic nonconvex case (non-linear least squares).
    \item Section \ref{sec:moreStochStCOnvResults}: stochastic strongly convex case ($\ell_2$ regularized logistic regression).
    \item Section \ref{sec:moreStochNonCOnvResults}: stochastic  nonconvex (non-linear least squares).
\end{enumerate}


Moreover, we investigated the sensitivity of \SONIA{} to its associated hyper-parameters (i.e., the memory size $m$, and the truncation parameter $\epsilon$). Sections \ref{sec:sensMMR} and \ref{sec:sensEPS} show sensitivity results for deterministic logistic regression problems. The key take-aways are that \SONIA{} is robust with respect to $m$ and $\epsilon$ (the variation in performance is small for different choices of the hyper-parameters) and under reasonable choices of these hyper-parameters ($m\in [0,d]$ and $\epsilon>0$), the \SONIA{} algorithm always converges, albeit at a slower rate for some choices. This of course is in contrast to certain methods that may diverge is the hyper-parameters are not chosen appropriately (e.g., the learning rate for the SGD method).


\clearpage
\subsection{Additional Numerical Experiments: Deterministic Strongly Convex Functions}\label{sec:moreStConvResults}\vspace{-10pt}
\begin{figure*}[htb]
	\centering
	{	\includegraphics[trim=10 110 10 110,clip, width=0.6\textwidth]{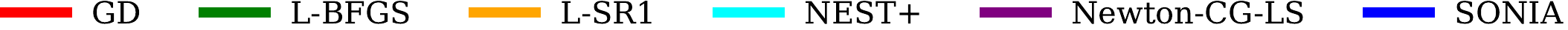}
}\includegraphics[width=0.8\textwidth]{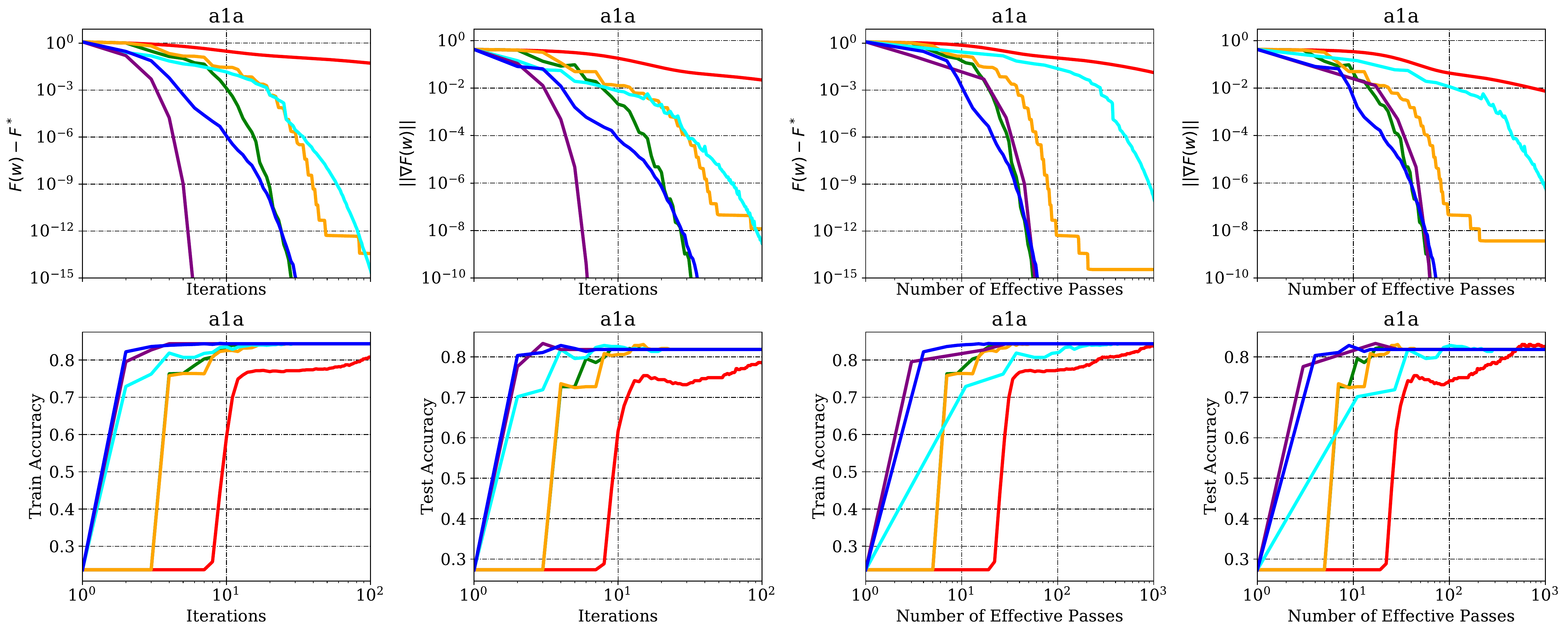}\vspace{-7pt}
	\caption{\texttt{a1a}: Deterministic Logistic Regression with $\lambda = 10^{-3}$.}
	\label{fig:a1a_3}
\end{figure*}
\vspace{-10pt}
\begin{figure*}[htb]
	\centering
	\includegraphics[width=0.8\textwidth]{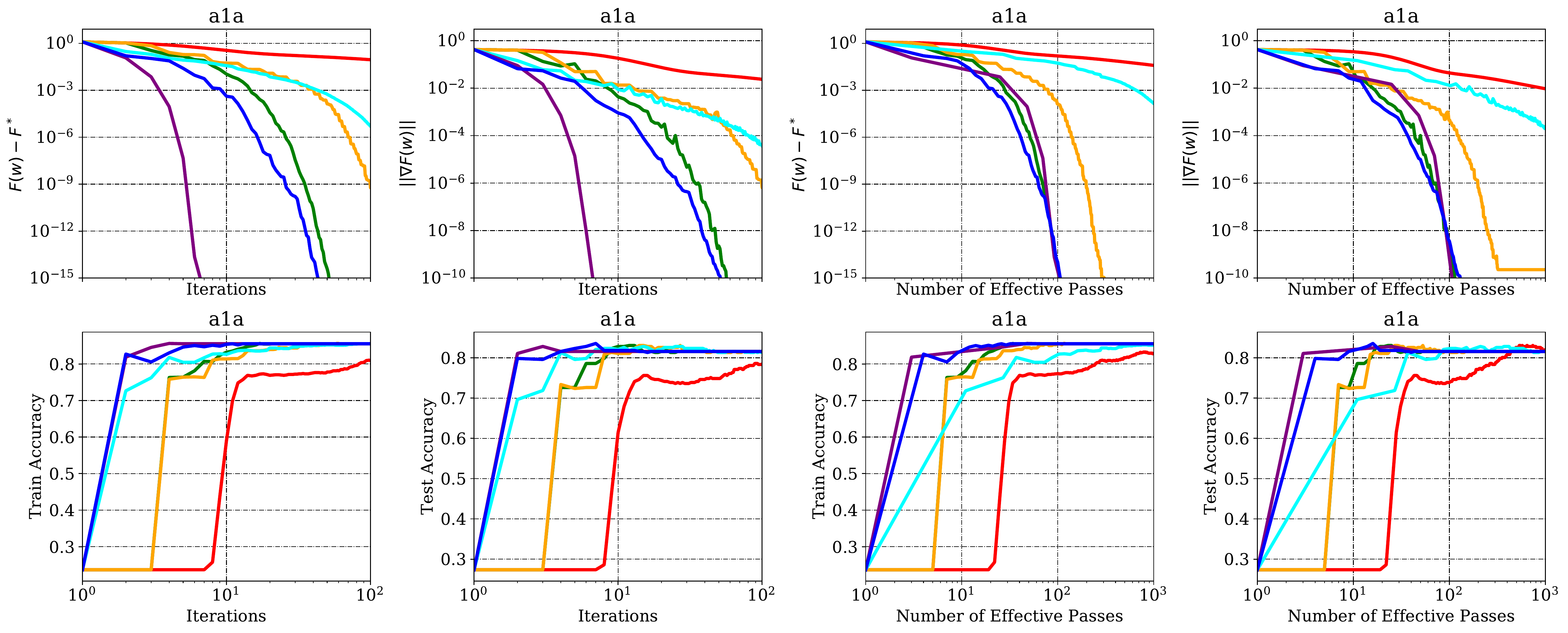}\vspace{-7pt}
	\caption{\texttt{a1a}: Deterministic Logistic Regression with $\lambda = 10^{-4}$.}
	\label{fig:a1a_4}
\end{figure*}
\vspace{-10pt}
\begin{figure*}[htb]
	\centering
	\includegraphics[width=0.8\textwidth]{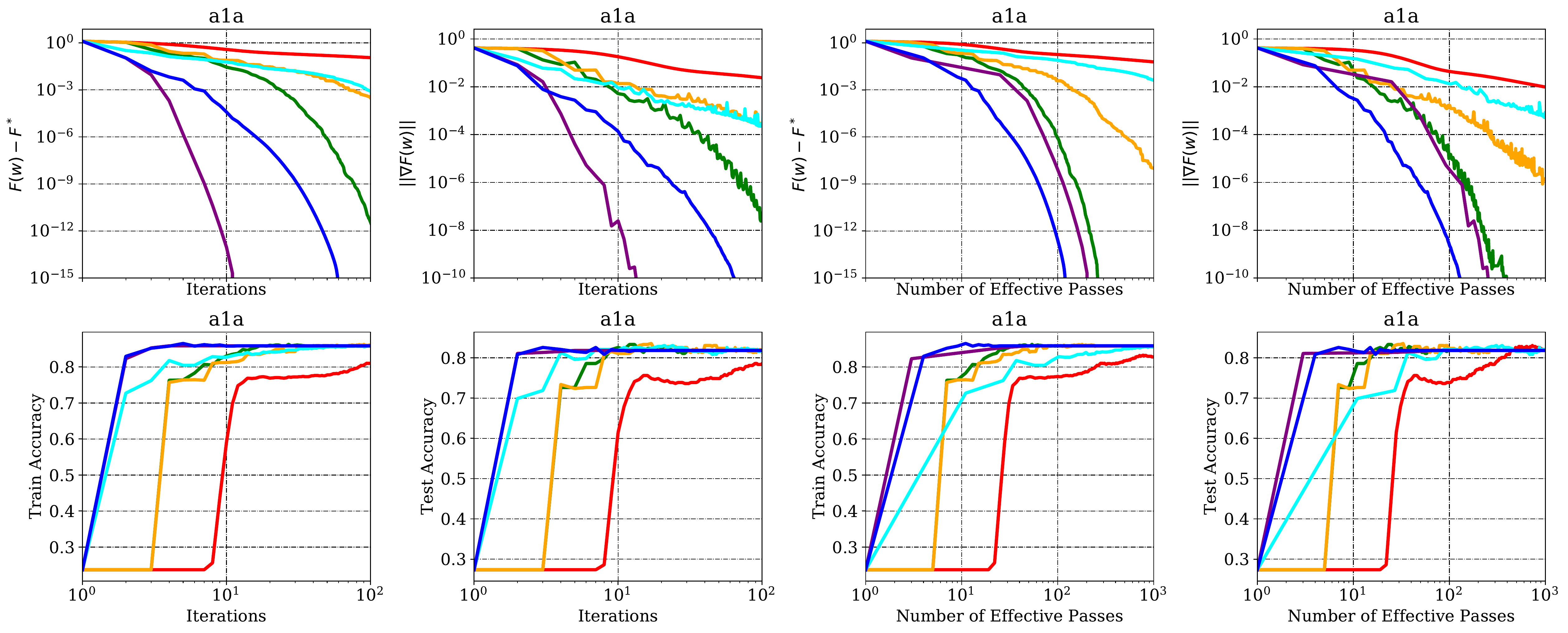}\vspace{-7pt}
	\caption{\texttt{a1a}: Deterministic Logistic Regression with $\lambda = 10^{-5}$.}
	\label{fig:a1a_5}
\end{figure*}
\vspace{-10pt}
\begin{figure*}[!h]
	\centering
	\includegraphics[width=0.8\textwidth]{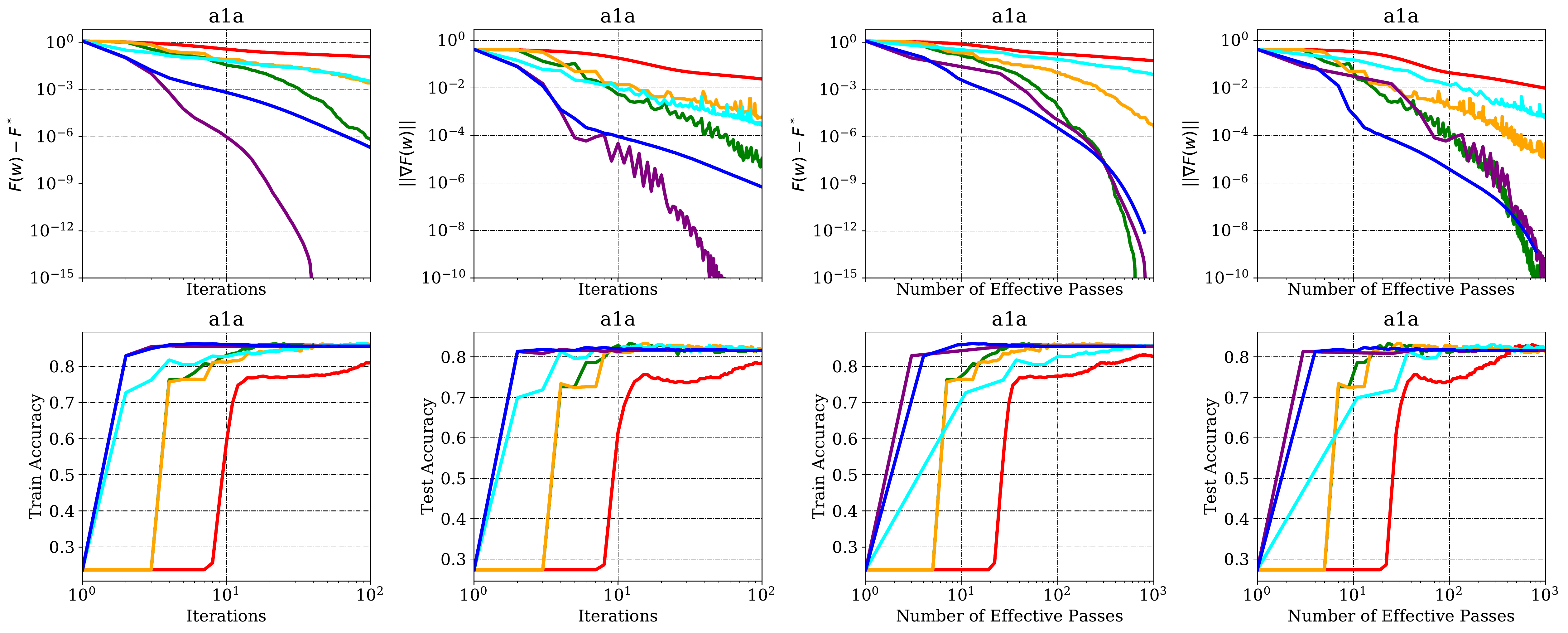}\vspace{-7pt}
	\caption{\texttt{a1a}: Deterministic Logistic Regression with $\lambda = 10^{-6}$.}
	\label{fig:a1a_6}
\end{figure*}

\begin{figure*}[ht]
	\centering
	{	\includegraphics[trim=10 110 10 110,clip, width=0.6\textwidth]{Figures/LIBSVM_LEGEND.pdf}
}
	\includegraphics[width=0.8\textwidth]{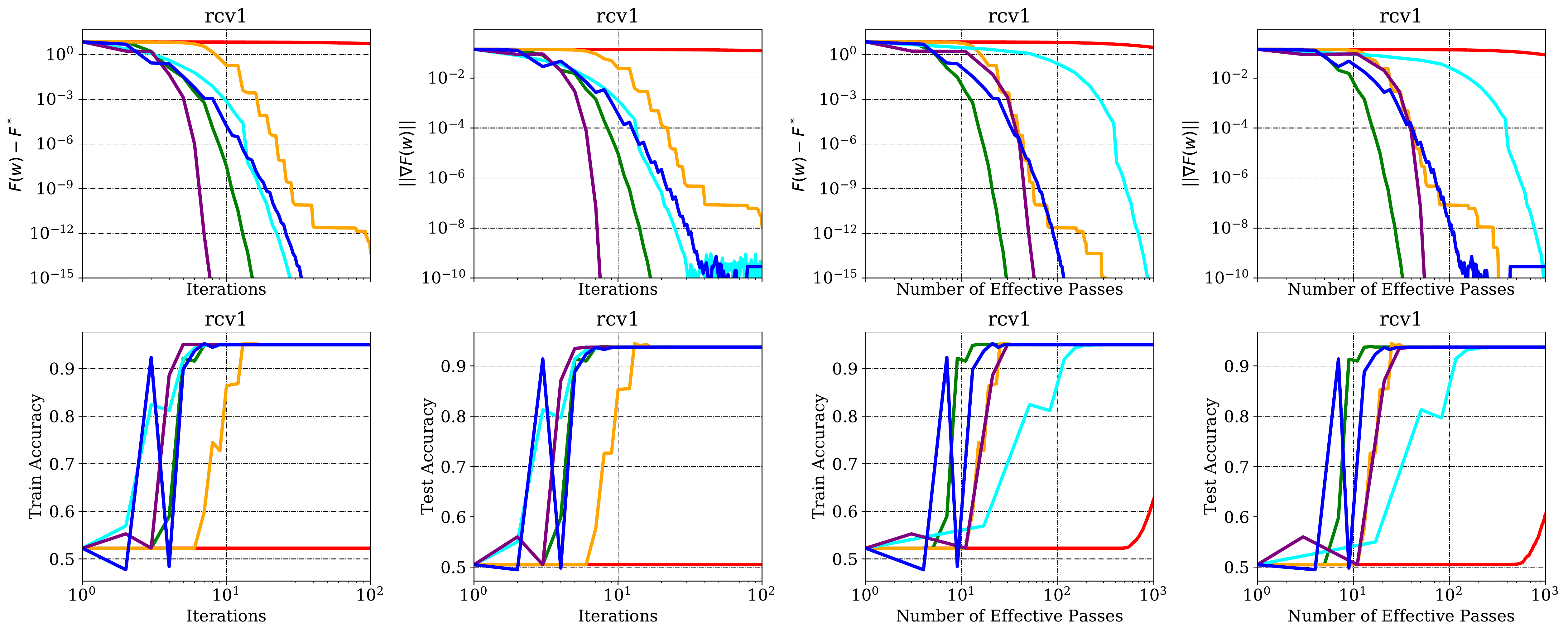}
	\caption{\texttt{rcv1}: Deterministic Logistic Regression with $\lambda = 10^{-3}$.}
	\label{fig:rcv1_train_3}
\end{figure*}
\vspace{-10pt}
\begin{figure*}[ht]
	\centering
	\includegraphics[width=0.8\textwidth]{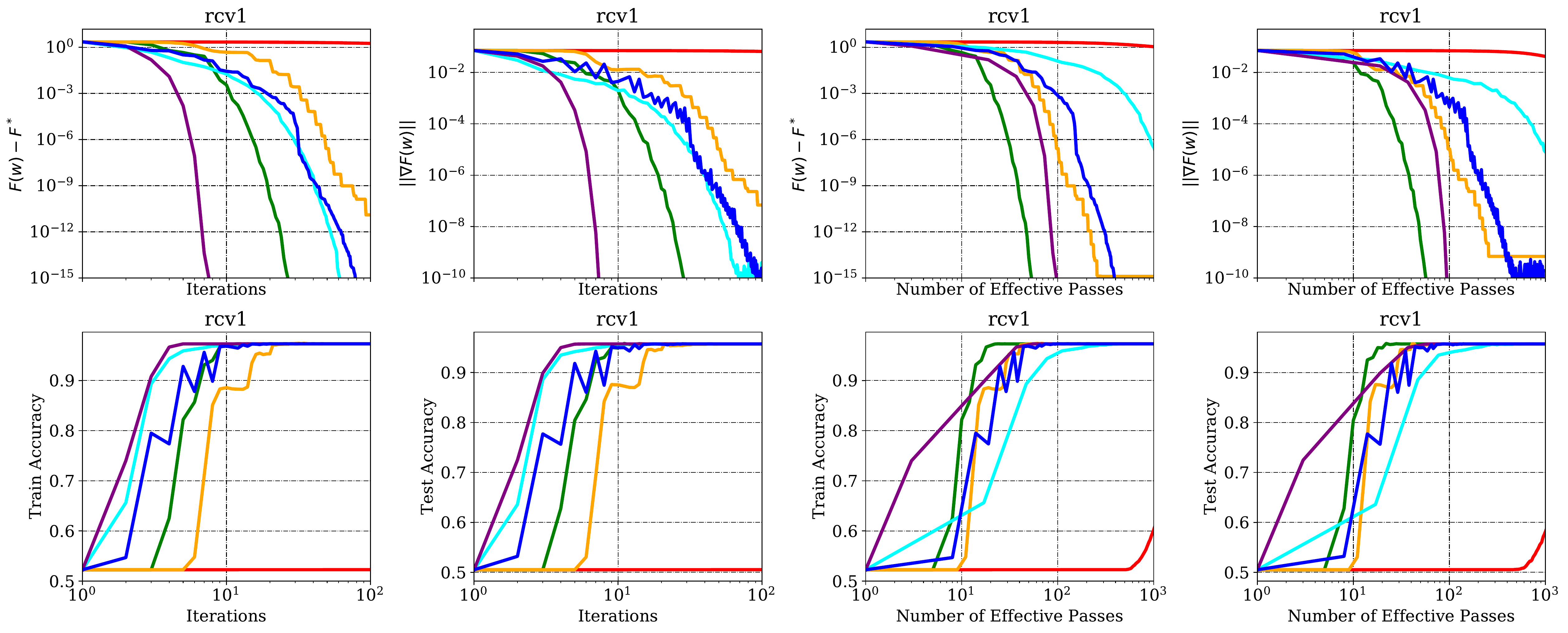}
	\caption{\texttt{rcv1}: Deterministic Logistic Regression with $\lambda = 10^{-4}$.}
\end{figure*}
\vspace{-10pt}
\begin{figure*}[ht]
	\centering
	\includegraphics[width=0.8\textwidth]{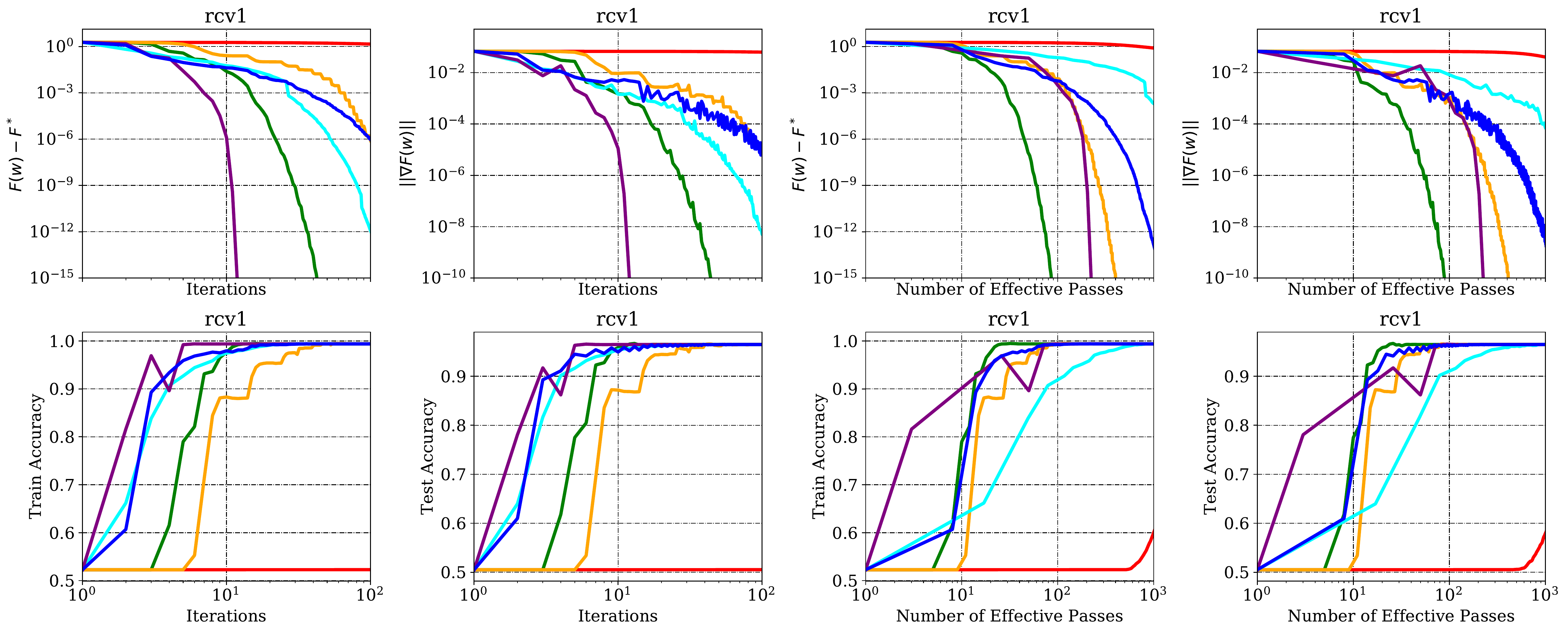}
	\caption{\texttt{rcv1}: Deterministic Logistic Regression with $\lambda = 10^{-5}$.}
\end{figure*}
\vspace{-10pt}
\begin{figure*}[htb]
	\centering
	\includegraphics[width=0.8\textwidth]{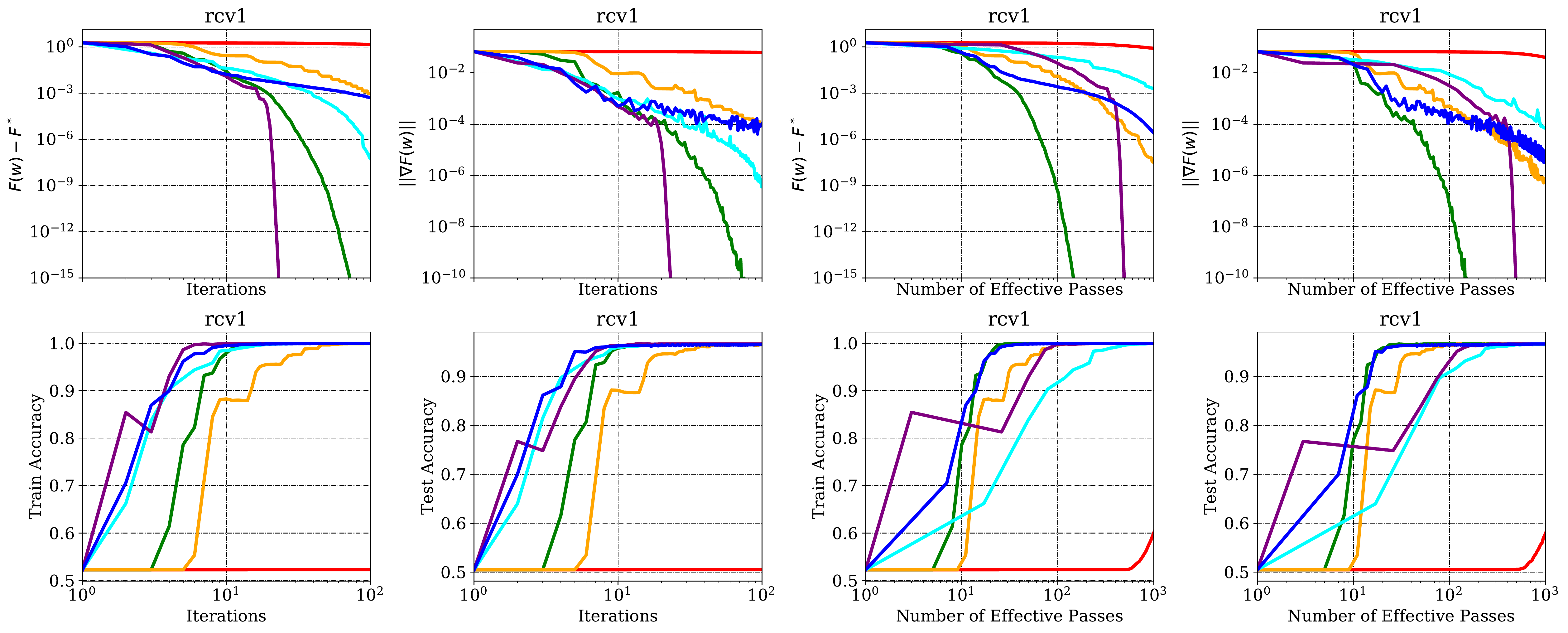}
	\caption{\texttt{rcv1}: Deterministic Logistic Regression with $\lambda = 10^{-6}$.}
	\label{fig:rcv1_train_6}
\end{figure*}

\begin{figure*}[ht]
	\centering
	{	\includegraphics[trim=10 110 10 110,clip, width=0.6\textwidth]{Figures/LIBSVM_LEGEND.pdf}
}
	\includegraphics[width=0.8\textwidth]{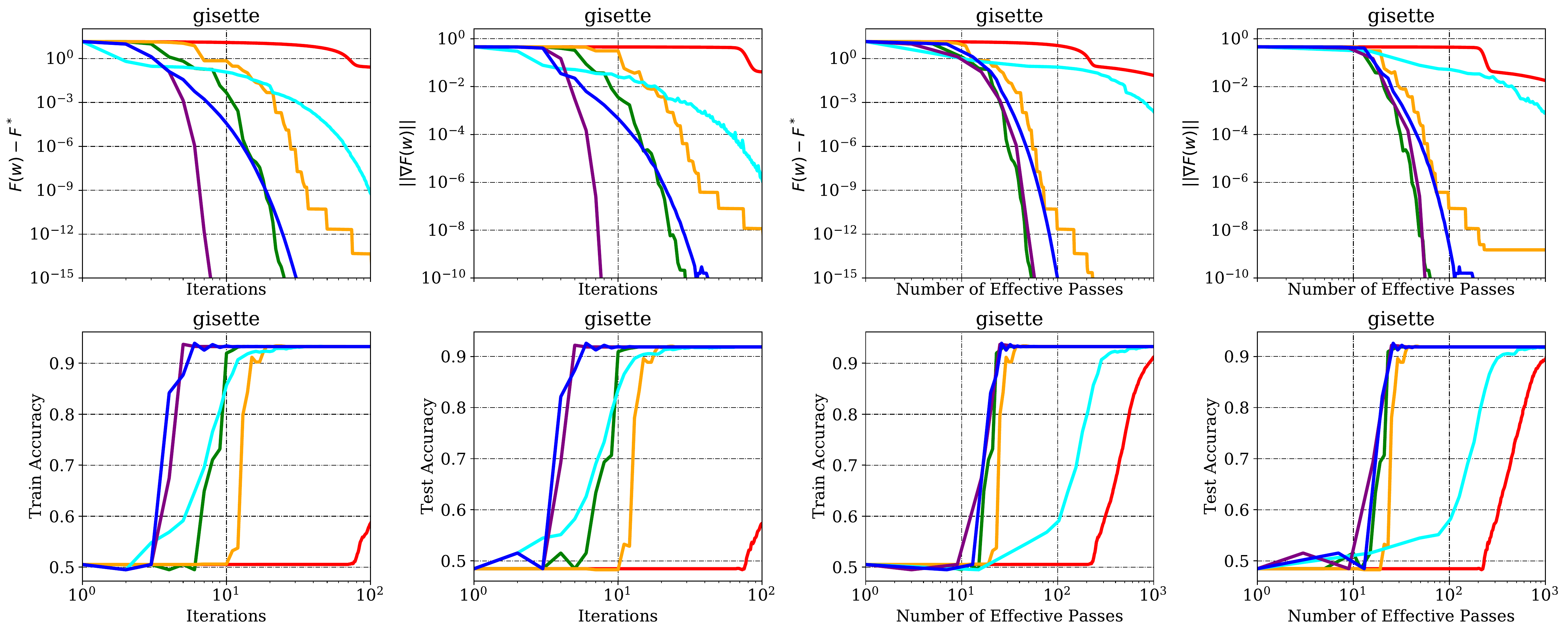}
    \caption{\texttt{gisette}: Deterministic Logistic Regression with $\lambda = 10^{-3}$.}
	\label{fig:gisette_scale_3}
\end{figure*}
\vspace{-10pt}
\begin{figure*}[ht]
	\centering
	\includegraphics[width=0.8\textwidth]{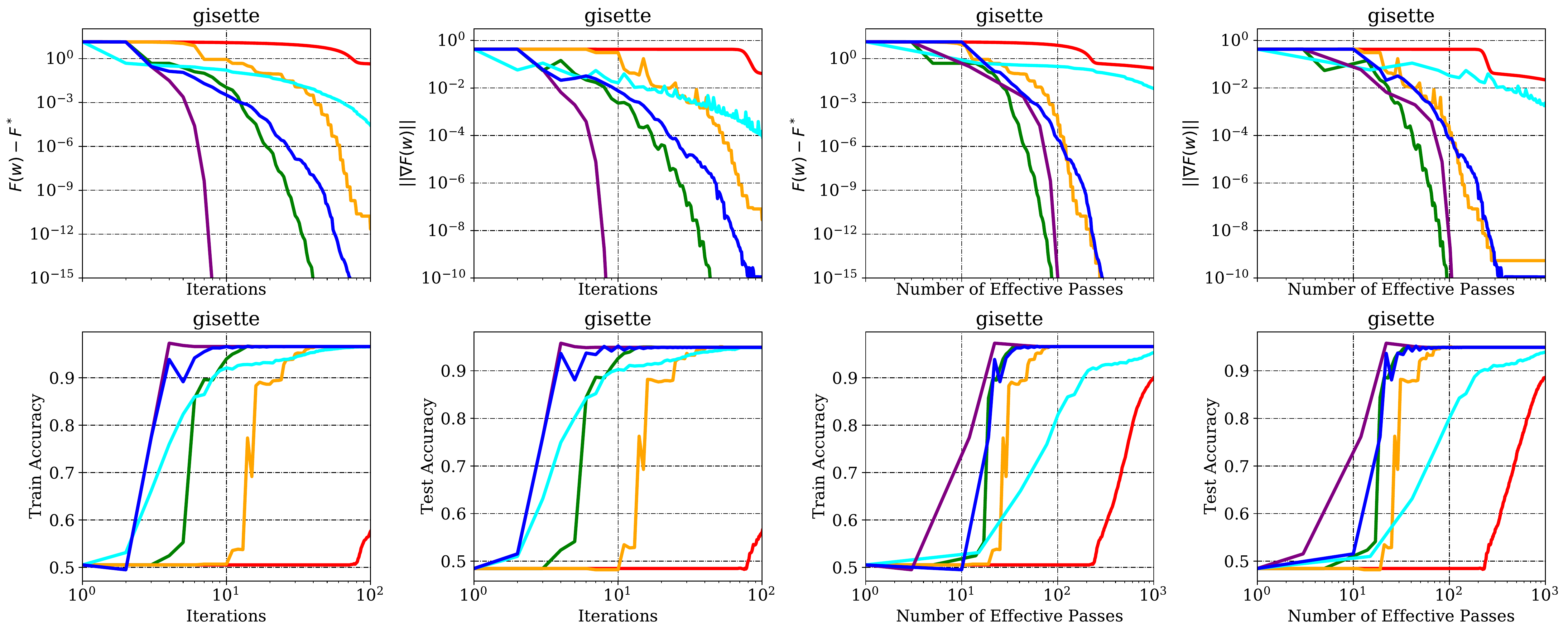}
	\caption{\texttt{gisette}: Deterministic Logistic Regression with $\lambda = 10^{-4}$.}
	\label{fig:gisette_scale_4}
\end{figure*}
\vspace{-10pt}
\begin{figure*}[ht]
	\centering
	\includegraphics[width=0.8\textwidth]{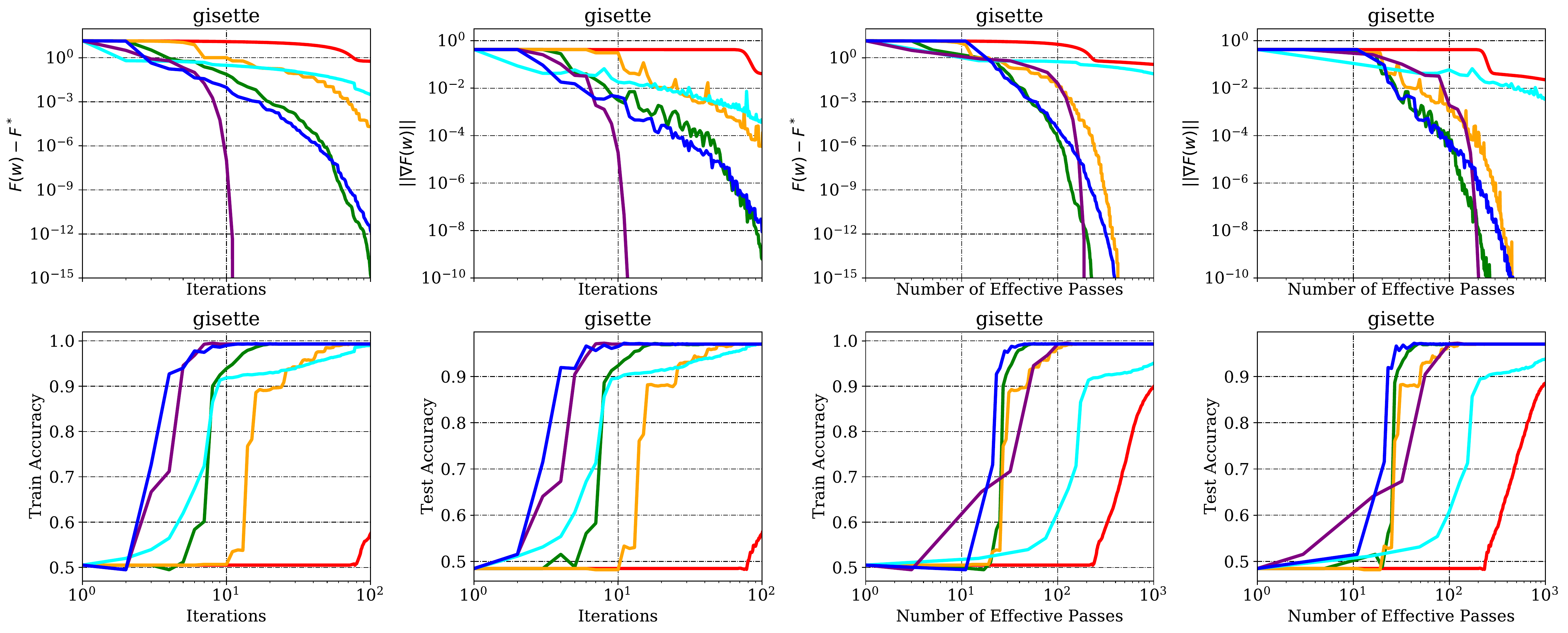}
	\caption{\texttt{gisette}: Deterministic Logistic Regression with $\lambda = 10^{-5}$.}
	\label{fig:gisette_scale_5}
\end{figure*}
\vspace{-10pt}
\begin{figure*}[ht]
	\centering
	\includegraphics[width=0.8\textwidth]{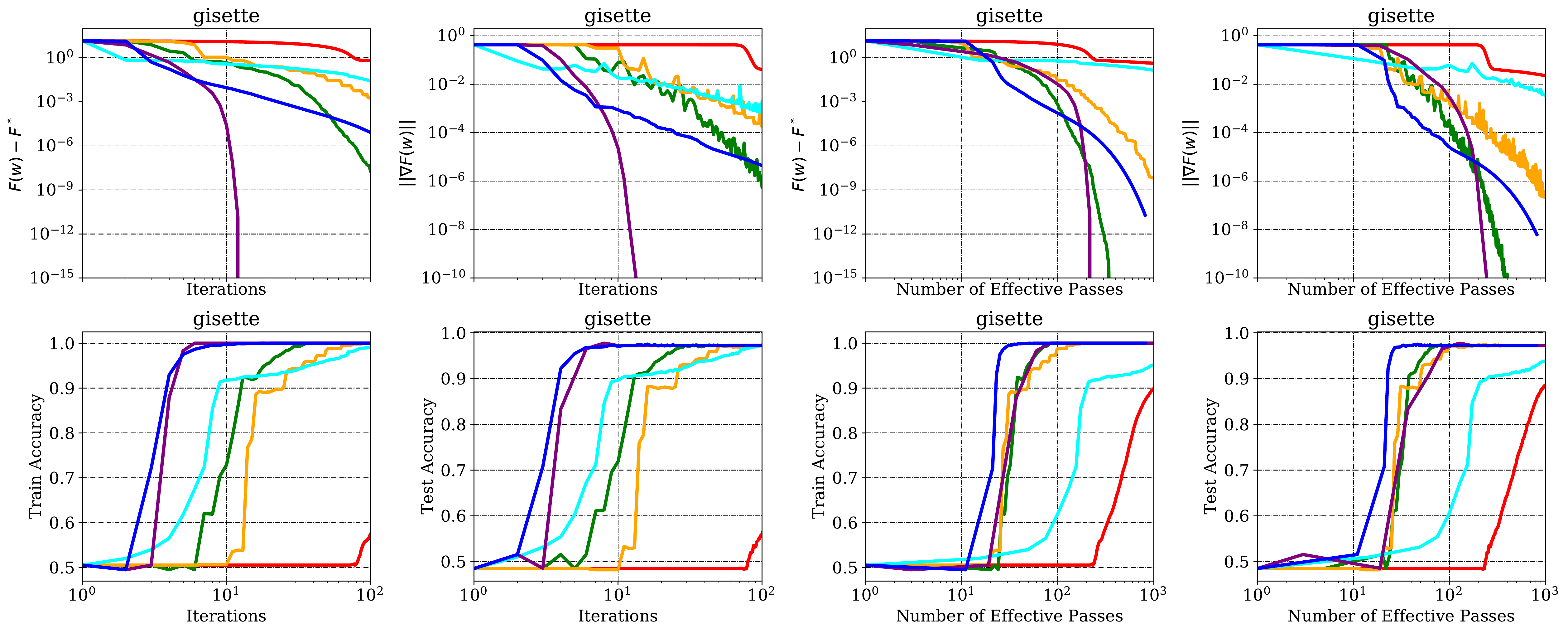}
	\caption{\texttt{gisette}: Deterministic Logistic Regression with $\lambda = 10^{-6}$.}
	\label{fig:gisette_scale_6}
\end{figure*}

\begin{figure*}[ht]
	\centering
	{	\includegraphics[trim=10 110 10 110,clip, width=0.6\textwidth]{Figures/LIBSVM_LEGEND.pdf}
}
	\includegraphics[width=0.8\textwidth]{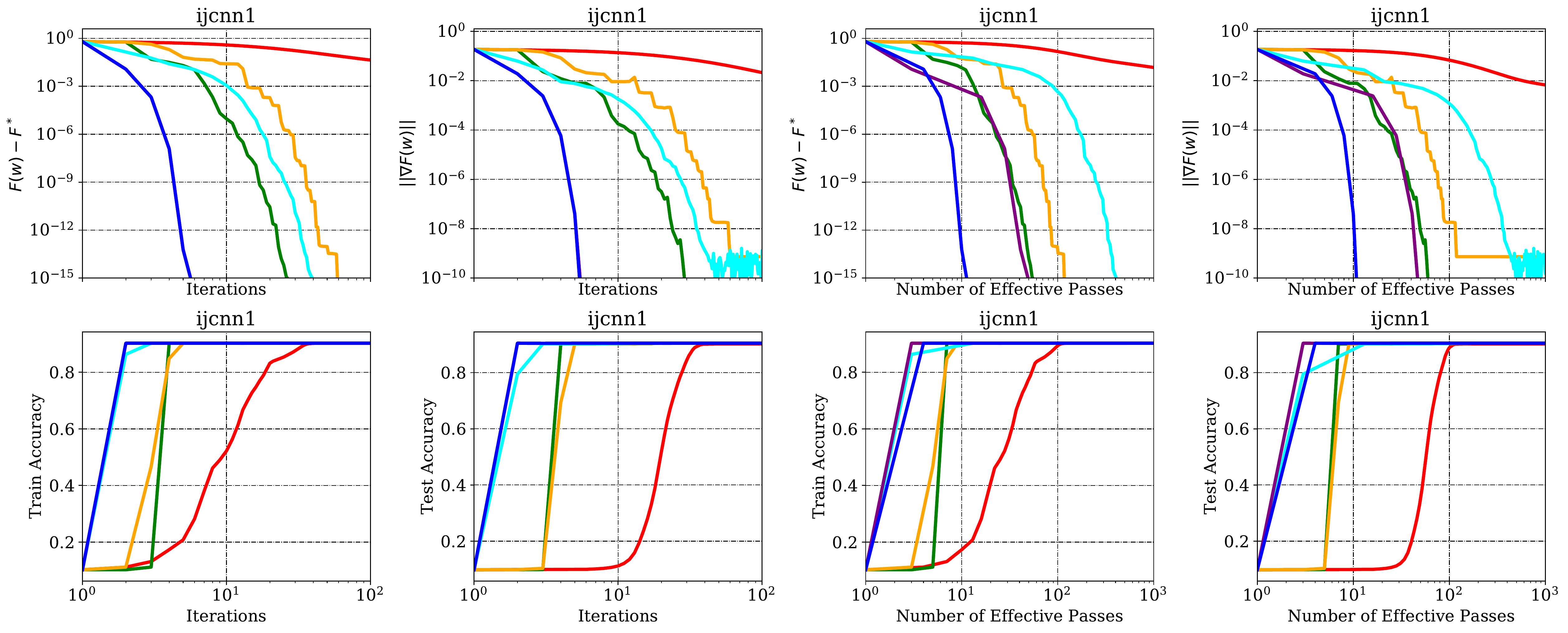}
	\caption{\texttt{ijcnn1}: Deterministic Logistic Regression with $\lambda = 10^{-3}$.}
	\label{fig:ijcnn1_3}
\end{figure*}
\vspace{-10pt}
\begin{figure*}[ht]
	\centering
	\includegraphics[width=0.8\textwidth]{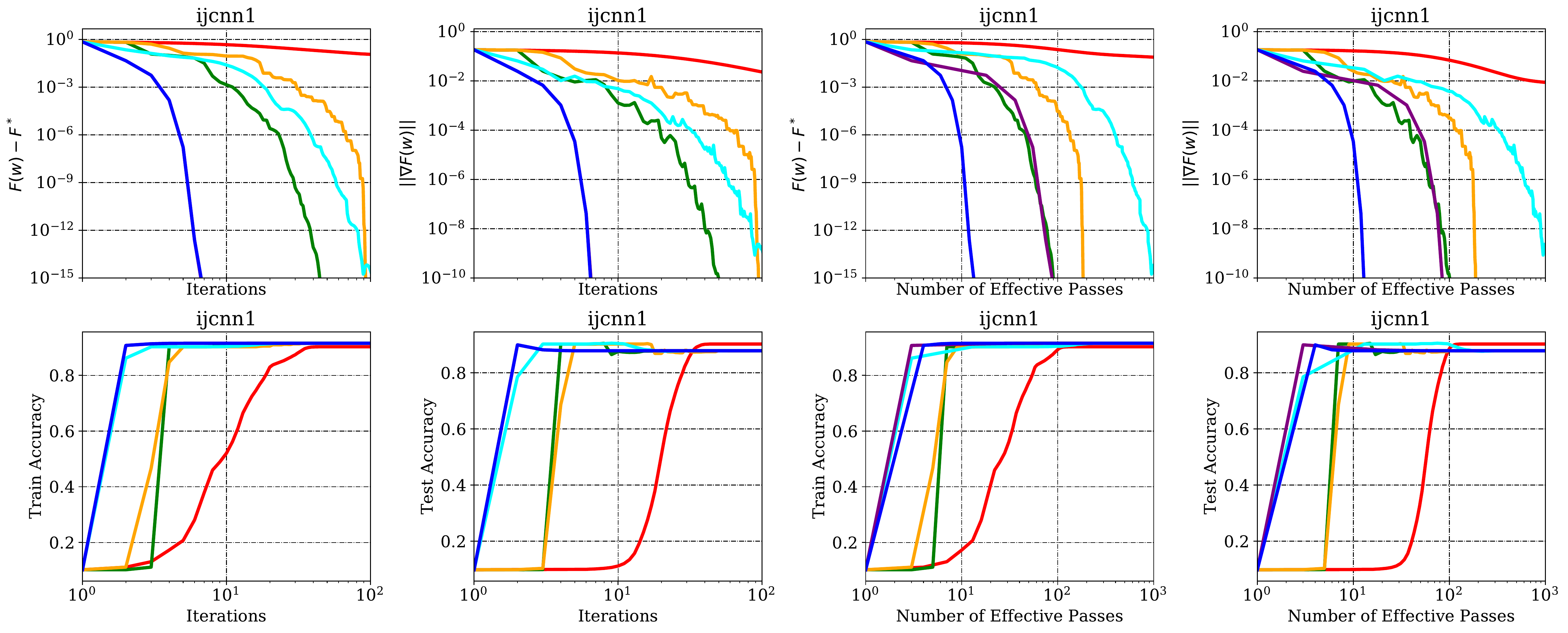}
	\caption{\texttt{ijcnn1}: Deterministic Logistic Regression with $\lambda = 10^{-4}$.}
	\label{fig:ijcnn1_4}
\end{figure*}
\vspace{-10pt}
\begin{figure*}[ht]
	\centering
	\includegraphics[width=0.8\textwidth]{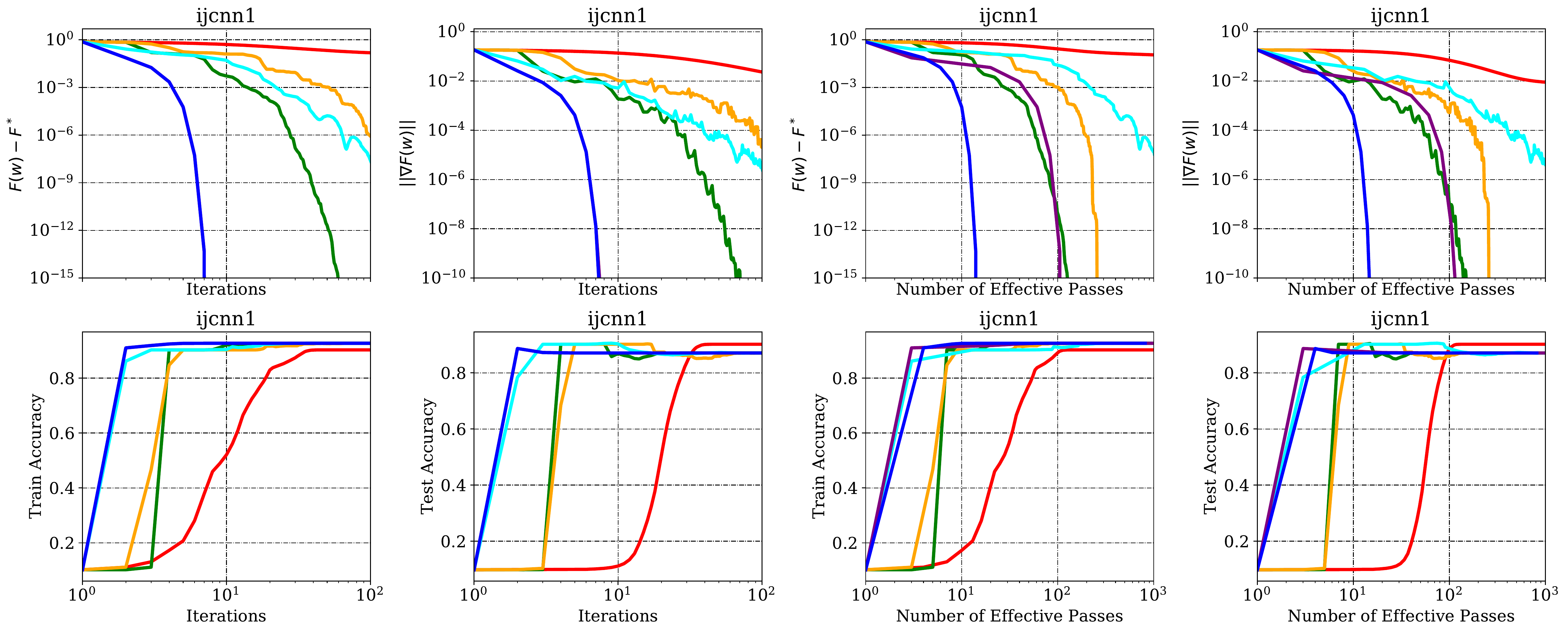}
	\caption{\texttt{ijcnn1}: Deterministic Logistic Regression with $\lambda = 10^{-5}$.}
	\label{fig:gisette_scale_5}
\end{figure*}
\vspace{-10pt}
\begin{figure*}[ht]
	\centering
	\includegraphics[width=0.8\textwidth]{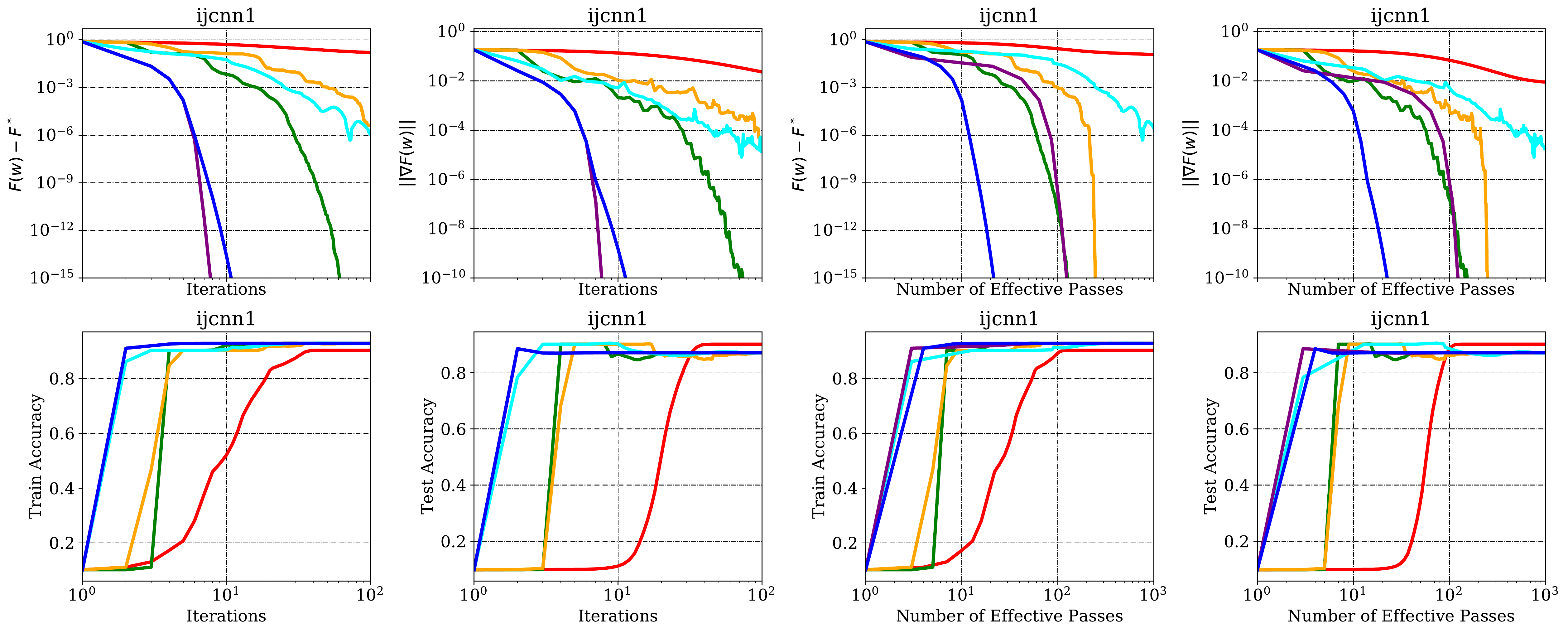}
	\caption{\texttt{ijcnn1}: Deterministic Logistic Regression with $\lambda = 10^{-6}$.}
	\label{fig:ijcnn1_6}
\end{figure*}

\clearpage

\subsubsection{Sensitivity of \SONIA{} to memory hyper-parameter}\label{sec:sensMMR}
\vspace{-10pt}
\begin{figure*}[ht]
	\centering
	\includegraphics[width=0.8\textwidth]{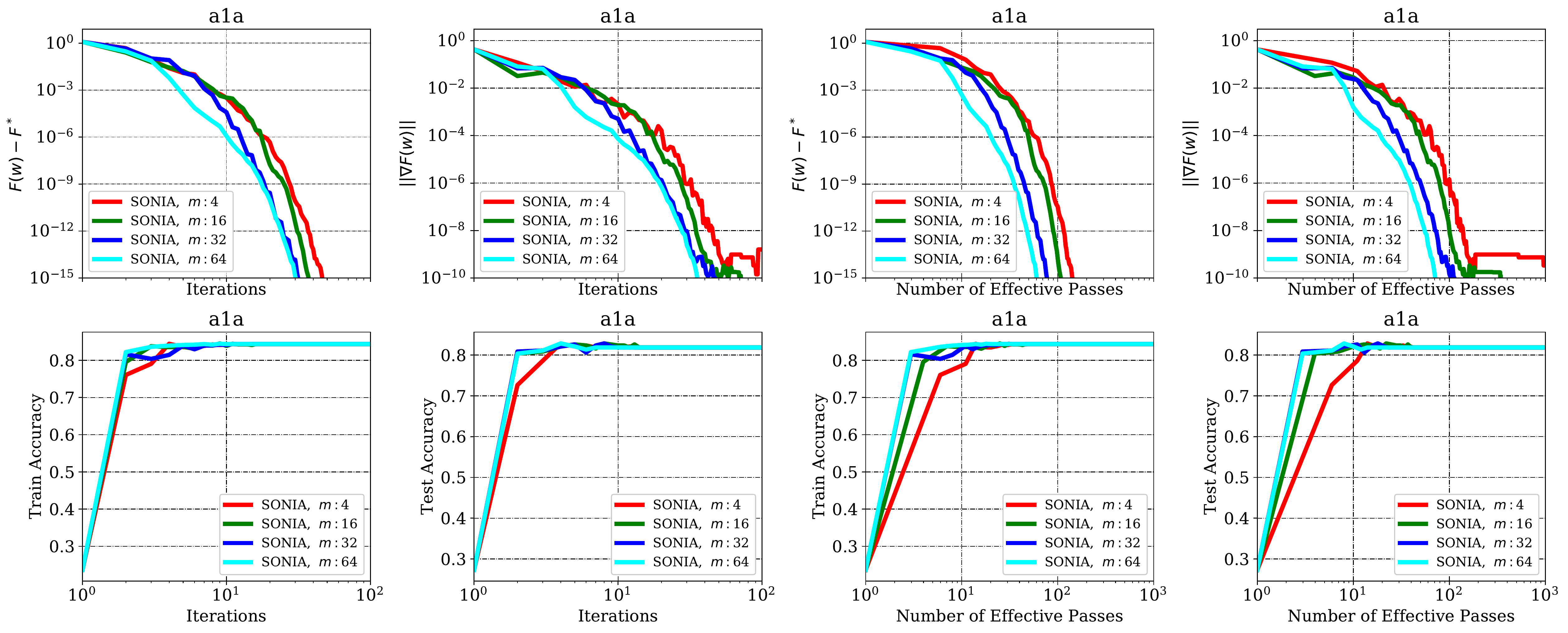}\vspace{-7pt}
	\caption{\texttt{a1a}: Sensitivity of \texttt{SONIA} to memory,  Deterministic Logistic Regression ($\lambda = 10^{-3}$).}
	\label{fig:a1a_3_FastSR1_sensitivityAnlysis_MMR}
\end{figure*}
\vspace{-10pt}
\begin{figure*}[ht]
	\centering
	\includegraphics[width=0.8\textwidth]{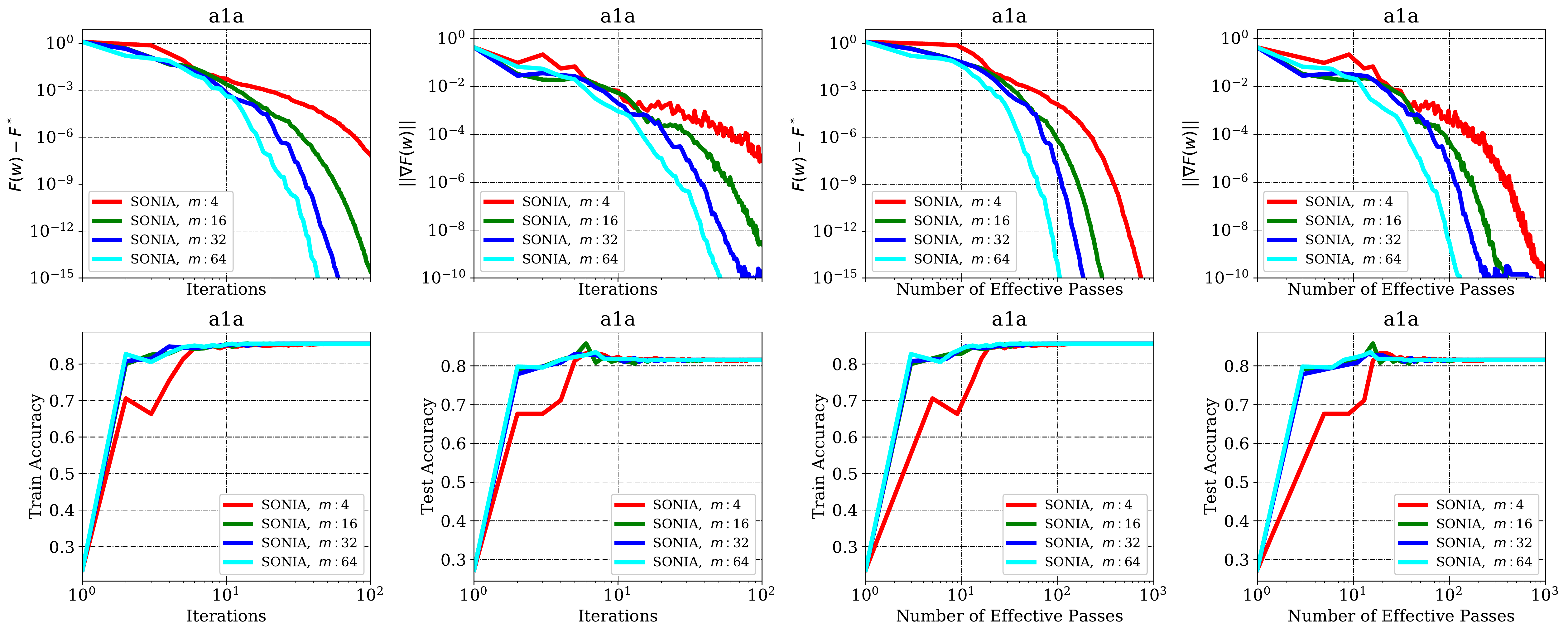}\vspace{-7pt}
	\caption{\texttt{a1a}: Sensitivity of \texttt{SONIA} to memory,  Deterministic Logistic Regression ($\lambda = 10^{-4}$).}
	\label{fig:a1a_4_FastSR1_sensitivityAnlysis_MMR}
\end{figure*}
\vspace{-10pt}
\begin{figure*}[ht]
	\centering
	\includegraphics[width=0.8\textwidth]{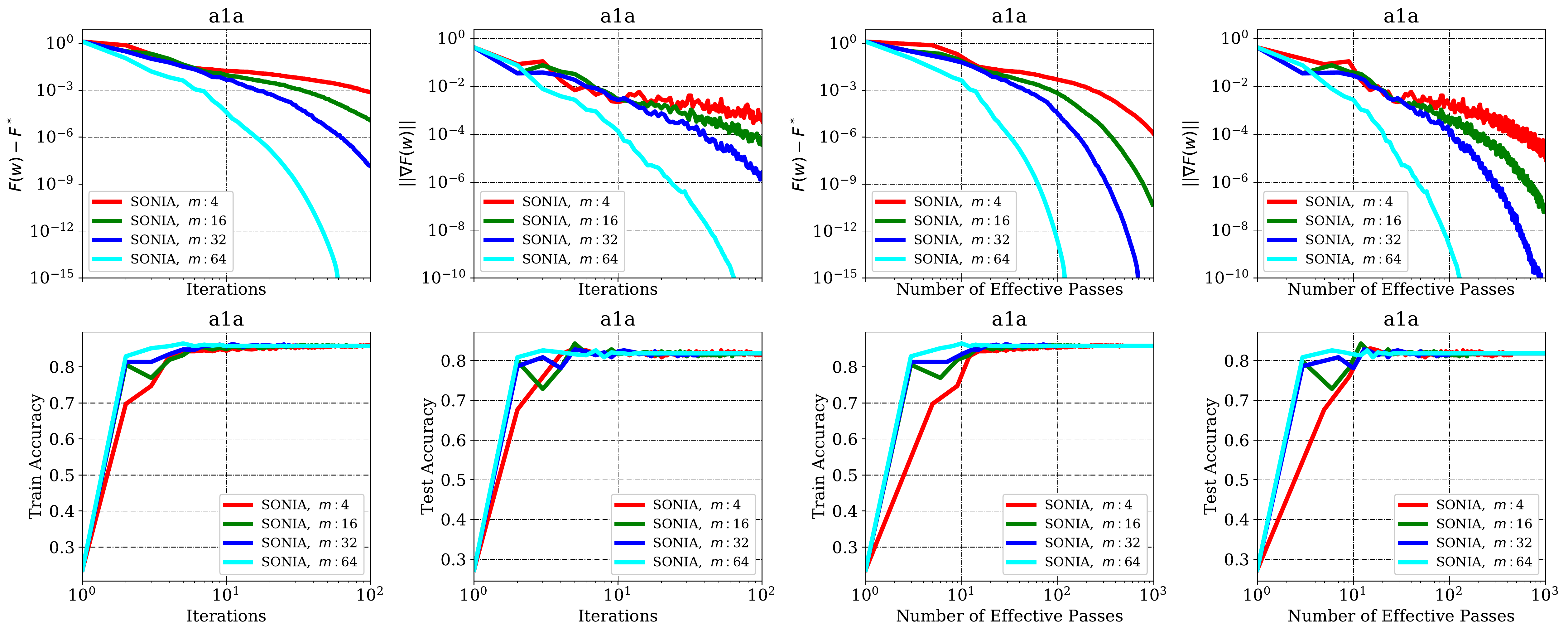}\vspace{-7pt}
	\caption{\texttt{a1a}: Sensitivity of \texttt{SONIA} to memory,  Deterministic Logistic Regression ($\lambda = 10^{-5}$).}
	\label{fig:a1a_5_FastSR1_sensitivityAnlysis_MMR}
\end{figure*}
\vspace{-10pt}
\begin{figure*}[!h]
	\centering
	\includegraphics[width=0.8\textwidth]{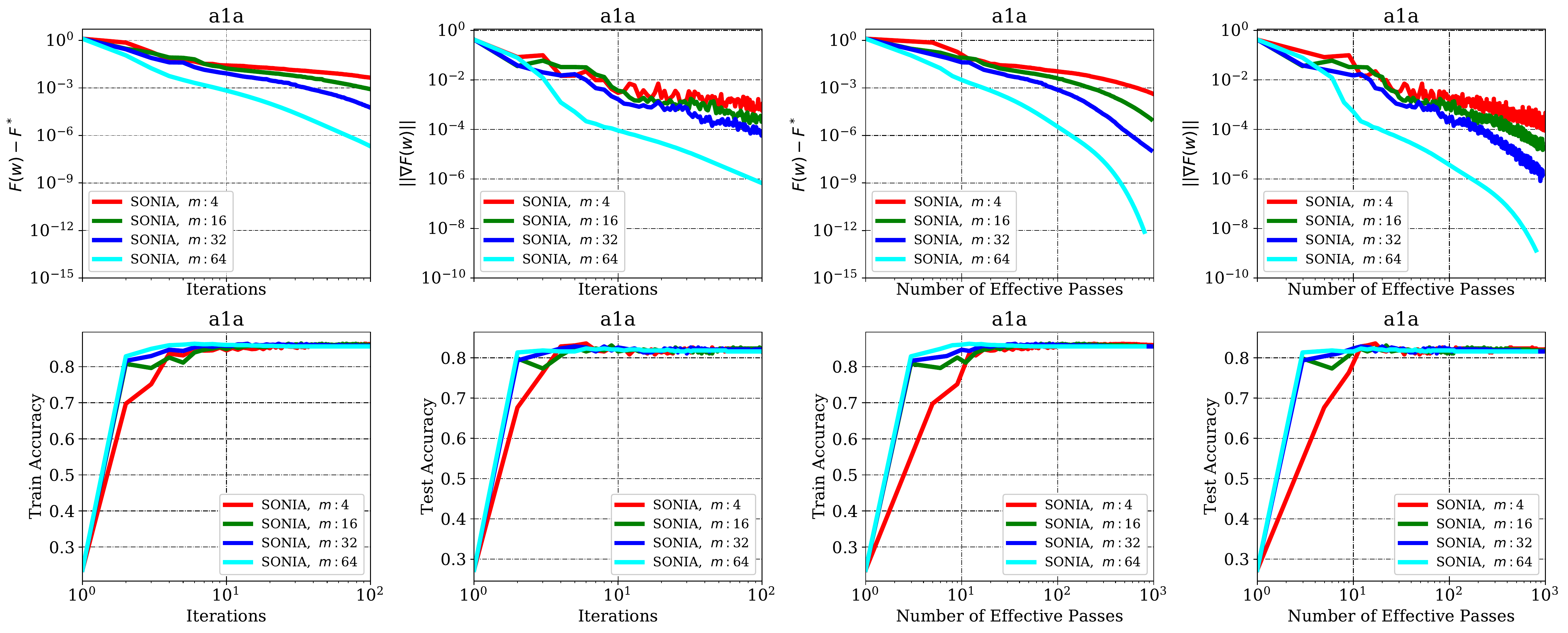}\vspace{-7pt}
	\caption{\texttt{a1a}: Sensitivity of \texttt{SONIA} to memory,  Deterministic Logistic Regression ($\lambda = 10^{-6}$).}
	\label{fig:a1a_6_FastSR1_sensitivityAnlysis_MMR}
\end{figure*}

\vspace{-10pt}
\begin{figure*}[ht]
	\centering
	\includegraphics[width=0.8\textwidth]{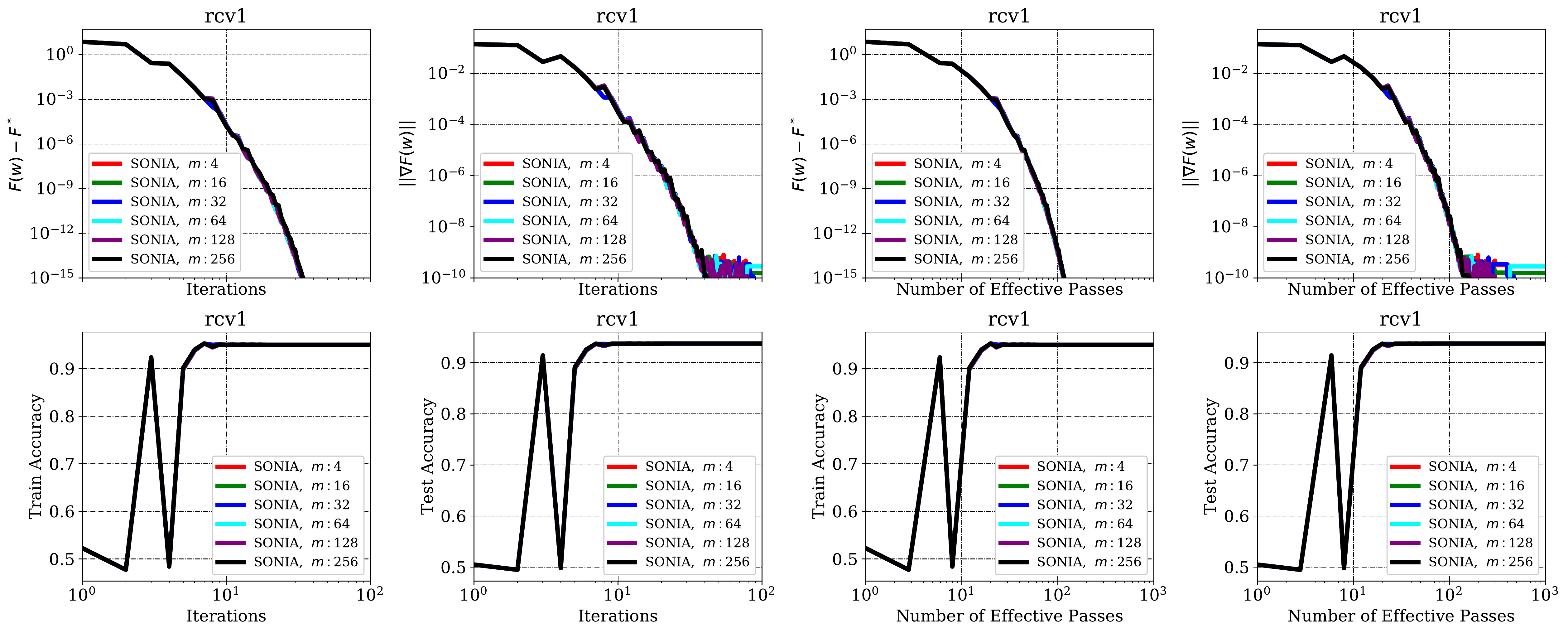}
	\caption{\texttt{rcv1}: Sensitivity of \texttt{SONIA} to memory,  Deterministic Logistic Regression ($\lambda = 10^{-3}$).}
	\label{fig:rcv1_train_3_FastSR1_sensitivityAnlysis_MMR}
\end{figure*}
\vspace{-10pt}
\begin{figure*}[ht]
	\centering
	\includegraphics[width=0.8\textwidth]{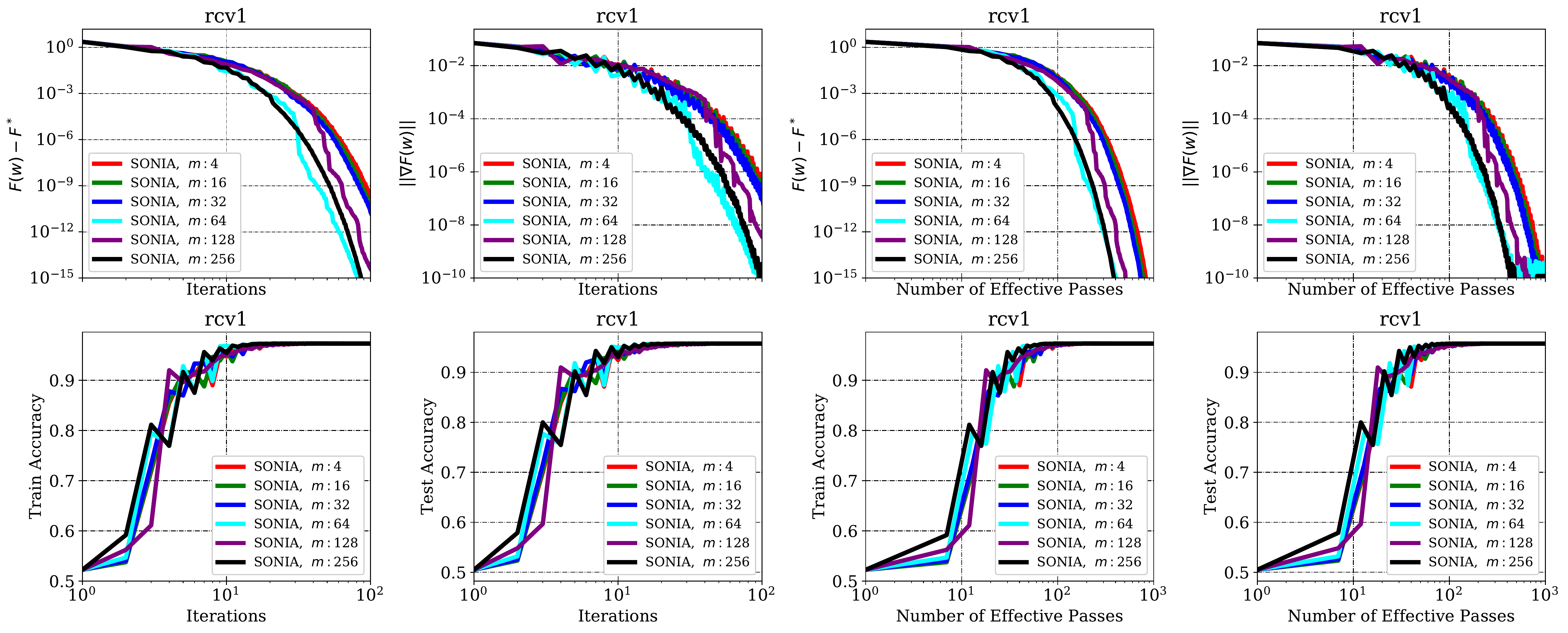}
	\caption{\texttt{rcv1}: Sensitivity of \texttt{SONIA} to memory,  Deterministic Logistic Regression ($\lambda = 10^{-4}$).}
	\label{fig:rcv1_train_4_FastSR1_sensitivityAnlysis_MMR}
\end{figure*}
\vspace{-10pt}
\begin{figure*}[ht]
	\centering
	\includegraphics[width=0.8\textwidth]{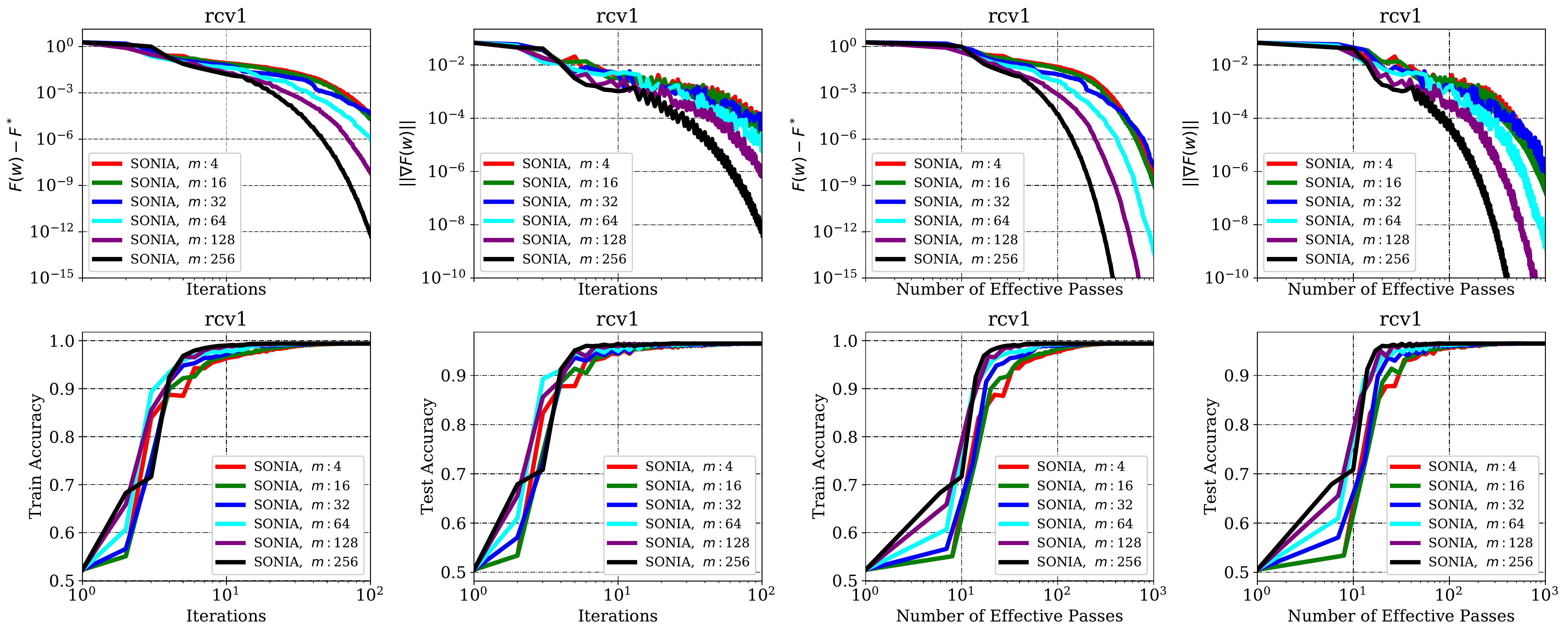}
	\caption{\texttt{rcv1}: Sensitivity of \texttt{SONIA} to memory,  Deterministic Logistic Regression ($\lambda = 10^{-5}$).}
	\label{fig:rcv1_train_5_FastSR1_sensitivityAnlysis_MMR}
\end{figure*}
\vspace{-10pt}
\begin{figure*}[!h]
	\centering
	\includegraphics[width=0.8\textwidth]{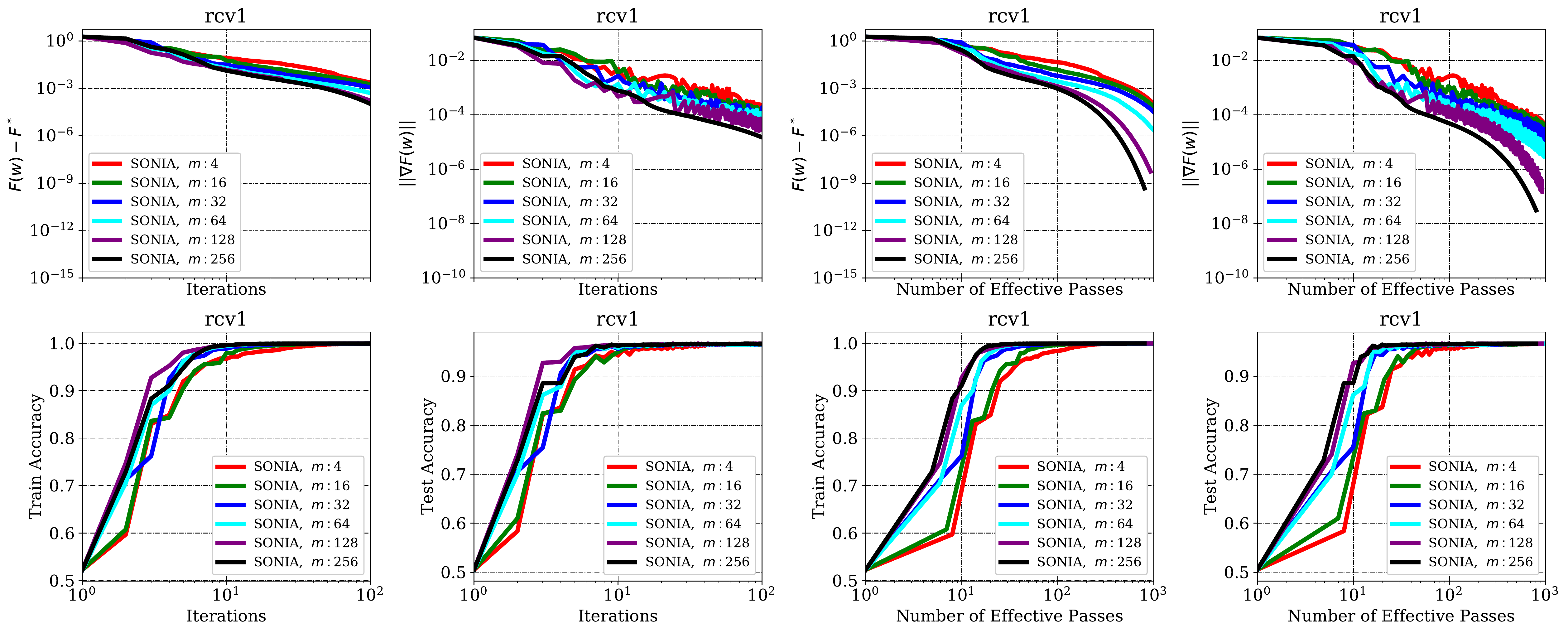}
	\caption{\texttt{rcv1}: Sensitivity of \texttt{SONIA} to memory,  Deterministic Logistic Regression ($\lambda = 10^{-6}$).}
	\label{fig:rcv1_train_6_FastSR1_sensitivityAnlysis_MMR}
\end{figure*}

\vspace{-10pt}
\begin{figure*}[ht]
	\centering
	\includegraphics[width=0.8\textwidth]{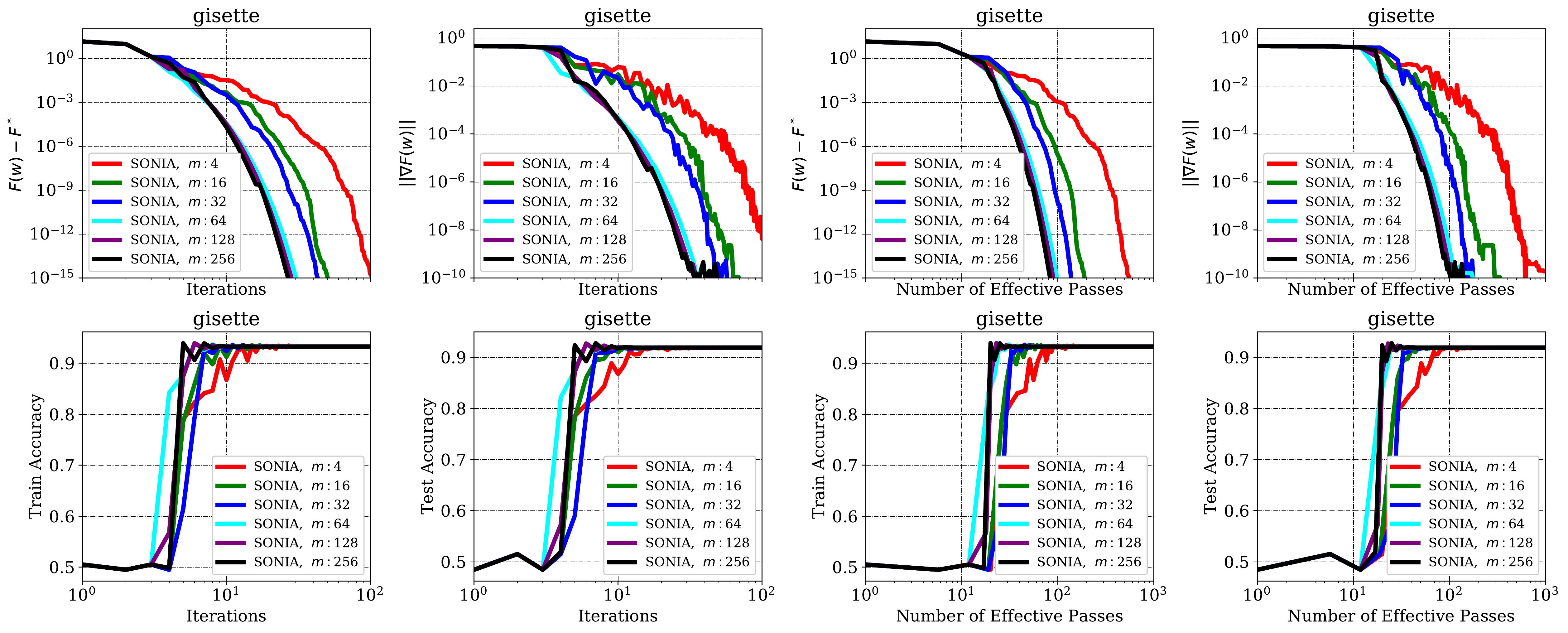}
	\caption{\texttt{gisette}: Sensitivity of \texttt{SONIA} to memory,  Deterministic Logistic Regression ($\lambda = 10^{-3}$).}
	\label{fig:gisette_scale_3_FastSR1_sensitivityAnlysis_MMR}
\end{figure*}
\vspace{-10pt}
\begin{figure*}[ht]
	\centering
	\includegraphics[width=0.8\textwidth]{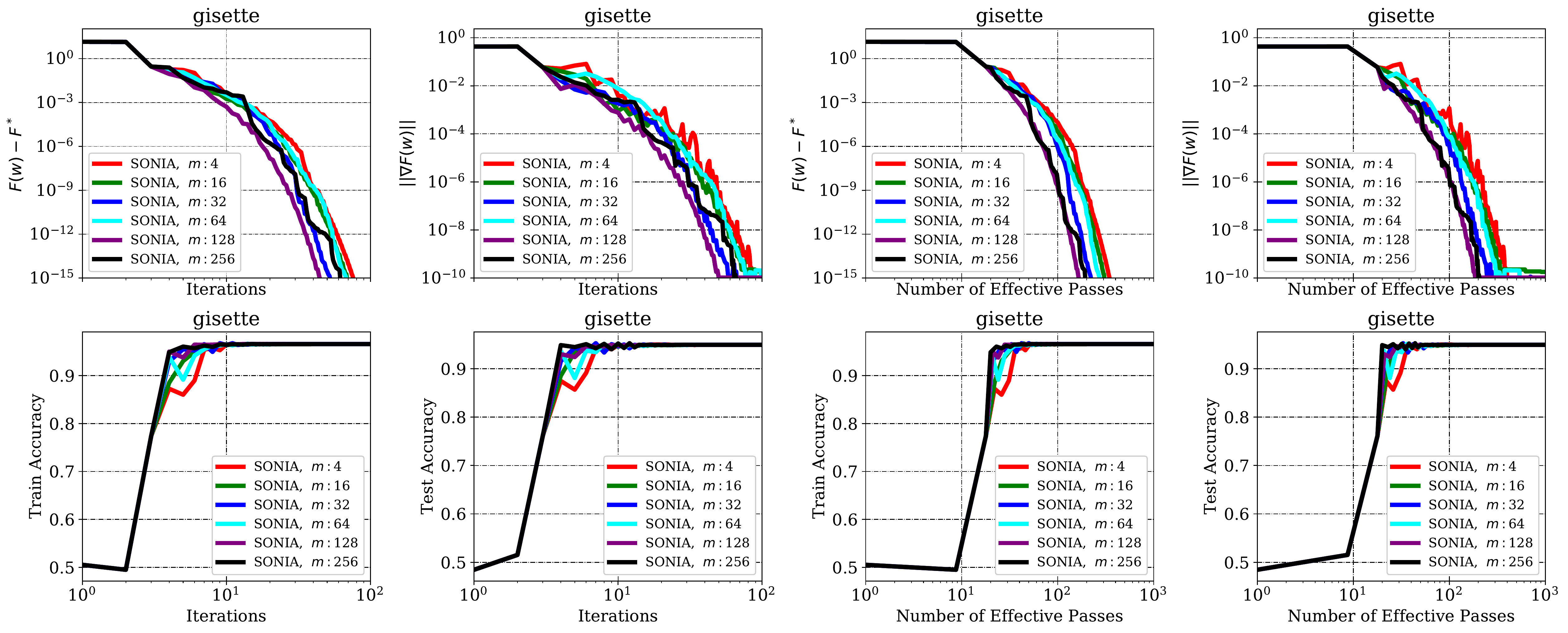}
	\caption{\texttt{gisette}: Sensitivity of \texttt{SONIA} to memory,  Deterministic Logistic Regression ($\lambda = 10^{-4}$).}
	\label{fig:gisette_scale_4_FastSR1_sensitivityAnlysis_MMR}
\end{figure*}
\vspace{-10pt}
\begin{figure*}[ht]
	\centering
	\includegraphics[width=0.8\textwidth]{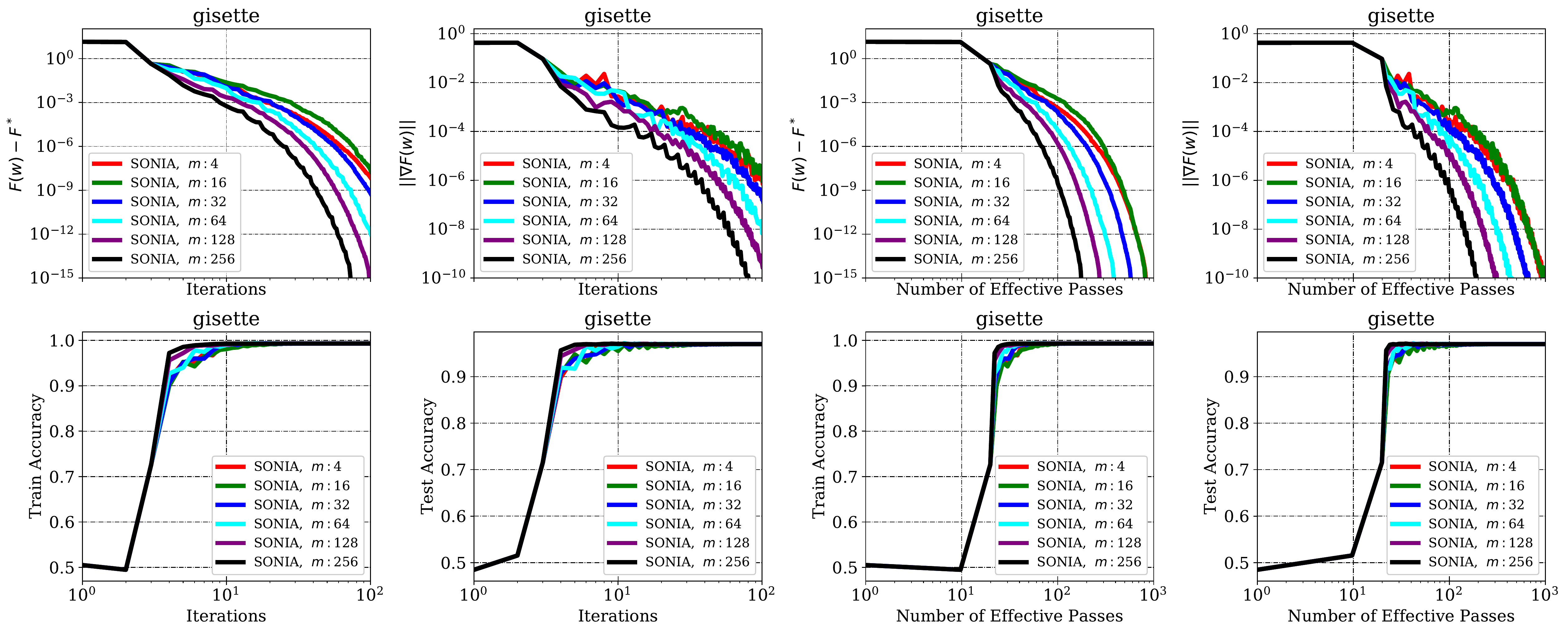}
	\caption{\texttt{gisette}: Sensitivity of \texttt{SONIA} to memory,  Deterministic Logistic Regression ($\lambda = 10^{-5}$).}
	\label{fig:gisettegisette_scale_5_FastSR1_sensitivityAnlysis_MMR}
\end{figure*}
\vspace{-10pt}
\begin{figure*}[htb]
	\centering
	\includegraphics[width=0.8\textwidth]{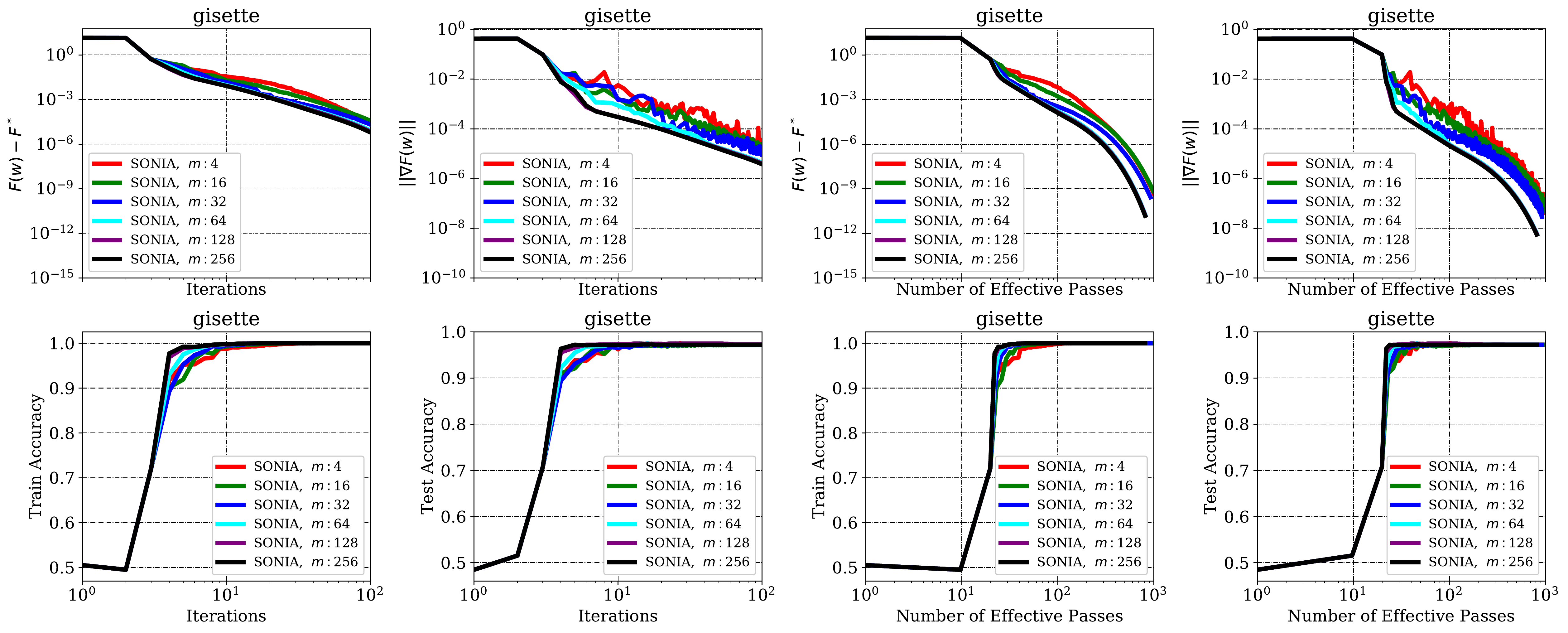}
	\caption{\texttt{gisette}: Sensitivity of \texttt{SONIA} to memory,  Deterministic Logistic Regression ($\lambda = 10^{-6}$).}
	\label{fig:gisette_scale_6_FastSR1_sensitivityAnlysis_MMR}
\end{figure*}
\vspace{-10pt}
\begin{figure*}[ht]
	\centering
	\includegraphics[width=0.8\textwidth]{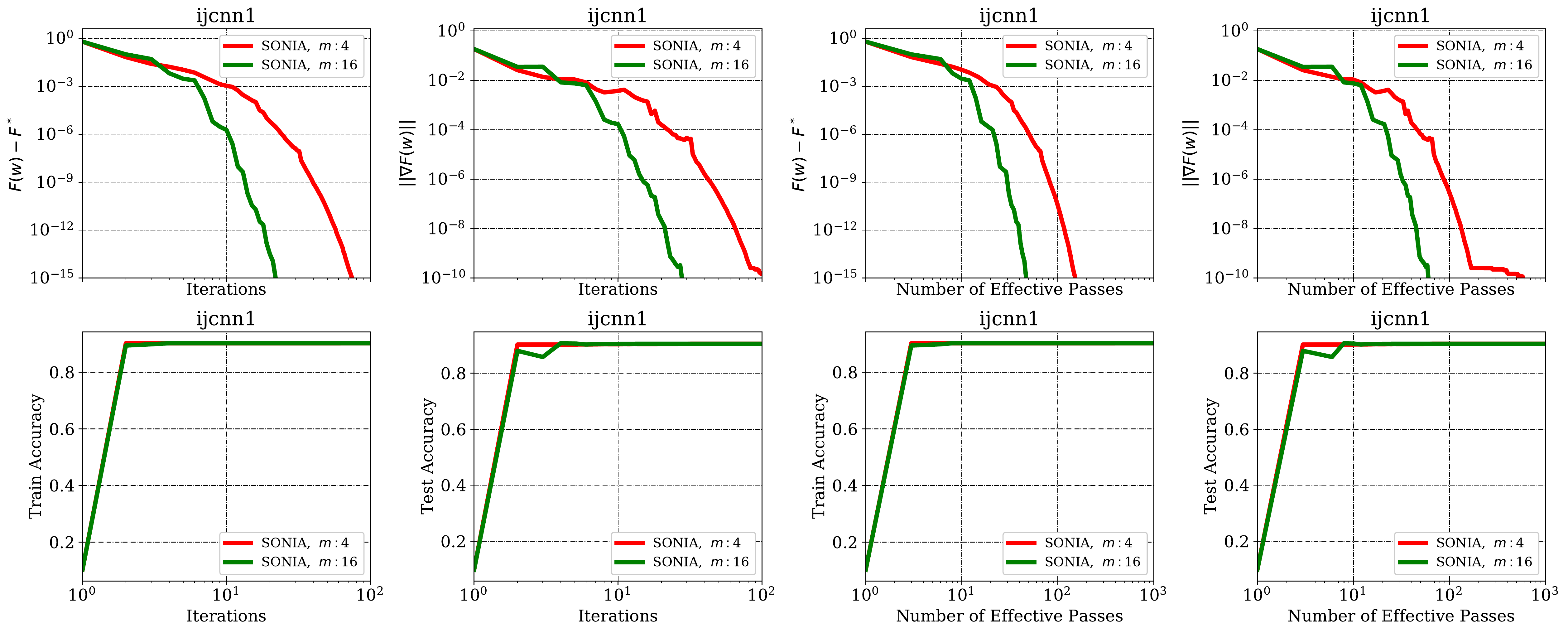}
	\caption{\texttt{ijcnn1}: Sensitivity of \texttt{SONIA} to memory,  Deterministic Logistic Regression ($\lambda = 10^{-3}$).}
	\label{fig:ijcnn1_3_FastSR1_sensitivityAnlysis_MMR}
\end{figure*}
\vspace{-10pt}
\begin{figure*}[ht]
	\centering
	\includegraphics[width=0.8\textwidth]{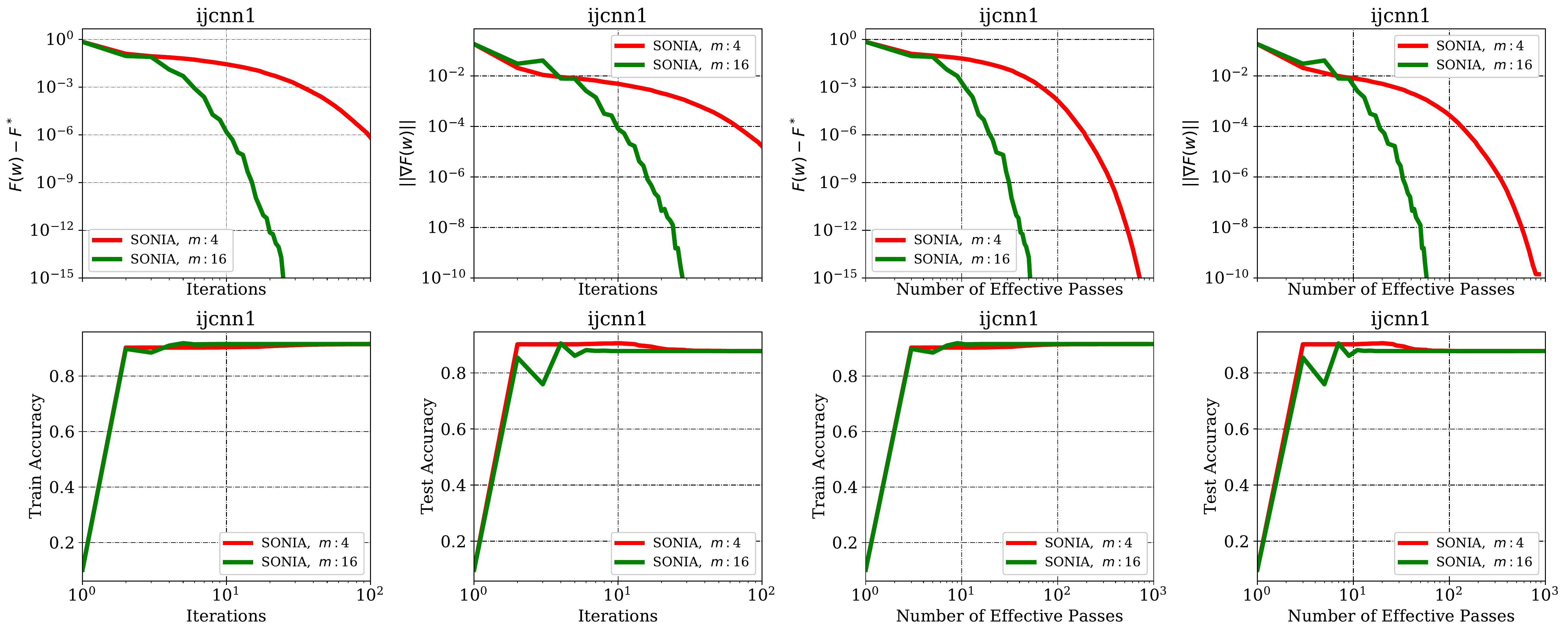}
	\caption{\texttt{ijcnn1}: Sensitivity of \texttt{SONIA} to memory,  Deterministic Logistic Regression ($\lambda = 10^{-4}$).}
	\label{fig:ijcnn1_4_FastSR1_sensitivityAnlysis_MMR}
\end{figure*}
\vspace{-10pt}
\begin{figure*}[ht]
	\centering
	\includegraphics[width=0.8\textwidth]{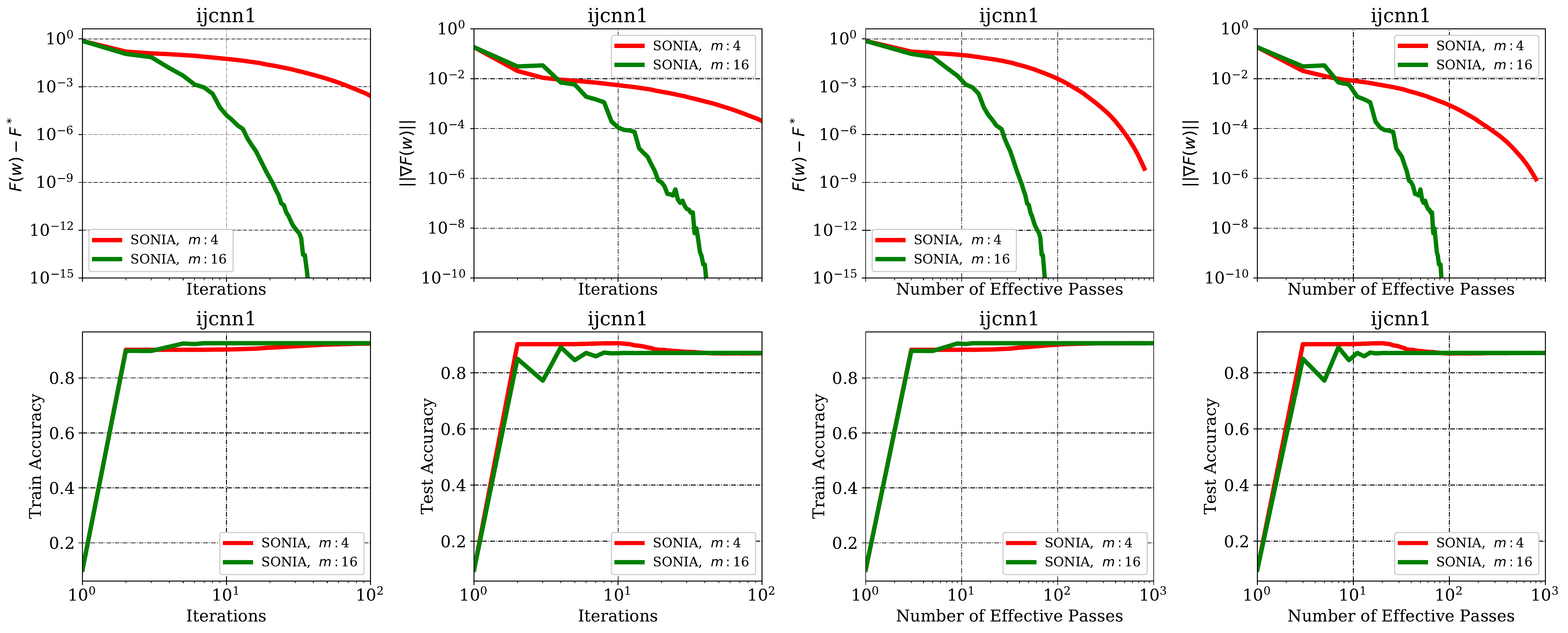}
	\caption{\texttt{ijcnn1}: Sensitivity of \texttt{SONIA} to memory,  Deterministic Logistic Regression ($\lambda = 10^{-5}$).}
	\label{fig:ijcnn1_5_FastSR1_sensitivityAnlysis_MMR}
\end{figure*}
\vspace{-10pt}
\begin{figure*}[htb]
	\centering
	\includegraphics[width=0.8\textwidth]{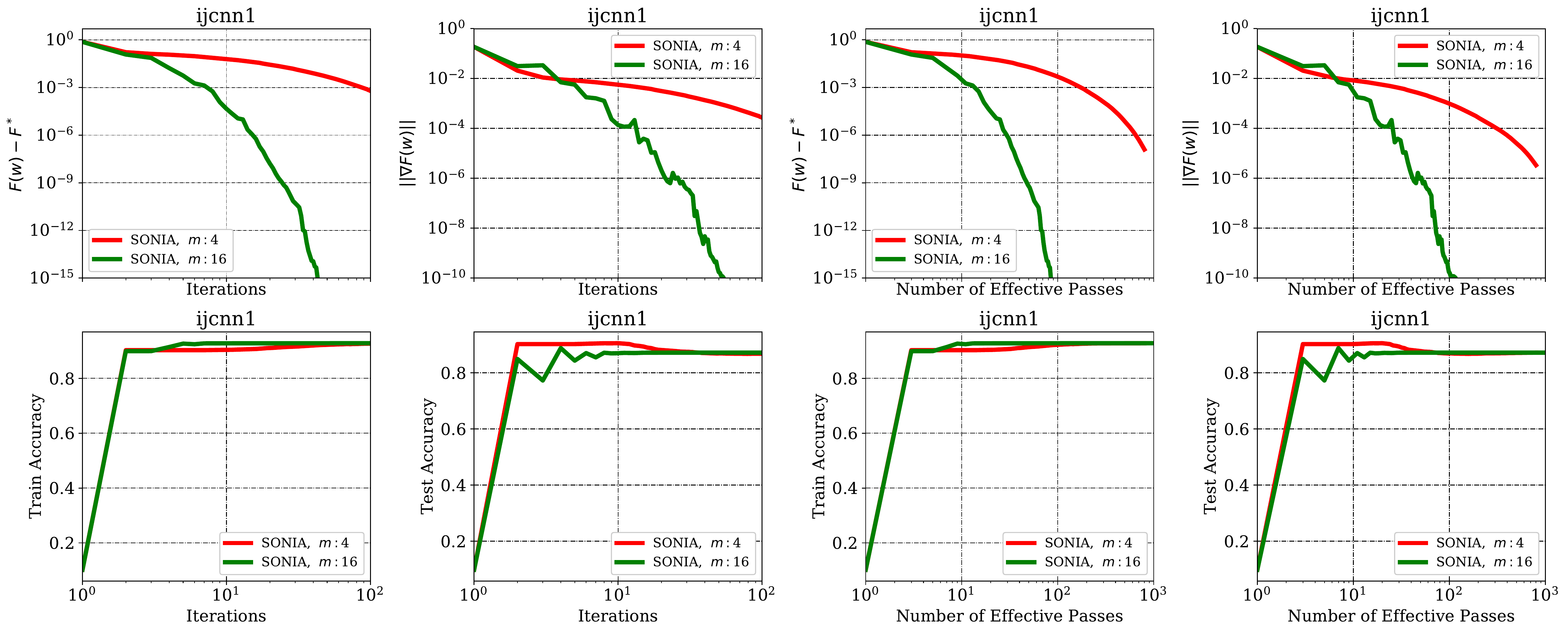}
	\caption{\texttt{ijcnn1}: Sensitivity of \texttt{SONIA} to memory,  Deterministic Logistic Regression ($\lambda = 10^{-6}$).}
	\label{fig:ijcnn1_6_FastSR1_sensitivityAnlysis_MMR}
\end{figure*}

\clearpage

\subsubsection{Sensitivity of \SONIA{} to truncation hyper-parameter}\label{sec:sensEPS}
\vspace{-10pt}
\begin{figure*}[htb]
	\centering
	\includegraphics[width=0.8\textwidth]{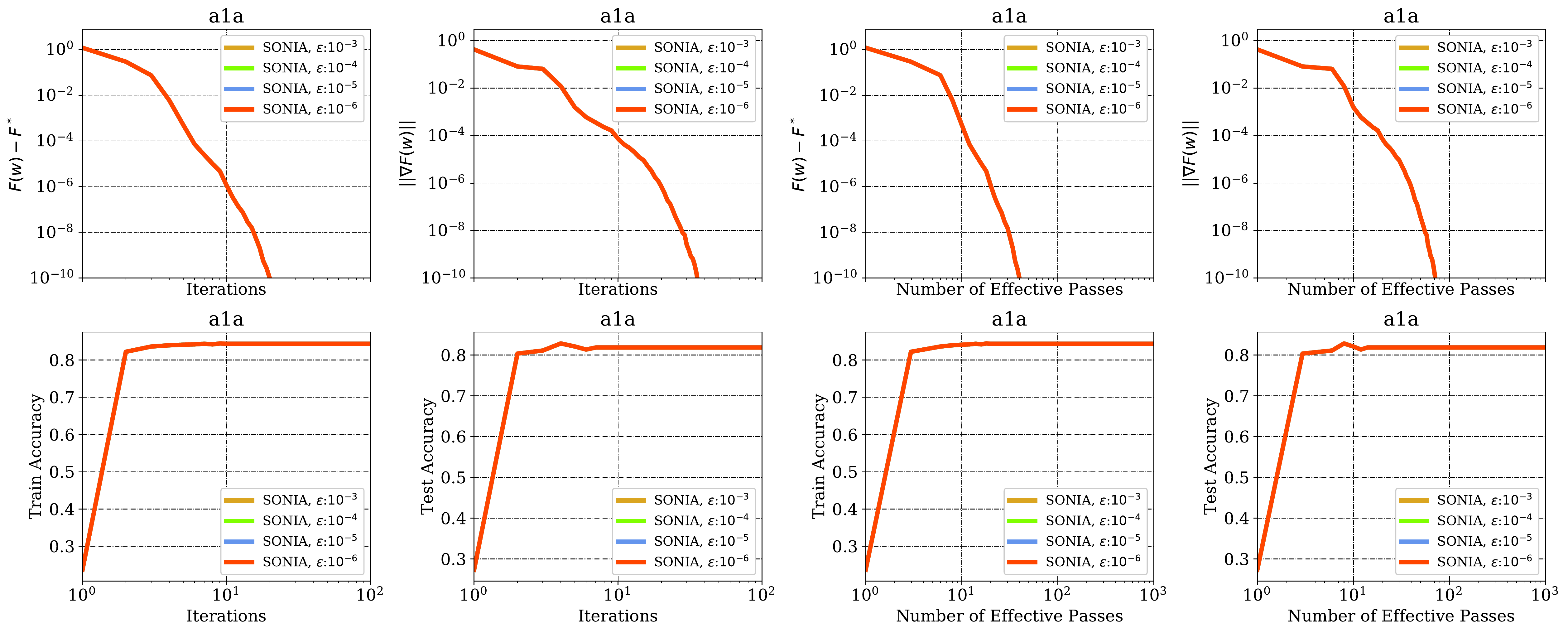}\vspace{-7pt}
	\caption{\texttt{a1a}: Sensitivity of \texttt{SONIA} to $\epsilon$,  Deterministic Logistic Regression ($\lambda = 10^{-3}$).}
	\label{fig:a1a_3}
\end{figure*}
\vspace{-10pt}
\begin{figure*}[htb]
	\centering
	\includegraphics[width=0.8\textwidth]{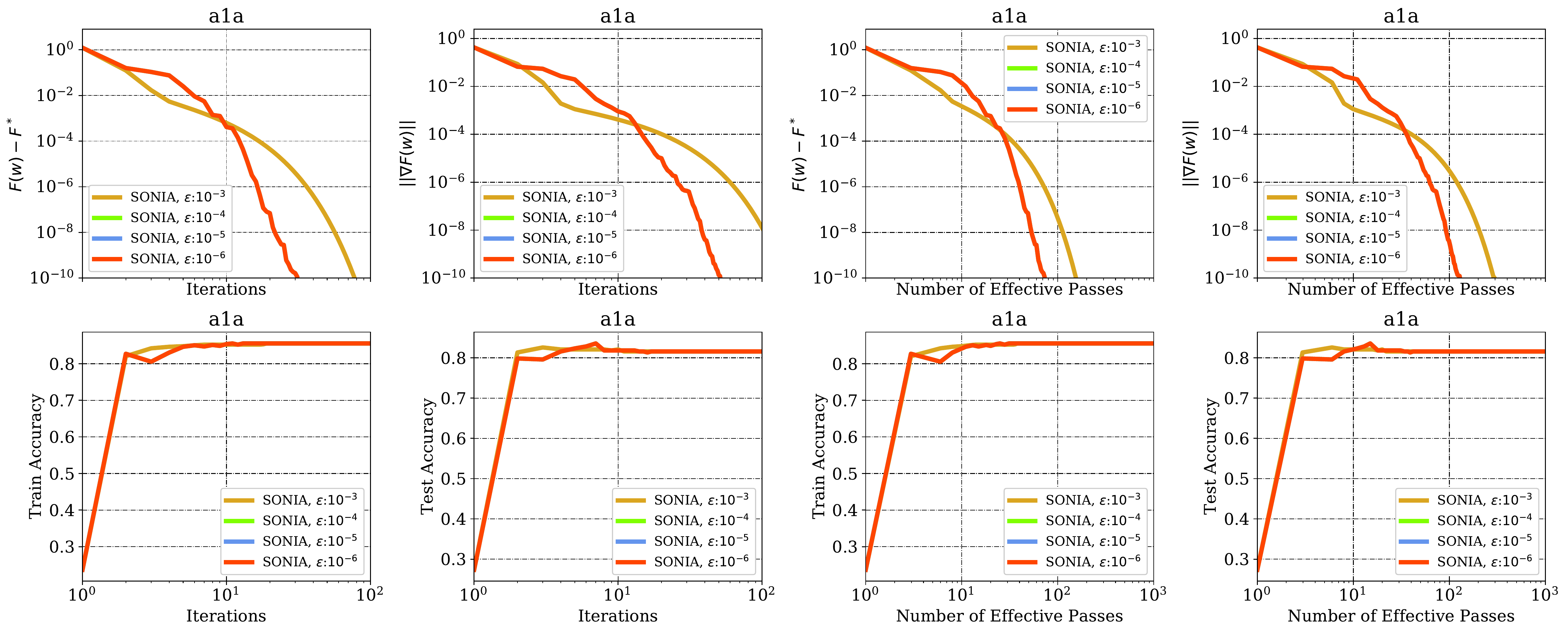}\vspace{-7pt}
	\caption{\texttt{a1a}: Sensitivity of \texttt{SONIA} to $\epsilon$,  Deterministic Logistic Regression ($\lambda = 10^{-4}$).}
	\label{fig:a1a_4}
\end{figure*}
\vspace{-10pt}
\begin{figure*}[htb]
	\centering
	\includegraphics[width=0.8\textwidth]{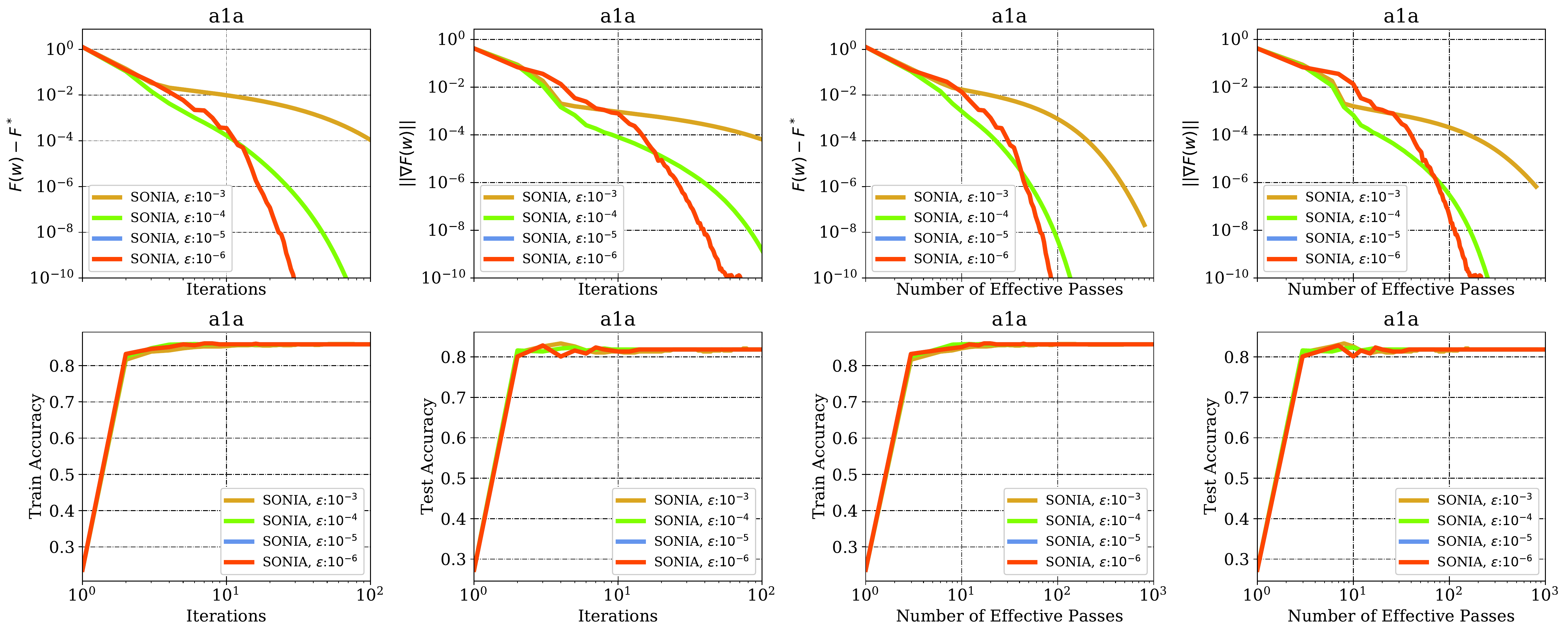}\vspace{-7pt}
	\caption{\texttt{a1a}: Sensitivity of \texttt{SONIA} to $\epsilon$,  Deterministic Logistic Regression ($\lambda = 10^{-5}$).}
	\label{fig:a1a_5}
\end{figure*}
\vspace{-10pt}
\begin{figure*}[!h]
	\centering
	\includegraphics[width=0.8\textwidth]{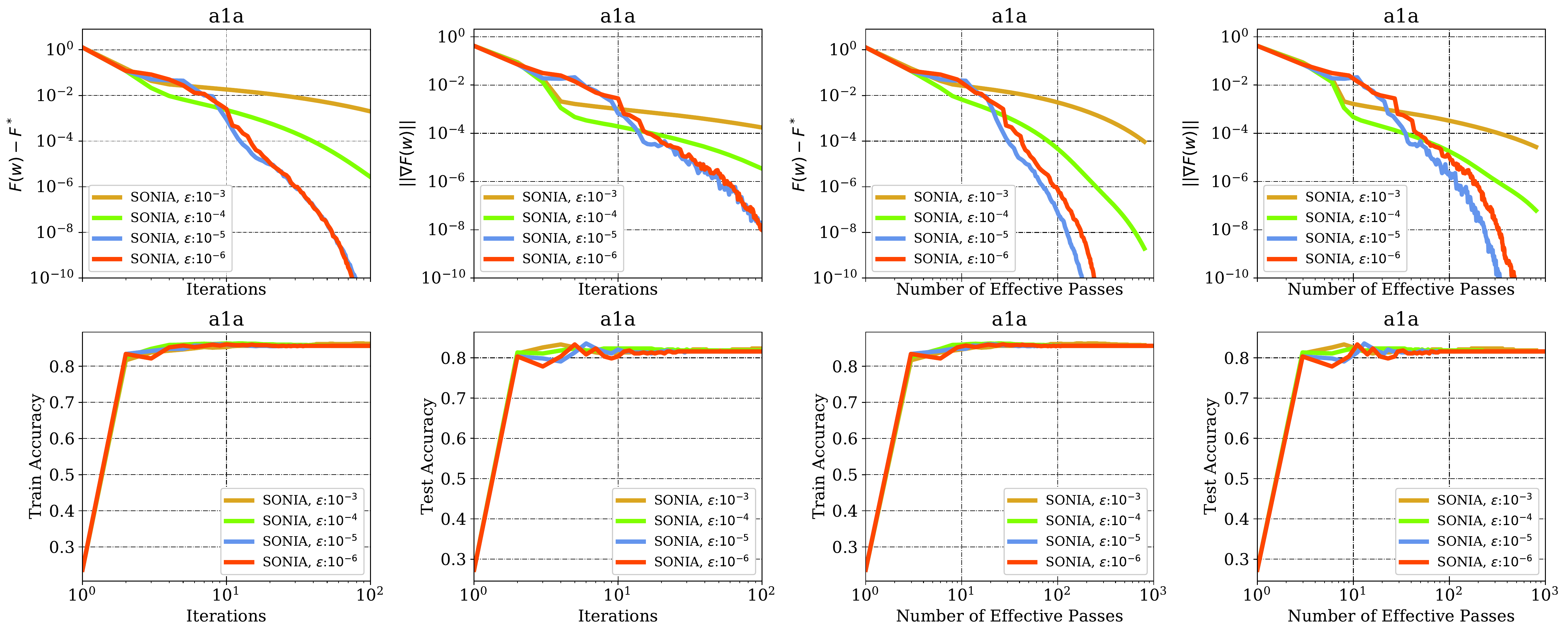}\vspace{-7pt}
	\caption{\texttt{a1a}: Sensitivity of \texttt{SONIA} to $\epsilon$,  Deterministic Logistic Regression ($\lambda = 10^{-6}$).}
	\label{fig:a1a_6}
\end{figure*}

\vspace{-10pt}
\begin{figure*}[ht]
	\centering
	\includegraphics[width=0.8\textwidth]{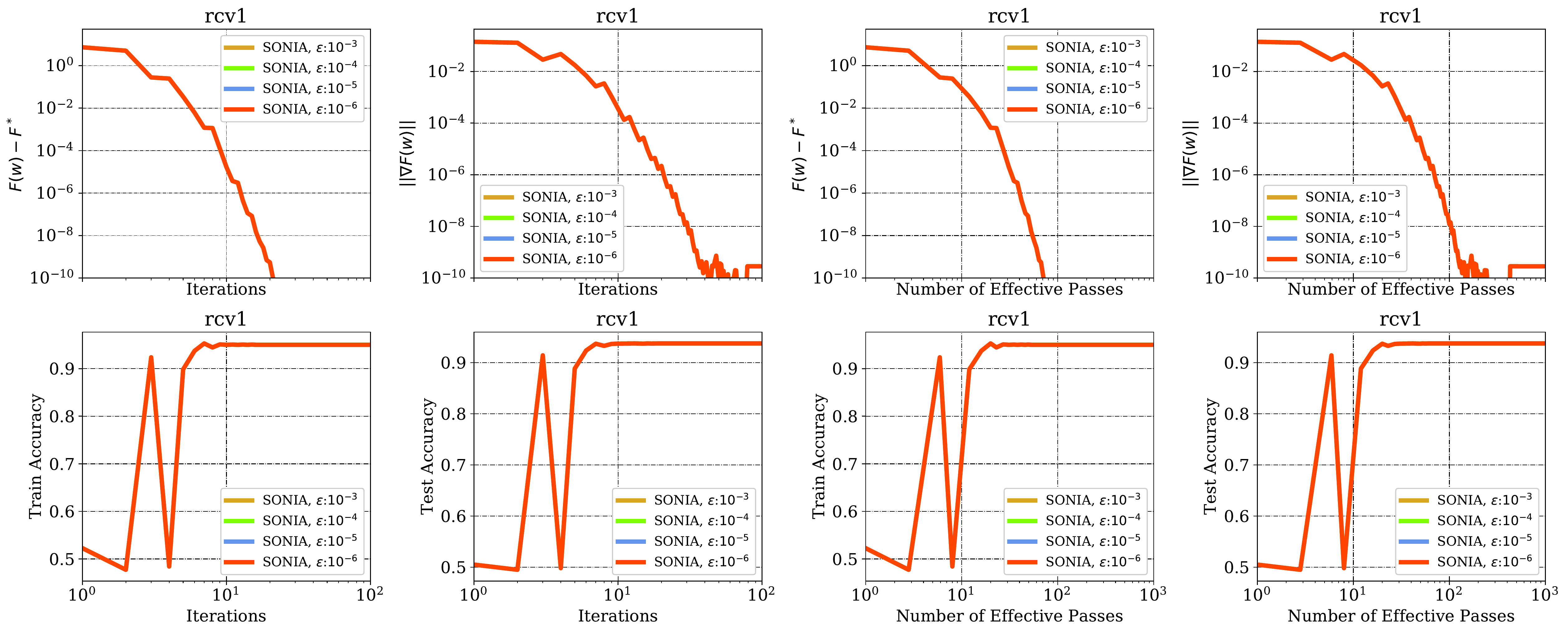}
	\caption{\texttt{rcv1}: Sensitivity of \texttt{SONIA} to $\epsilon$,  Deterministic Logistic Regression ($\lambda = 10^{-3}$).}
	\label{fig:rcv1_train_3}
\end{figure*}
\vspace{-10pt}
\begin{figure*}[ht]
	\centering
	\includegraphics[width=0.8\textwidth]{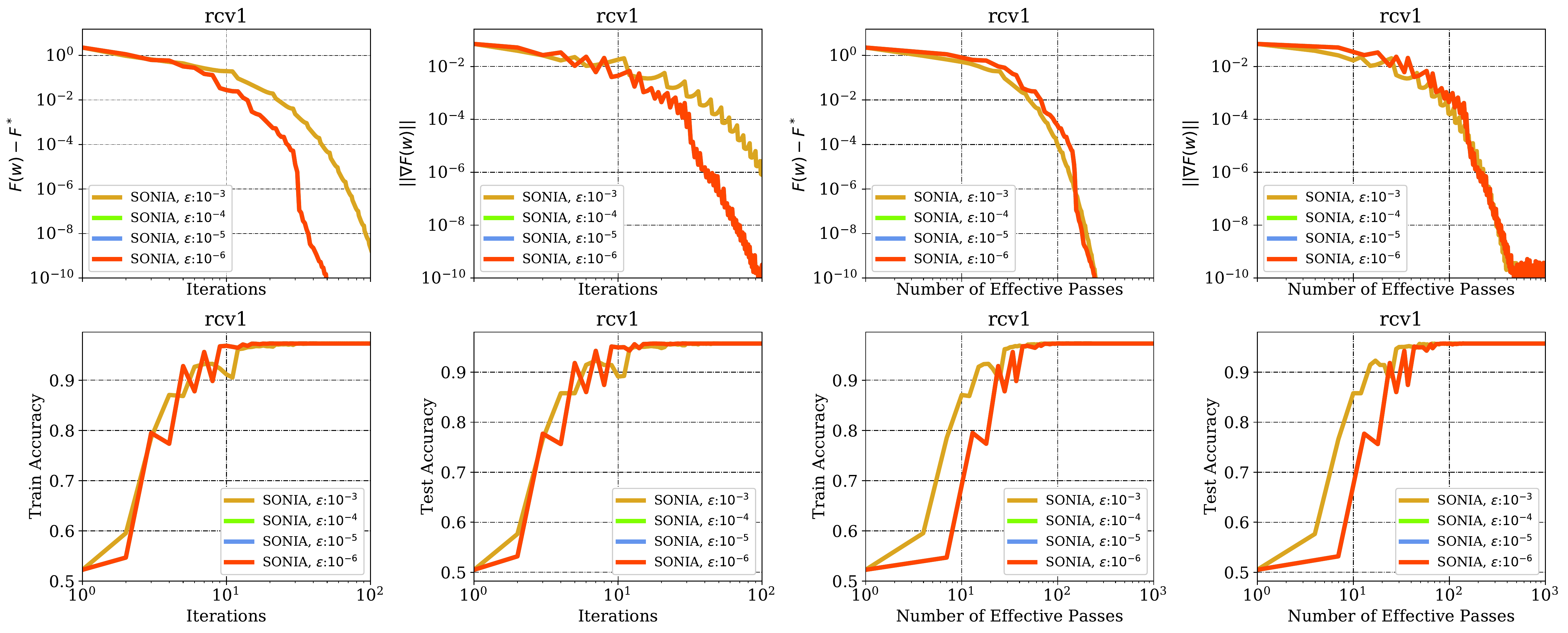}
	\caption{\texttt{rcv1}: Sensitivity of \texttt{SONIA} to $\epsilon$,  Deterministic Logistic Regression ($\lambda = 10^{-4}$).}
	\label{fig:rcv1_train_4}
\end{figure*}
\vspace{-10pt}
\begin{figure*}[ht]
	\centering
	\includegraphics[width=0.8\textwidth]{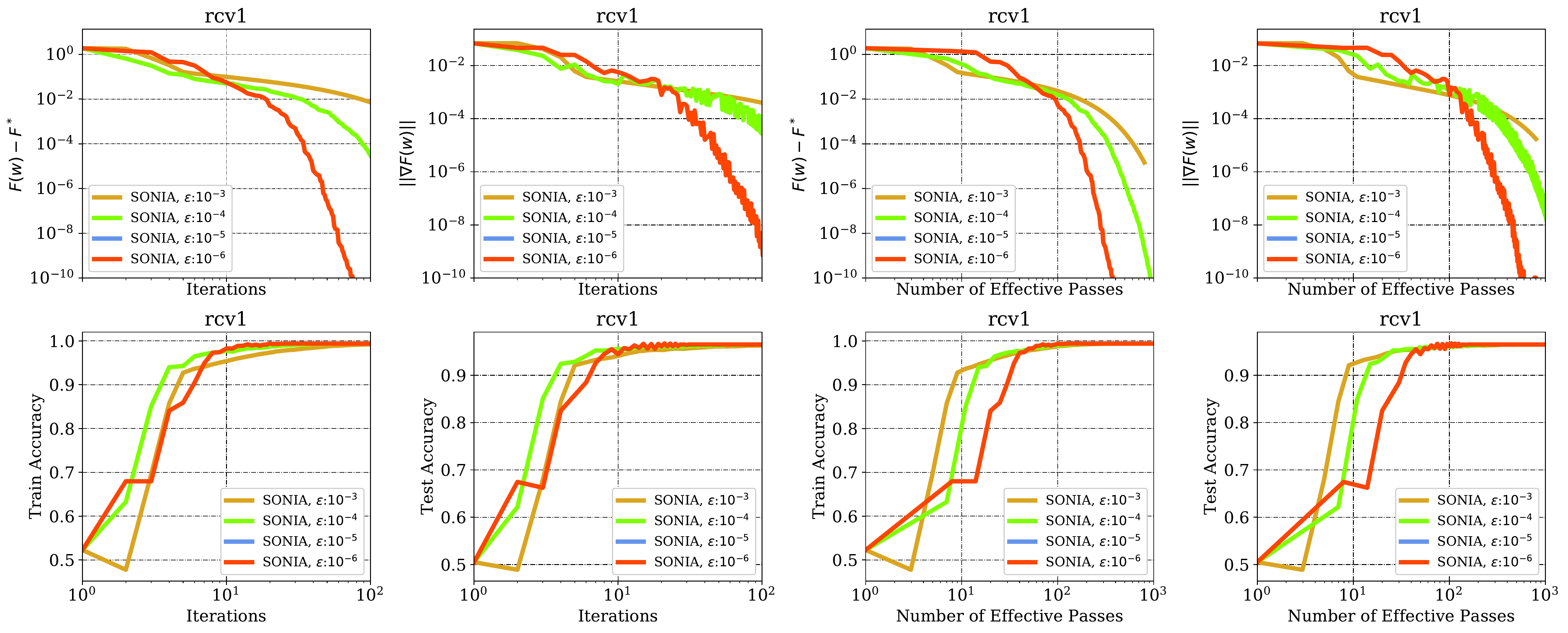}
	\caption{\texttt{rcv1}: Sensitivity of \texttt{SONIA} to $\epsilon$,  Deterministic Logistic Regression ($\lambda = 10^{-5}$).}
	\label{fig:gisettercv1_train_5}
\end{figure*}
\vspace{-10pt}
\begin{figure*}[htb]
	\centering
	\includegraphics[width=0.8\textwidth]{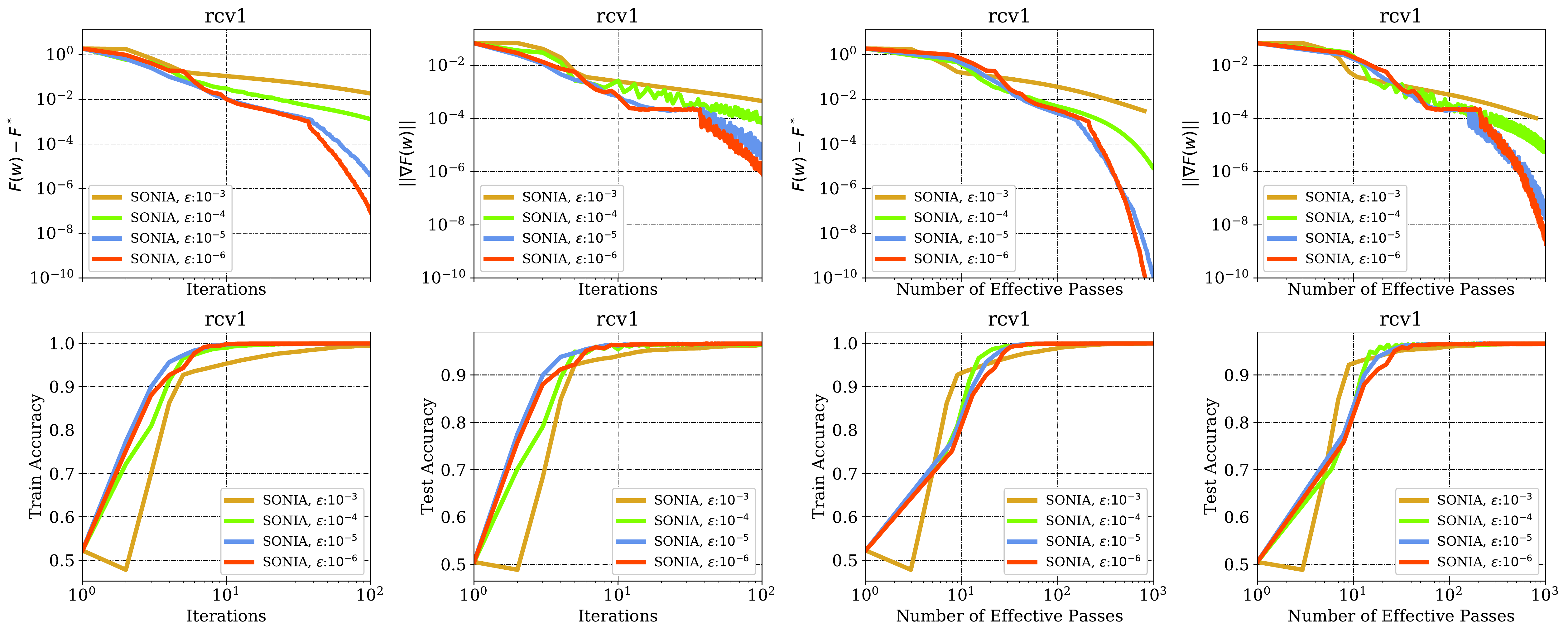}
	\caption{\texttt{rcv1}: Sensitivity of \texttt{SONIA} to $\epsilon$,  Deterministic Logistic Regression ($\lambda = 10^{-6}$).}
	\label{fig:rcv1_train_6}
\end{figure*}


\vspace{-10pt}
\begin{figure*}[ht]
	\centering
	\includegraphics[width=0.8\textwidth]{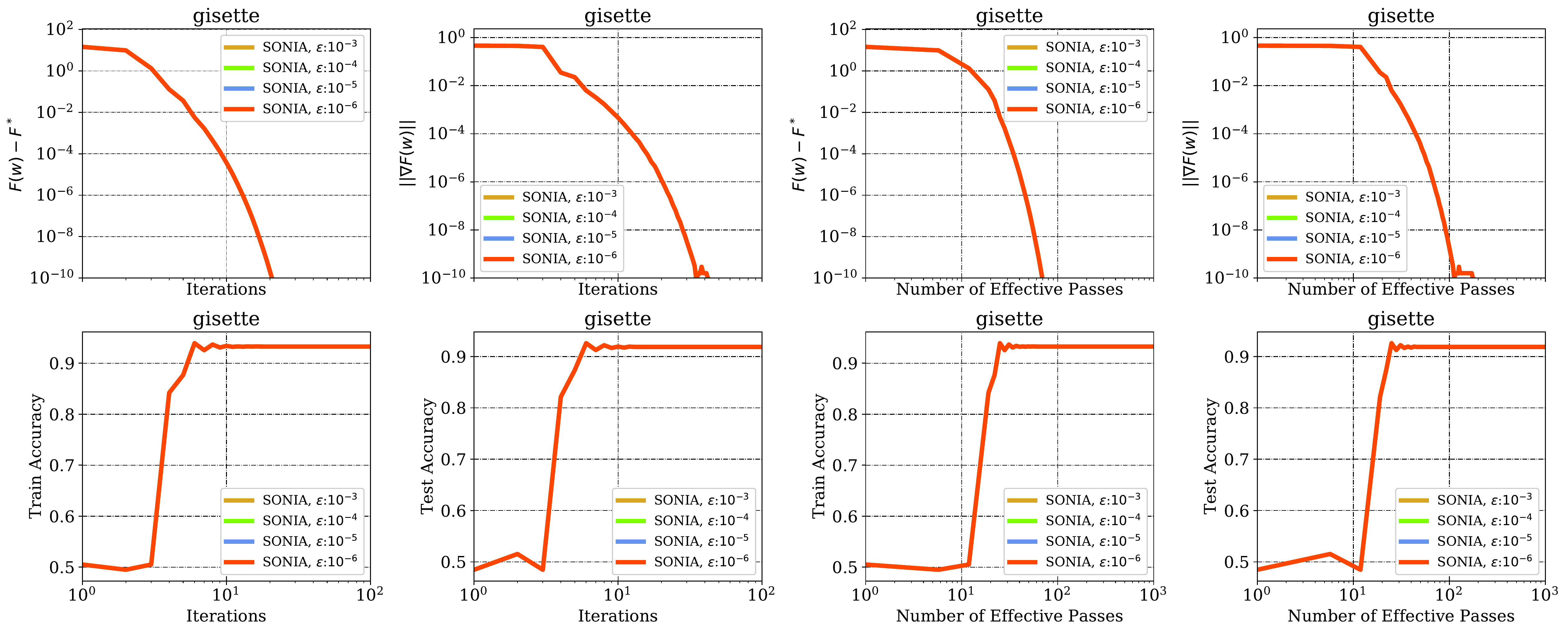}
	\caption{\texttt{gisette}: Sensitivity of \texttt{SONIA} to $\epsilon$,  Deterministic Logistic Regression ($\lambda = 10^{-3}$).}
	\label{fig:gisette_scale_3}
\end{figure*}
\vspace{-10pt}
\begin{figure*}[ht]
	\centering
	\includegraphics[width=0.8\textwidth]{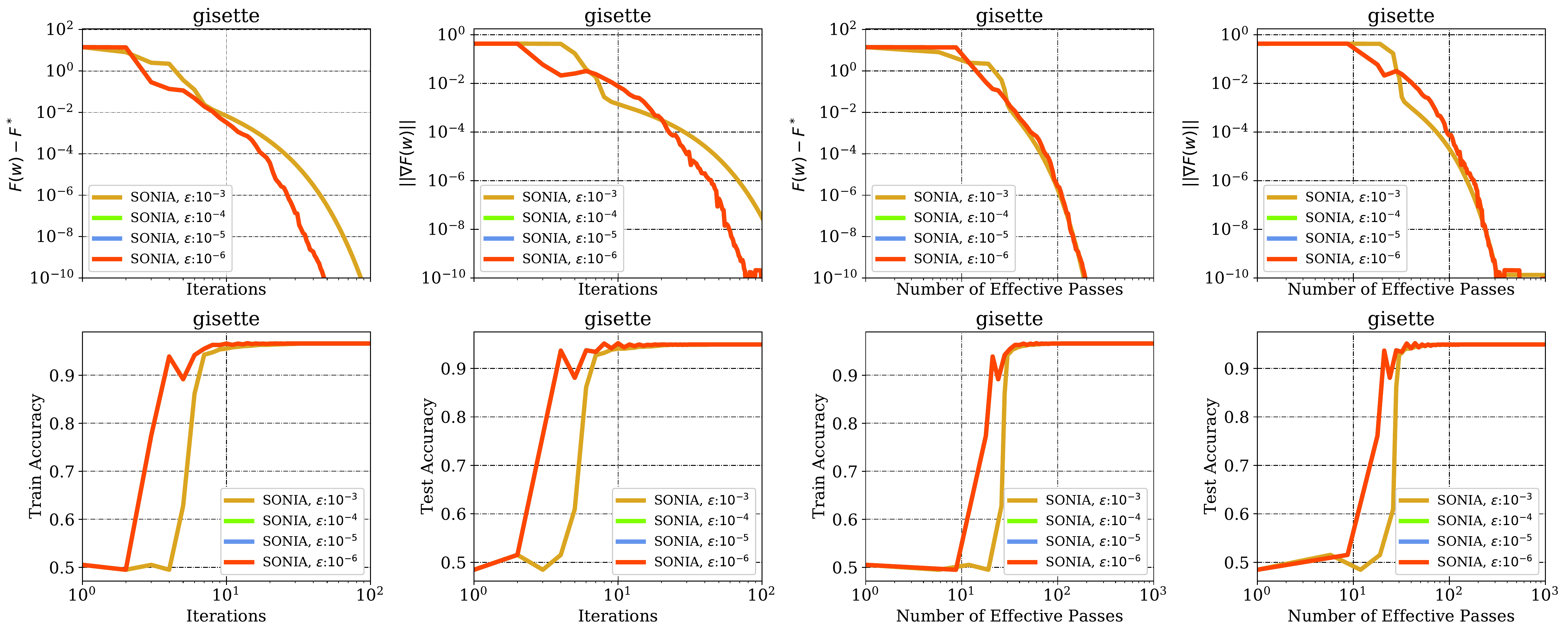}
	\caption{\texttt{gisette}: Sensitivity of \texttt{SONIA} to $\epsilon$,  Deterministic Logistic Regression ($\lambda = 10^{-4}$).}
	\label{fig:gisette_scale_4}
\end{figure*}
\vspace{-10pt}
\begin{figure*}[ht]
	\centering
	\includegraphics[width=0.8\textwidth]{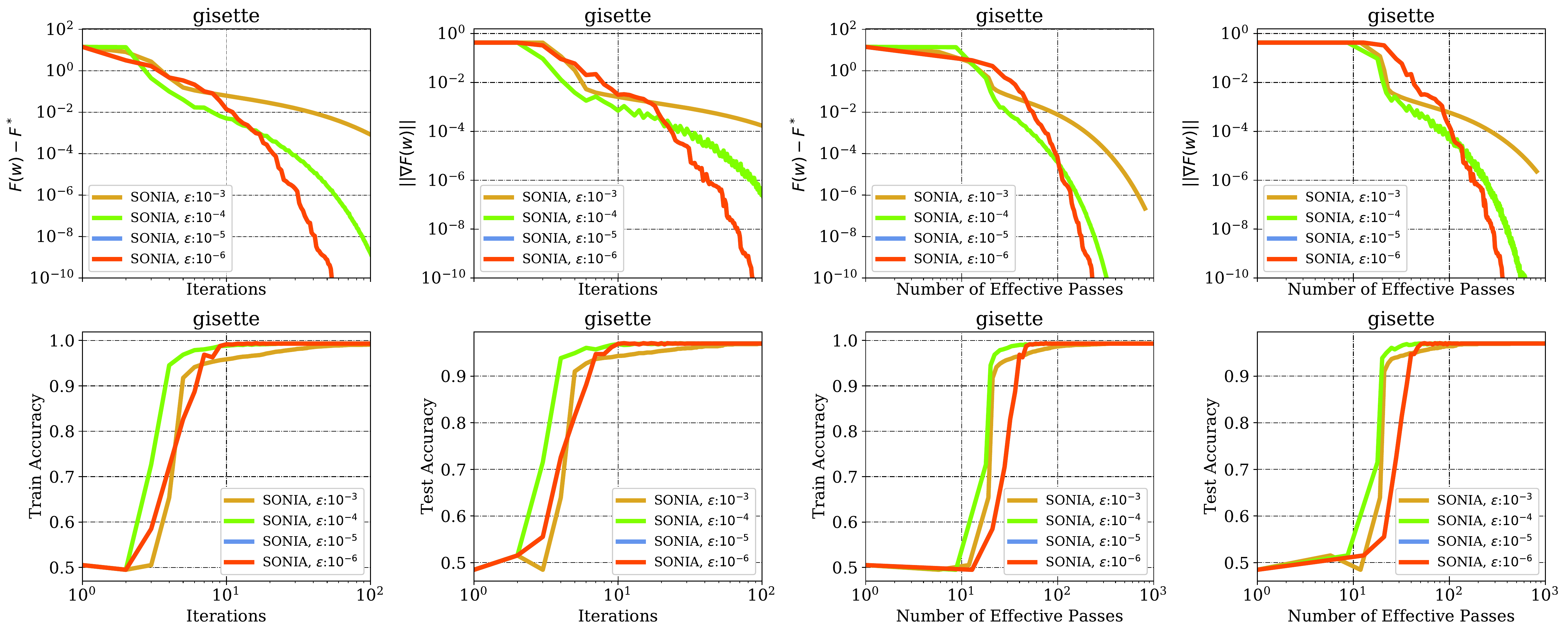}
	\caption{\texttt{gisette}: Sensitivity of \texttt{SONIA} to $\epsilon$,  Deterministic Logistic Regression ($\lambda = 10^{-5}$).}
	\label{fig:gisette_scale_5}
\end{figure*}
\vspace{-10pt}
\begin{figure*}[!h]
	\centering
	\includegraphics[width=0.8\textwidth]{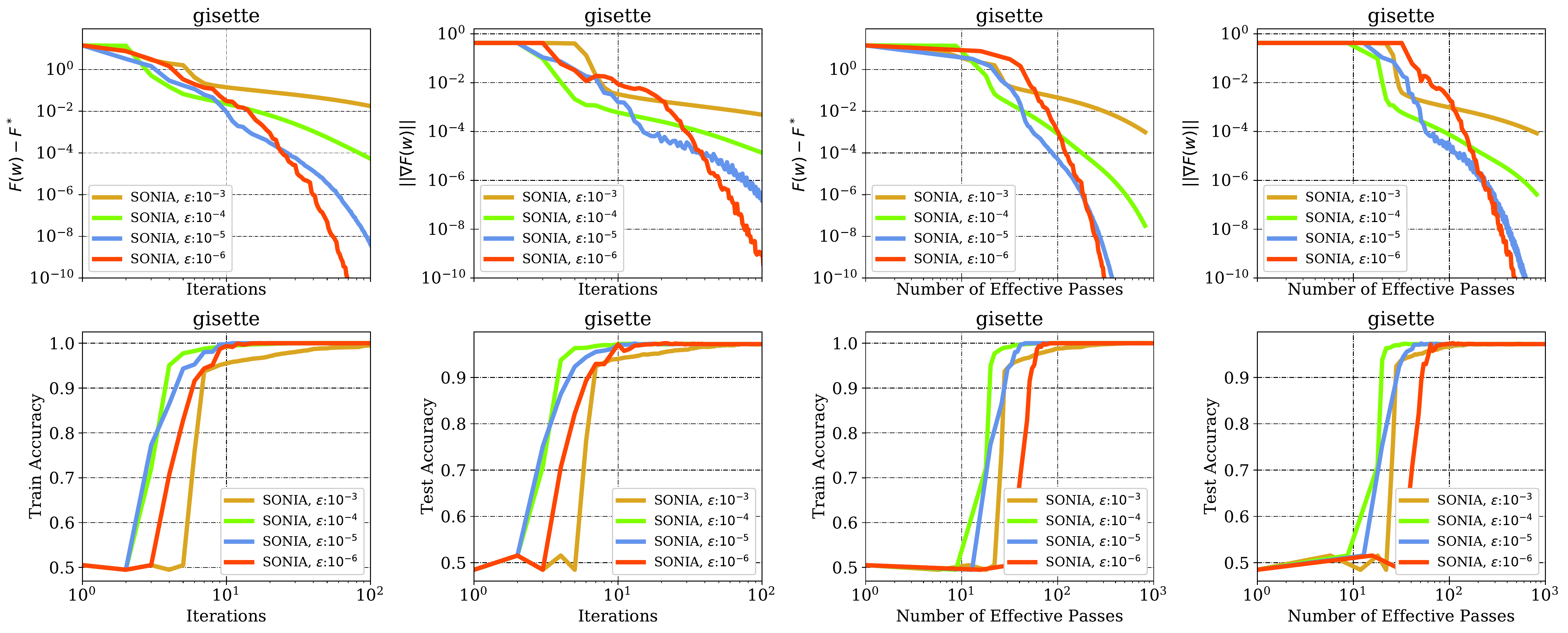}
	\caption{\texttt{gisette}: Sensitivity of \texttt{SONIA} to $\epsilon$,  Deterministic Logistic Regression ($\lambda = 10^{-6}$).}
	\label{fig:gisette_scale_6}
\end{figure*}

\vspace{-10pt}
\begin{figure*}[ht]
	\centering
	\includegraphics[width=0.8\textwidth]{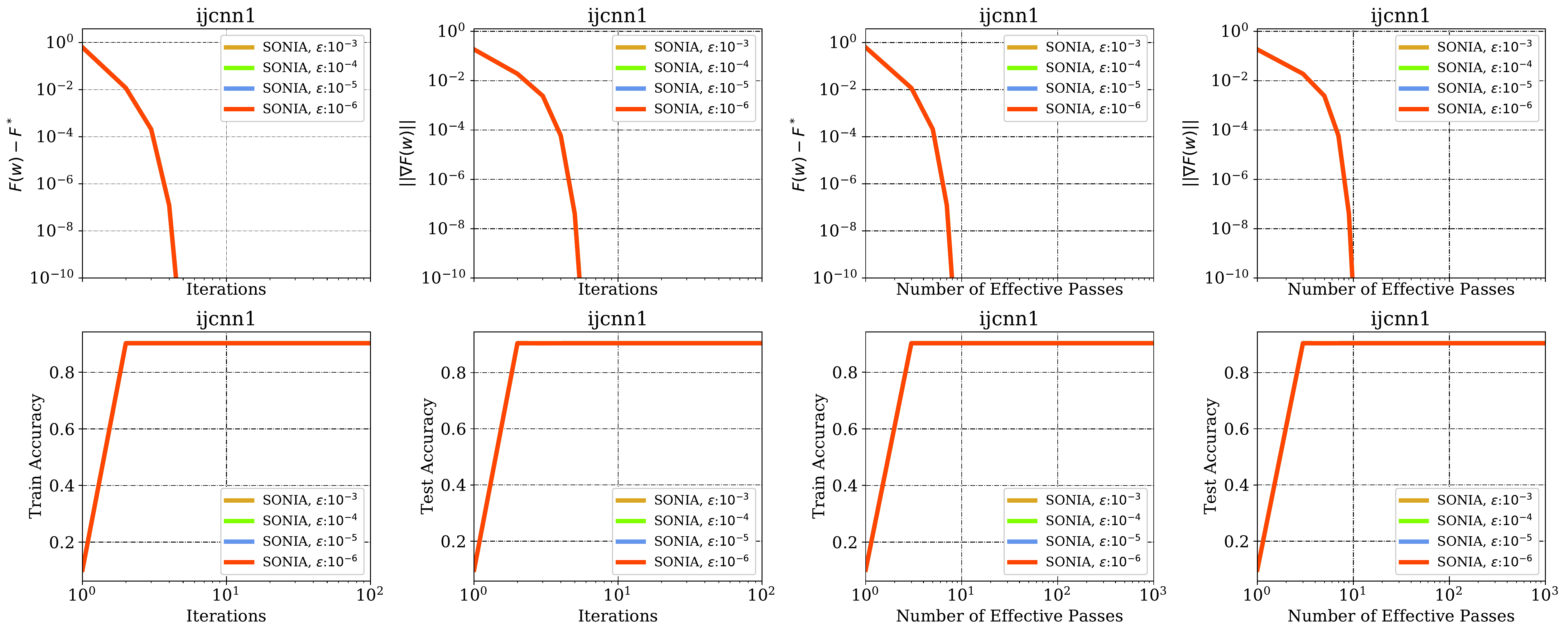}
	\caption{\texttt{ijcnn1}: Sensitivity of \texttt{SONIA} to $\epsilon$,  Deterministic Logistic Regression ($\lambda = 10^{-3}$).}
	\label{fig:ijcnn1_3}
\end{figure*}
\vspace{-10pt}
\begin{figure*}[ht]
	\centering
	\includegraphics[width=0.8\textwidth]{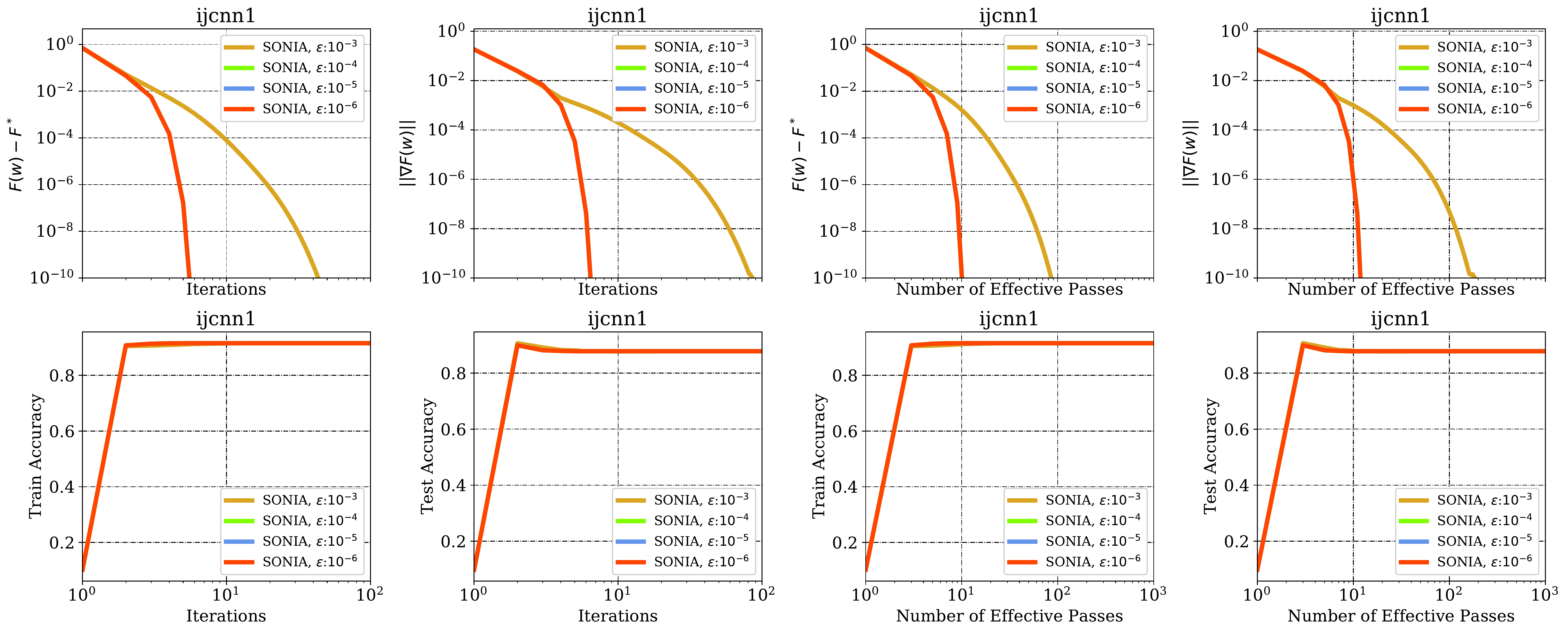}
	\caption{\texttt{ijcnn1}: Sensitivity of \texttt{SONIA} to $\epsilon$,  Deterministic Logistic Regression ($\lambda = 10^{-4}$).}
	\label{fig:ijcnn1_4}
\end{figure*}
\vspace{-10pt}
\begin{figure*}[ht]
	\centering
	\includegraphics[width=0.8\textwidth]{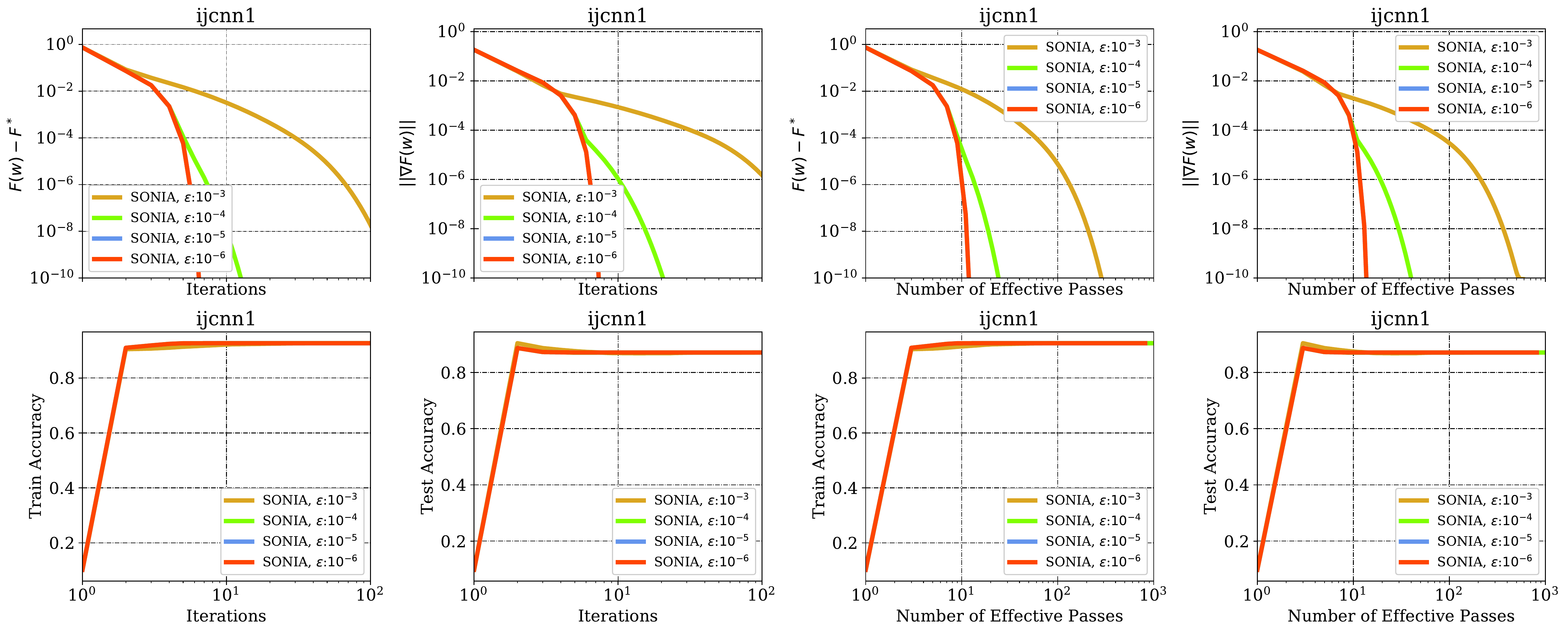}
	\caption{\texttt{ijcnn1}: Sensitivity of \texttt{SONIA} to $\epsilon$,  Deterministic Logistic Regression ($\lambda = 10^{-5}$).}
	\label{fig:gisette_scale_5}
\end{figure*}
\vspace{-10pt}
\begin{figure*}[ht]
	\centering
	\includegraphics[width=0.8\textwidth]{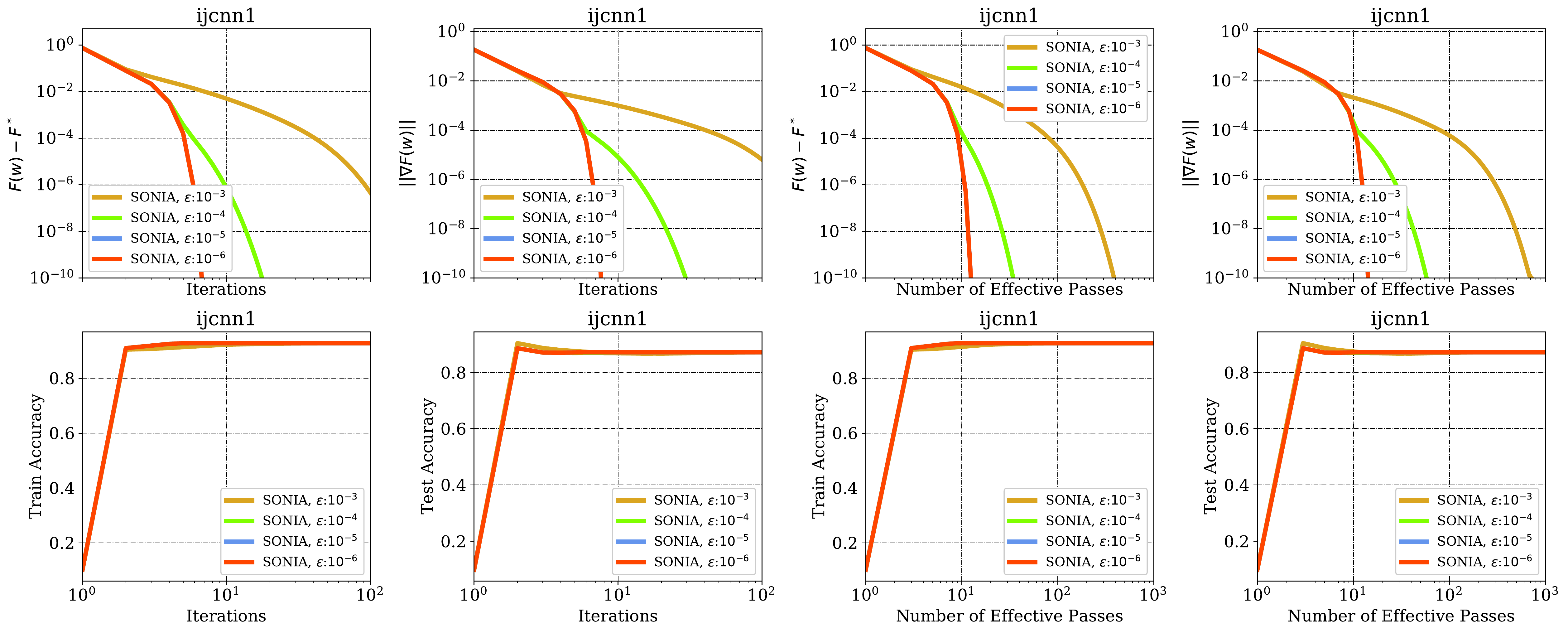}
	\caption{\texttt{ijcnn1}: Sensitivity of \texttt{SONIA} to $\epsilon$,  Deterministic Logistic Regression ($\lambda = 10^{-6}$).}
	\label{fig:ijcnn1_6}
\end{figure*}

\clearpage
\subsection{Additional Numerical Experiments: Deterministic Nonconvex Functions} \label{sec:moreNonConvResults}
\vspace{-10pt}
\begin{figure*}[htb]
	\centering
{	\includegraphics[trim=10 110 10 110,clip, width=0.6\textwidth]{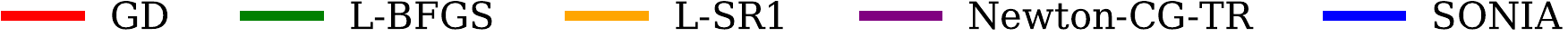}
}
	\includegraphics[width=0.8\textwidth]{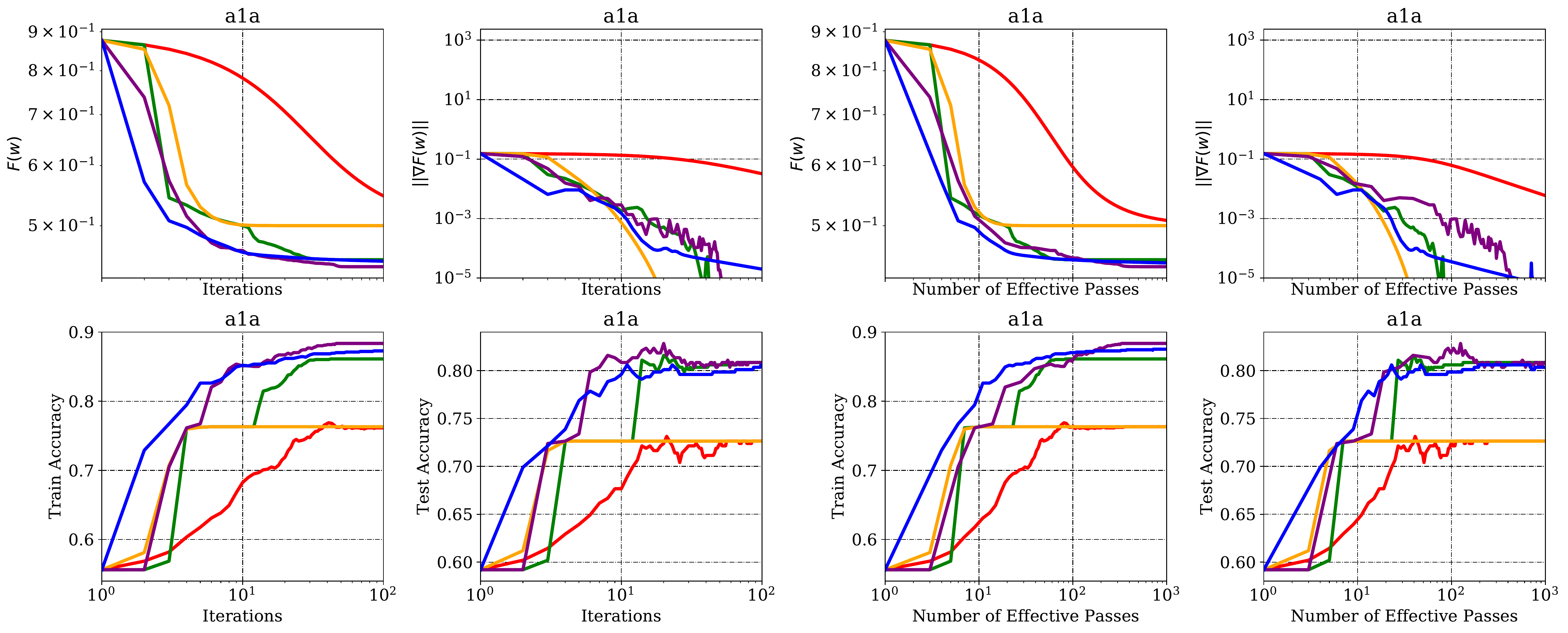}\vspace{-7pt}
	\caption{\texttt{a1a}: Deterministic Non-Linear Least Square}
	\label{fig:a1a_DET_NL_LS}
\end{figure*}
\vspace{-10pt}
\begin{figure*}[!h]
	\centering
	\includegraphics[width=0.8\textwidth]{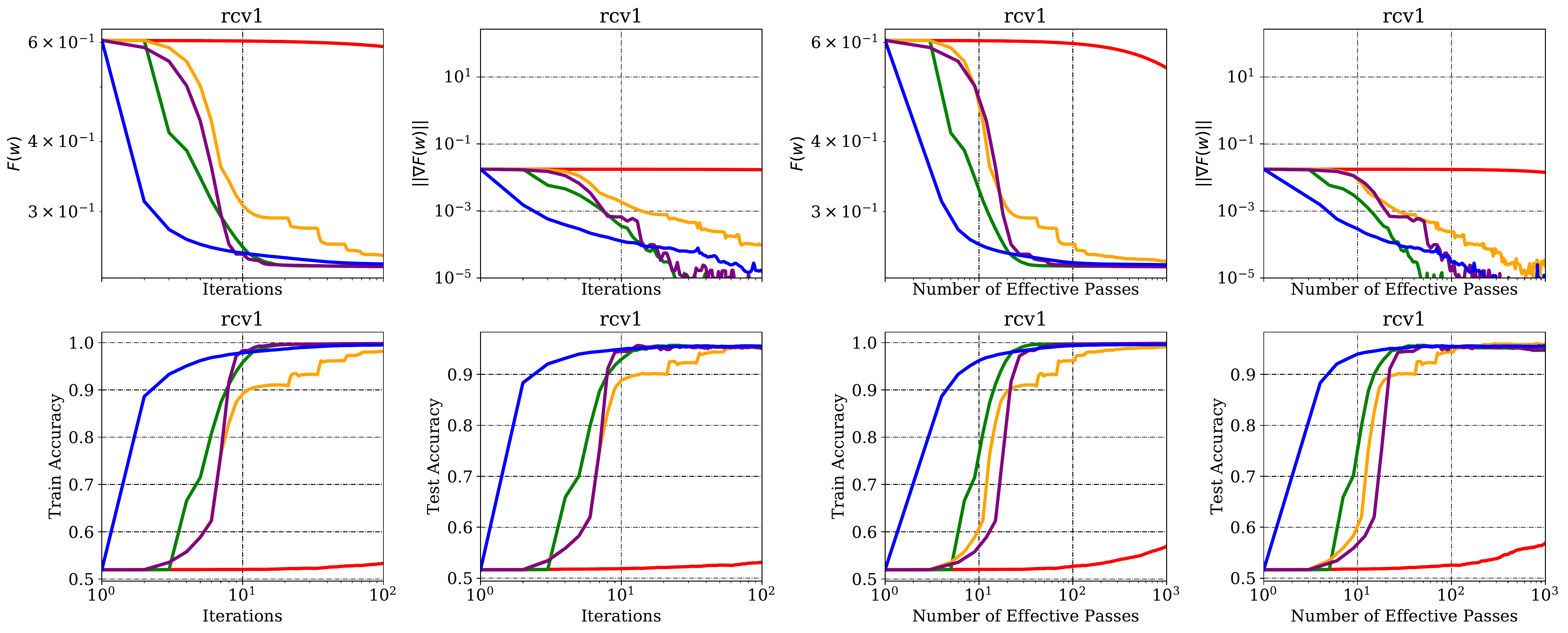}\vspace{-7pt}
	\caption{\texttt{rcv1}: Deterministic Non-Linear Least Square}
	\label{fig:rcv1_train_DET_NL_LS}
\end{figure*}
\vspace{-10pt}
\begin{figure*}[htb]
	\centering
	\includegraphics[width=0.8\textwidth]{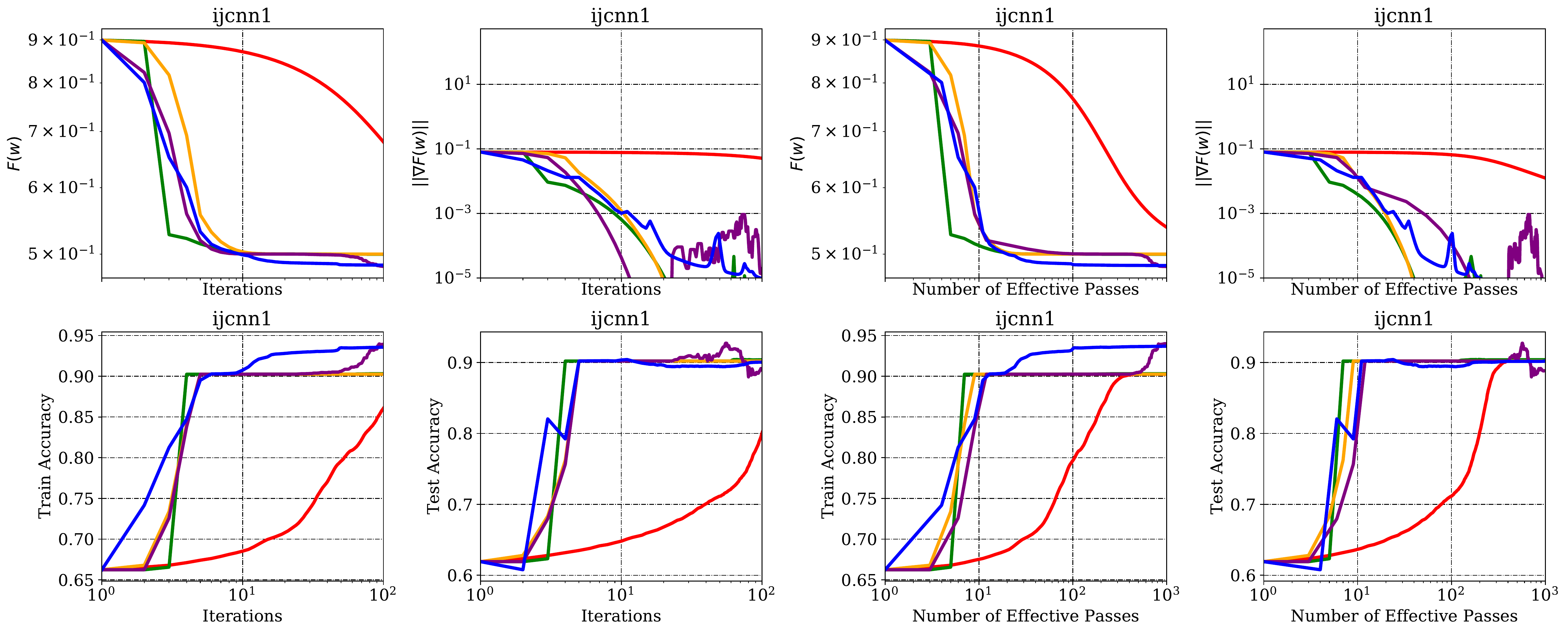}\vspace{-7pt}
	\caption{\texttt{ijcnn1}: Deterministic Non-Linear Least Square}
	\label{fig:ijcnn1_DET_NL_LS}
\end{figure*}
\vspace{-10pt}
\begin{figure*}[!h]
	\centering
	\includegraphics[width=0.8\textwidth]{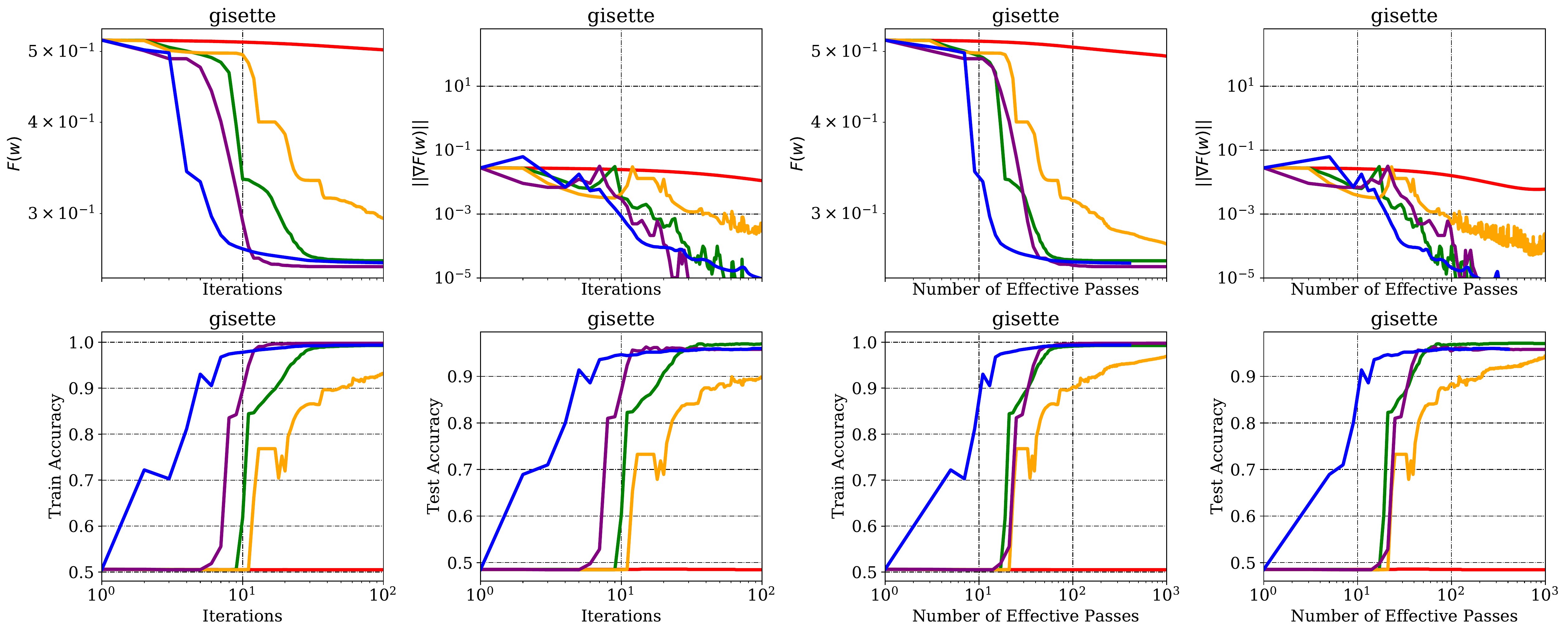}\vspace{-7pt}
	\caption{\texttt{gisette}: Deterministic Non-Linear Least Square}
	\label{fig:gisette_scale_DET_NL_LS}
\end{figure*}

\clearpage
\subsection{Additional Numerical Experiments: Stochastic Strongly Convex Functions} \label{sec:moreStochStCOnvResults}


\begin{figure*}[htb]
	\centering
{	\includegraphics[trim=10 100 10 110,clip, width=0.95\textwidth]{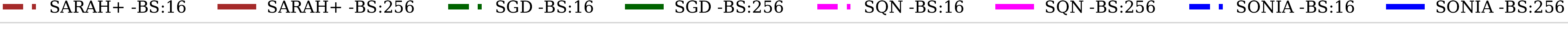}
}
	\includegraphics[width=0.99\textwidth]{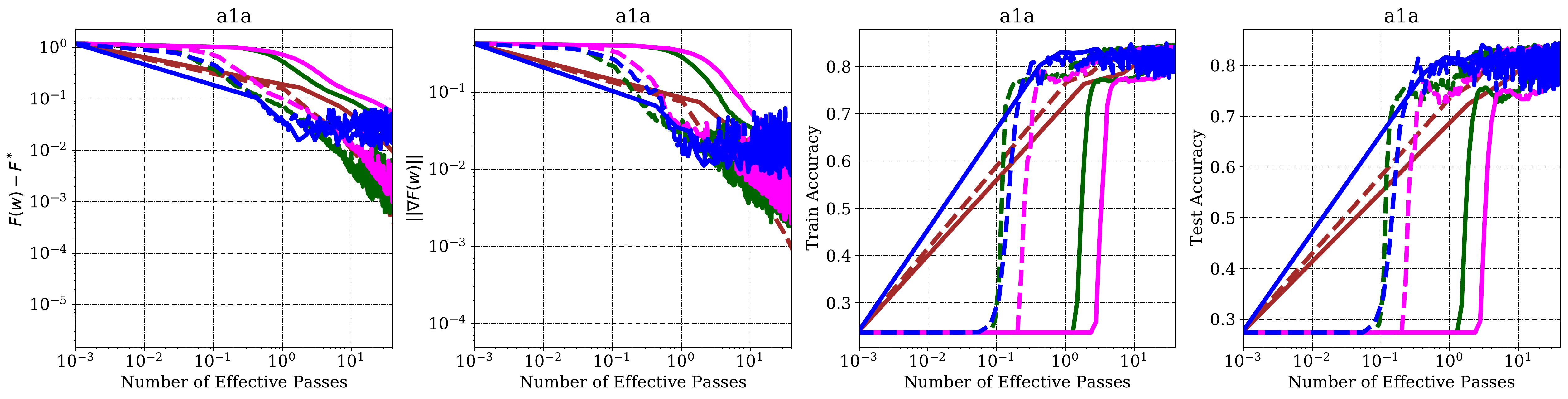}
	\caption{\texttt{a1a}: Stochastic Logistic Regression ($\lambda = 10^{-3}$).}
	\label{fig:stocha1a_3}
\end{figure*}

\begin{figure*}[htb]
	\centering
	\includegraphics[width=0.99\textwidth]{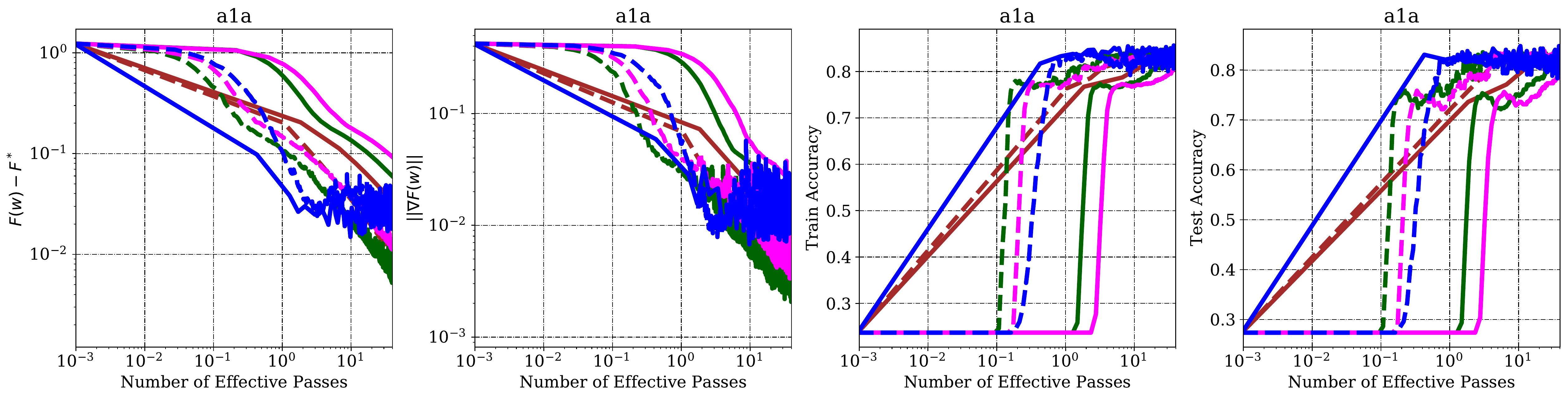}
	\caption{\texttt{a1a}: Stochastic Logistic Regression ($\lambda = 10^{-4}$).}
	\label{fig:stocha1a_4}
\end{figure*}

\begin{figure*}[htb]
	\includegraphics[width=0.99\textwidth]{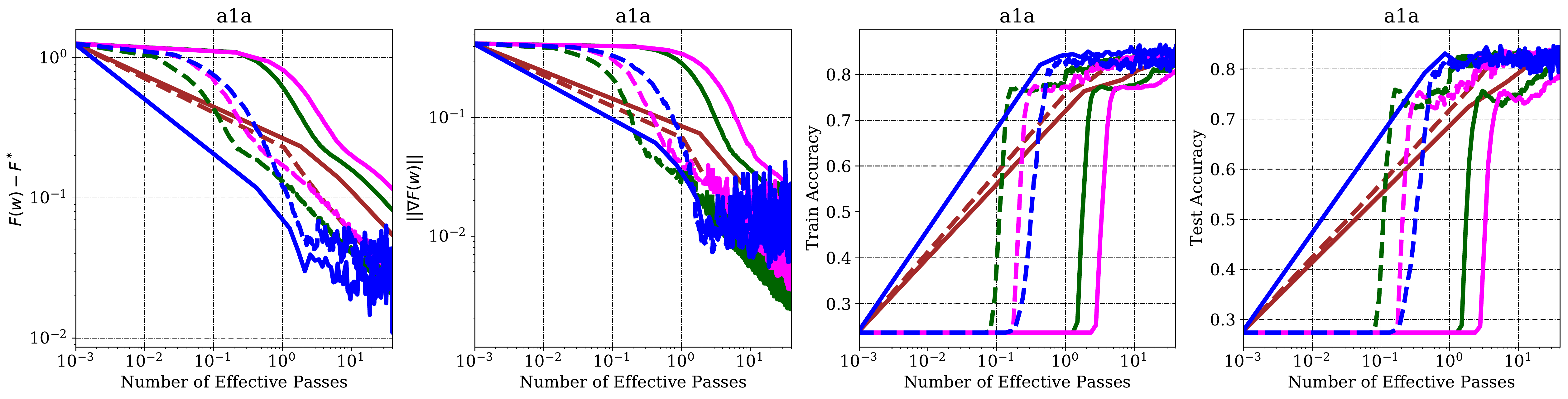}
	\caption{\texttt{a1a}: Stochastic Logistic Regression ($\lambda = 10^{-5}$).}
	\label{fig:stocha1a_5}
\end{figure*}

\begin{figure*}[!h]
	\centering
	\includegraphics[width=0.99\textwidth]{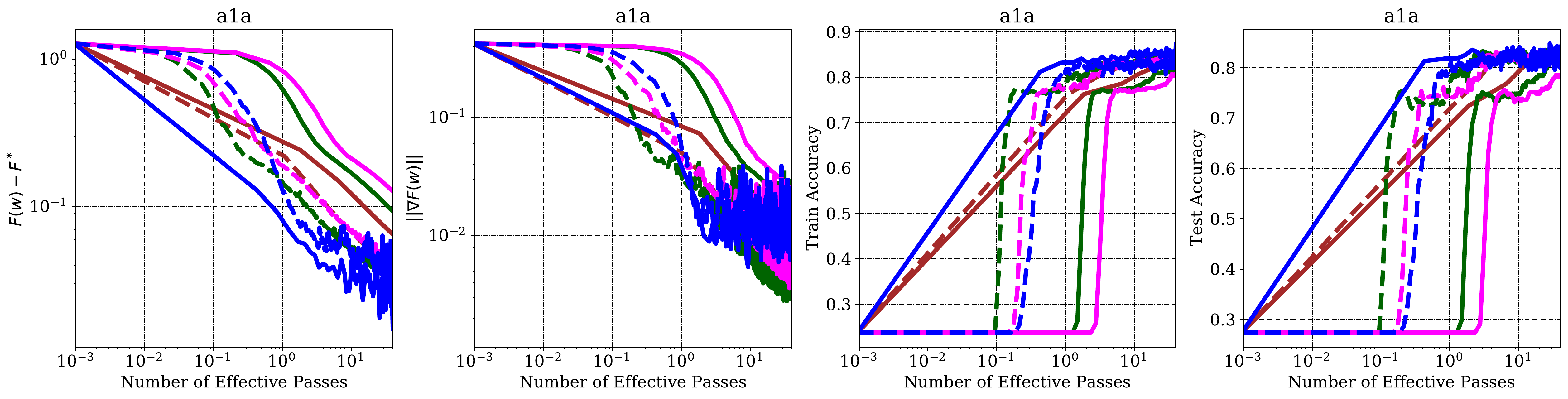}
	\caption{\texttt{a1a}: Stochastic Logistic Regression ($\lambda = 10^{-6}$).}
	\label{fig:stocha1a_6}
\end{figure*}


\begin{figure*}[ht]
	\centering
{	\includegraphics[trim=10 100 10 110,clip, width=0.95\textwidth]{Figures/a1a_comp_stoch_NLLS_leg.pdf}
}
	\includegraphics[width=0.99\textwidth]{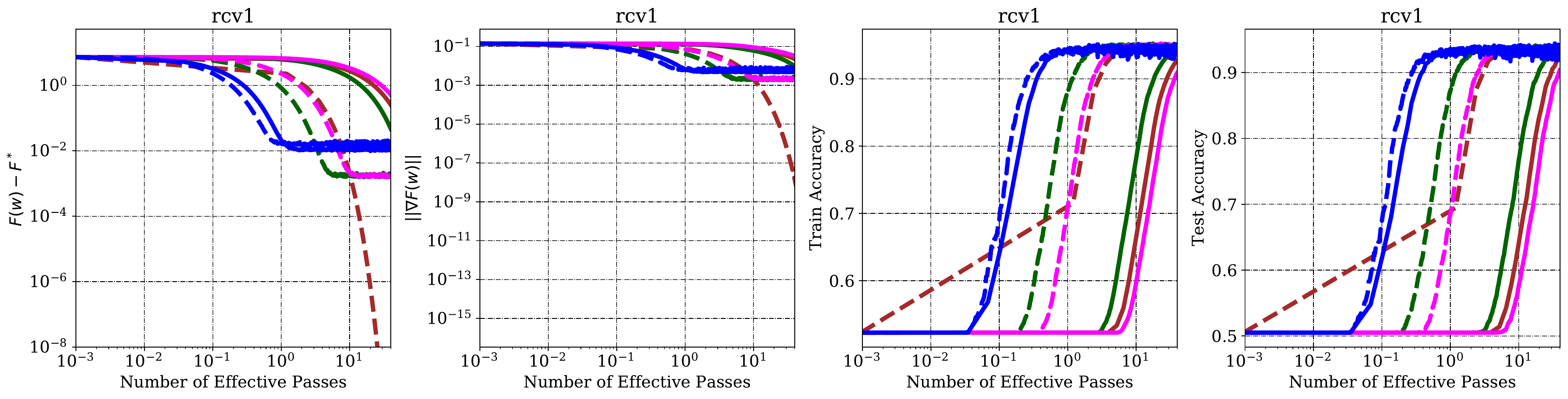}
	\caption{\texttt{rcv1}: Stochastic Logistic Regression ($\lambda = 10^{-3}$).}
	\label{fig:stochrcv1_train_3}
\end{figure*}

\begin{figure*}[ht]
	\centering
	\includegraphics[width=0.99\textwidth]{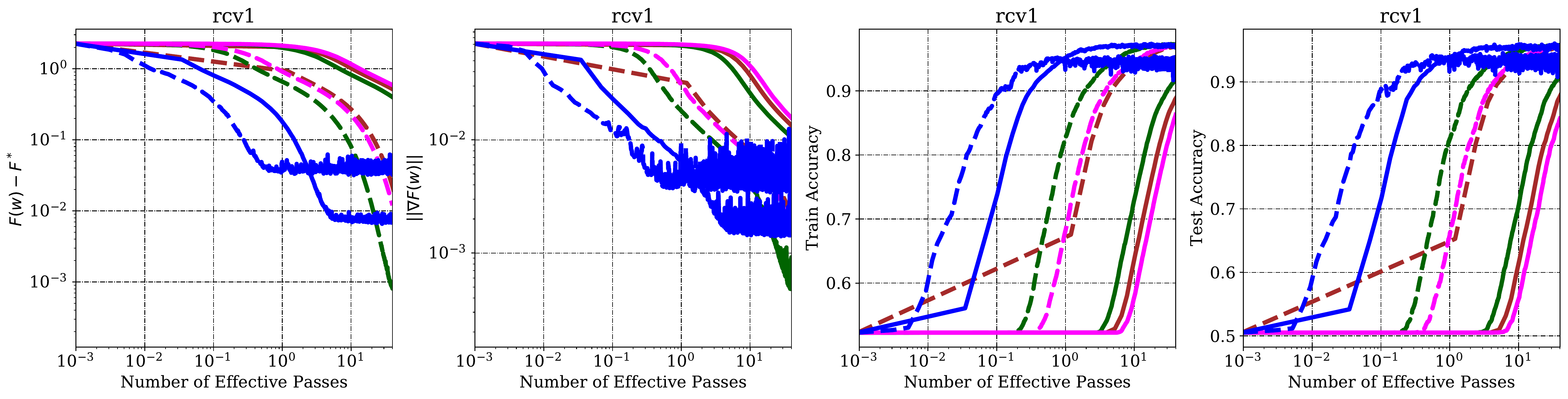}
	\caption{\texttt{rcv1}: Stochastic Logistic Regression ($\lambda = 10^{-4}$).}
	\label{fig:stochrcv1_train_4}
\end{figure*}

\begin{figure*}[ht]
	\centering
	\includegraphics[width=0.99\textwidth]{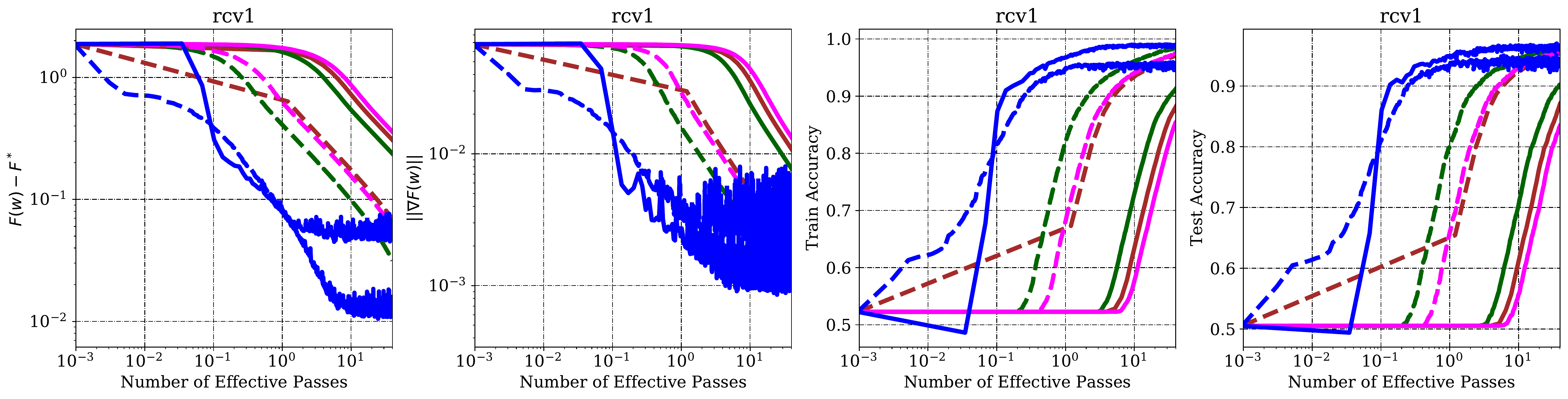}
	\caption{\texttt{rcv1}: Stochastic Logistic Regression ($\lambda = 10^{-5}$).}
	\label{fig:stochrcv1_train_5}
\end{figure*}

\begin{figure*}[htb]
	\centering
	\includegraphics[width=0.99\textwidth]{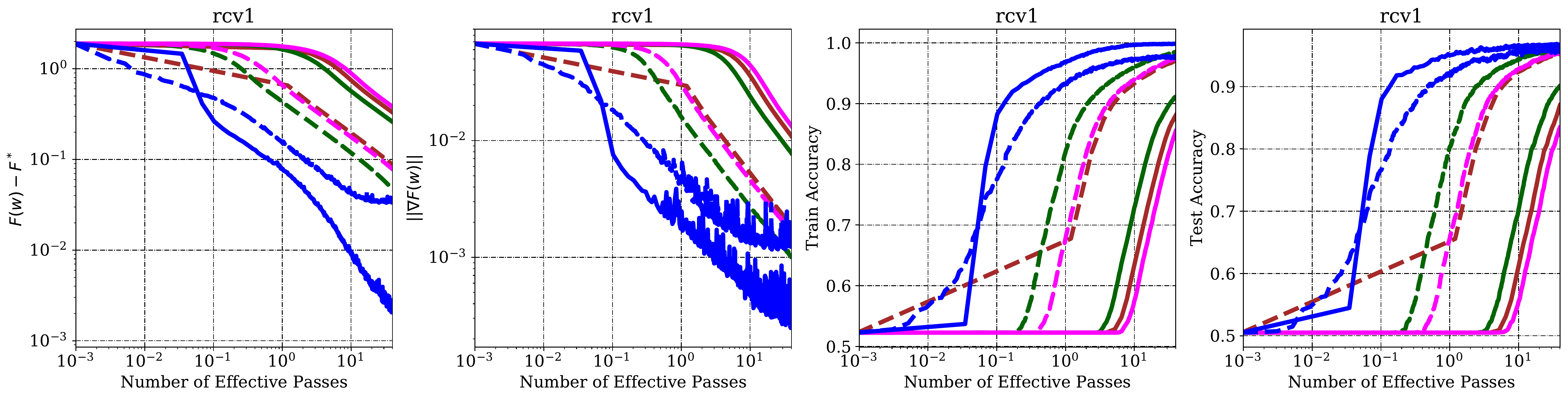}
	\caption{\texttt{rcv1}: Stochastic Logistic Regression ($\lambda = 10^{-6}$).}
	\label{fig:stochrcv1_train_6}
\end{figure*}
\clearpage
\begin{figure*}[ht]
{	\includegraphics[trim=10 100 10 110,clip, width=0.95\textwidth]{Figures/a1a_comp_stoch_NLLS_leg.pdf}
}
	\includegraphics[width=0.99\textwidth]{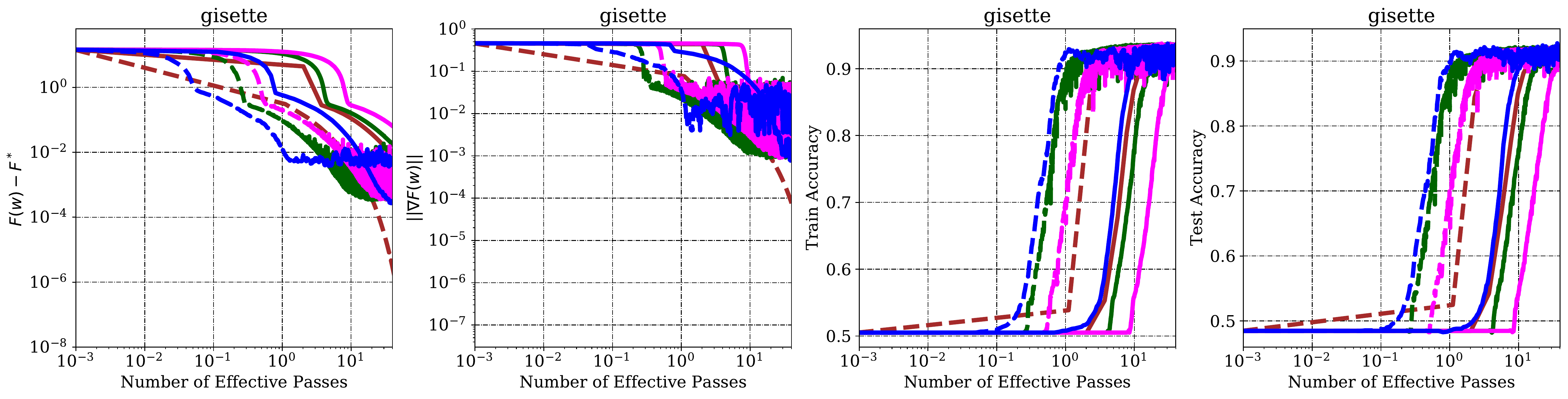}
	\caption{\texttt{gisette}: Stochastic Logistic Regression ($\lambda = 10^{-3}$).}
	\label{fig:gisette_scale_3}
\end{figure*}

\begin{figure*}[ht]
	\centering
	\includegraphics[width=0.99\textwidth]{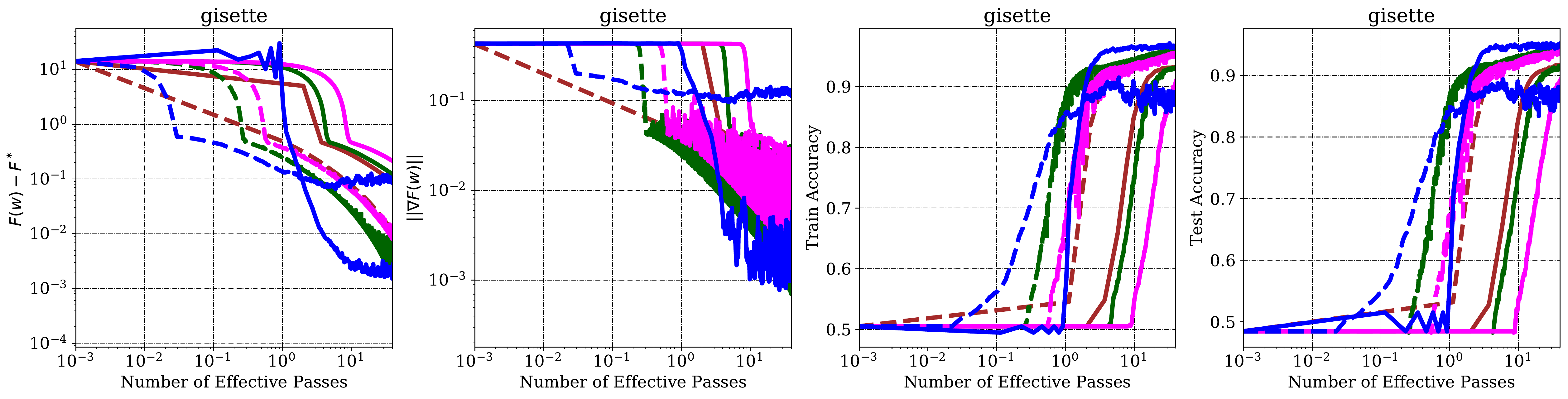}
	\caption{\texttt{gisette}: Stochastic Logistic Regression ($\lambda = 10^{-4}$).}
	\label{fig:stochgisette_scale_4}
\end{figure*}

\begin{figure*}[ht]
	\centering
	\includegraphics[width=0.99\textwidth]{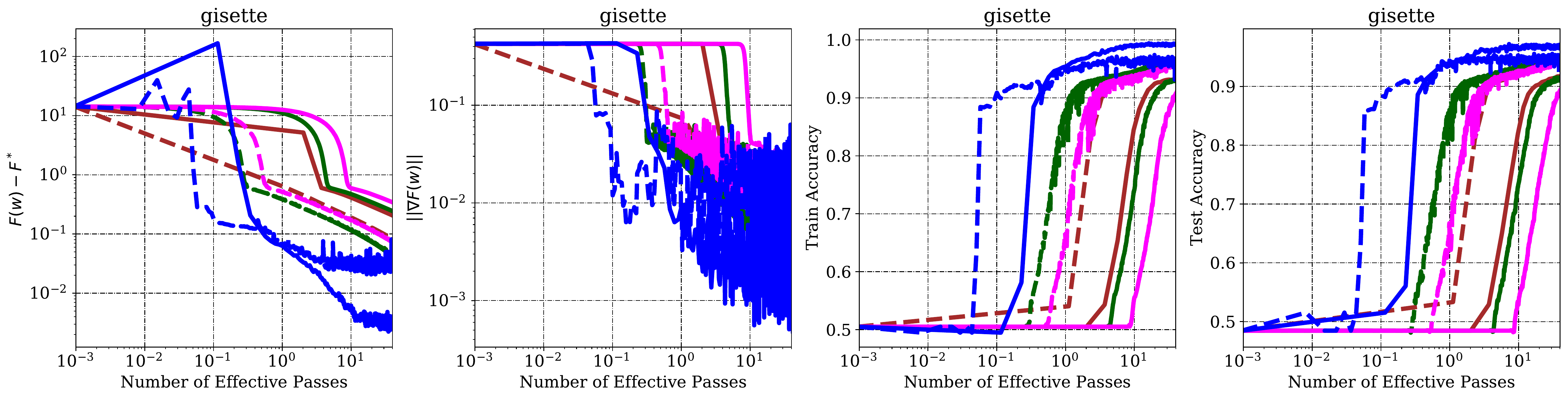}
	\caption{\texttt{gisette}: Stochastic Logistic Regression ($\lambda = 10^{-5}$).}
	\label{fig:stochgisette_scale_5}
\end{figure*}

\begin{figure*}[!h]
	\centering
	\includegraphics[width=0.99\textwidth]{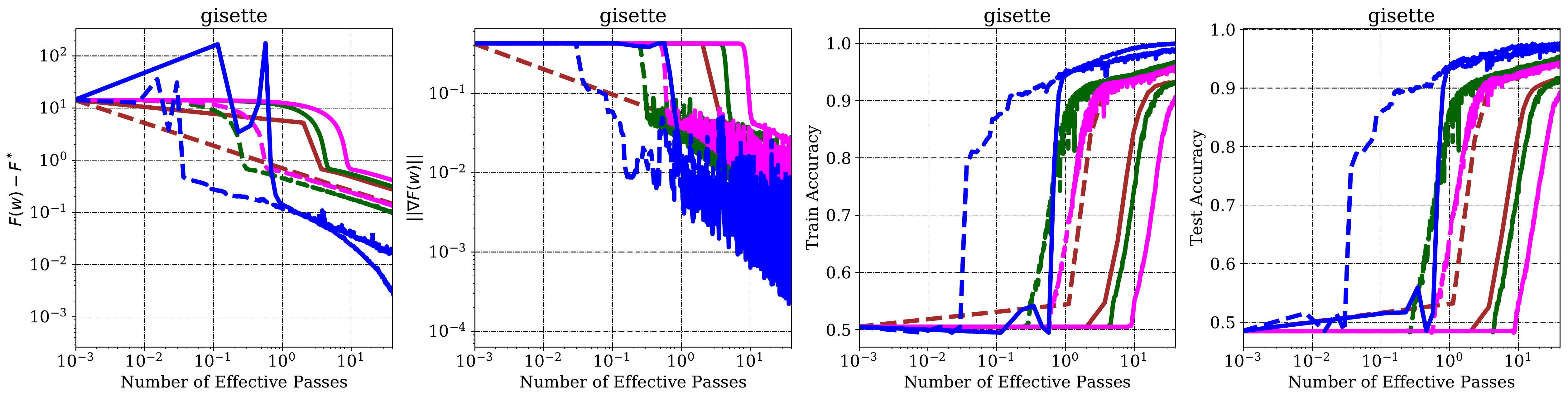}
	\caption{\texttt{gisette}: Stochastic Logistic Regression ($\lambda = 10^{-6}$).}
	
	\label{fig:stochgisette_scale_6}
\end{figure*}
\clearpage
\begin{figure*}[ht]
	\centering
{	\includegraphics[trim=10 100 10 110,clip, width=0.95\textwidth]{Figures/a1a_comp_stoch_NLLS_leg.pdf}
}
	\includegraphics[width=0.99\textwidth]{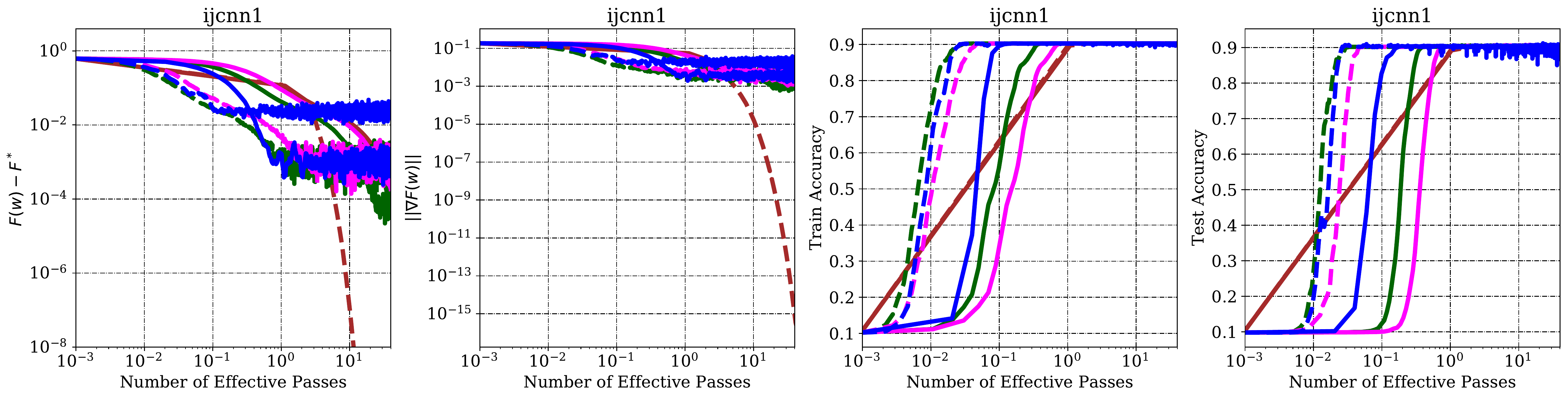}
	\caption{\texttt{ijcnn1}: Stochastic Logistic Regression ($\lambda = 10^{-3}$).}
	\label{fig:stochijcnn1_3}
\end{figure*}

\begin{figure*}[ht]
	\centering
	\includegraphics[width=0.99\textwidth]{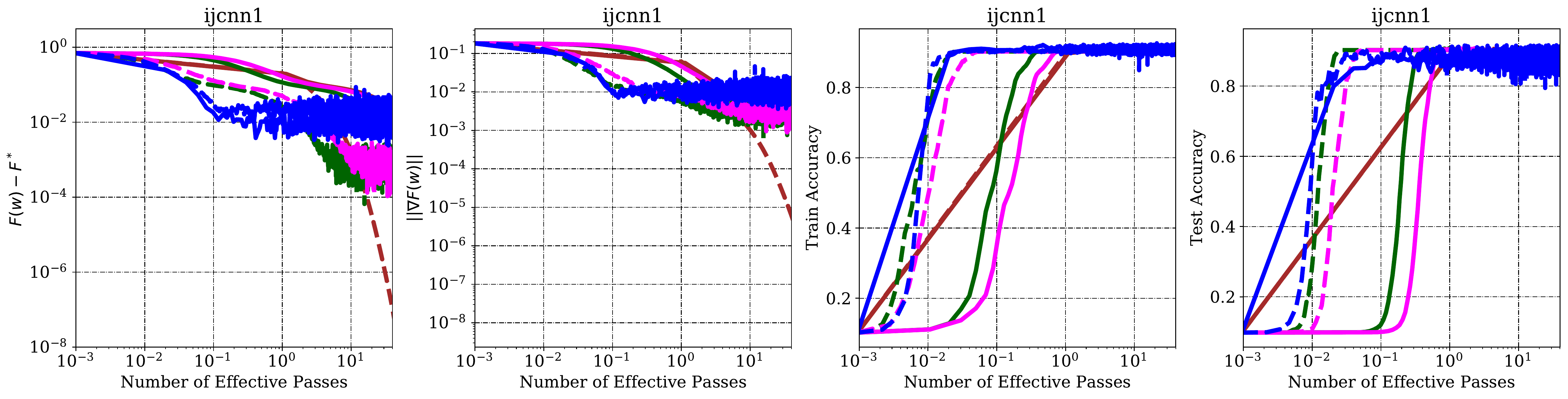}
	\caption{\texttt{ijcnn1}: Stochastic Logistic Regression ($\lambda = 10^{-4}$).}
	\label{fig:stochijcnn1_4}
\end{figure*}

\begin{figure*}[ht]
	\centering
	\includegraphics[width=0.99\textwidth]{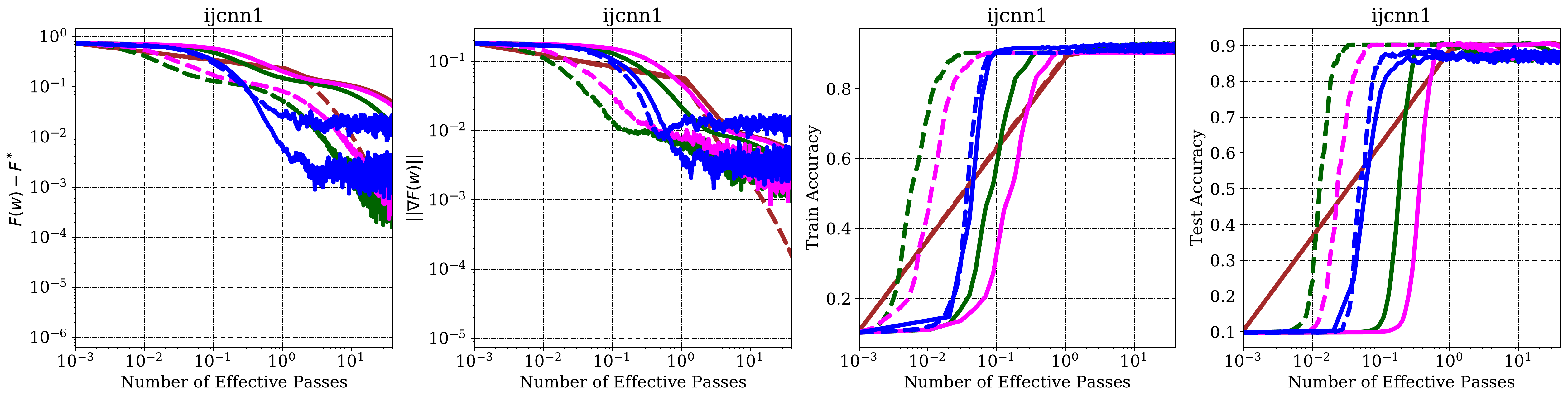}
	\caption{\texttt{ijcnn1}: Stochastic Logistic Regression ($\lambda = 10^{-5}$).}
	\label{fig:stochgisette_scale_5}
\end{figure*}

\begin{figure*}[!h]
	\centering
	\includegraphics[width=0.99\textwidth]{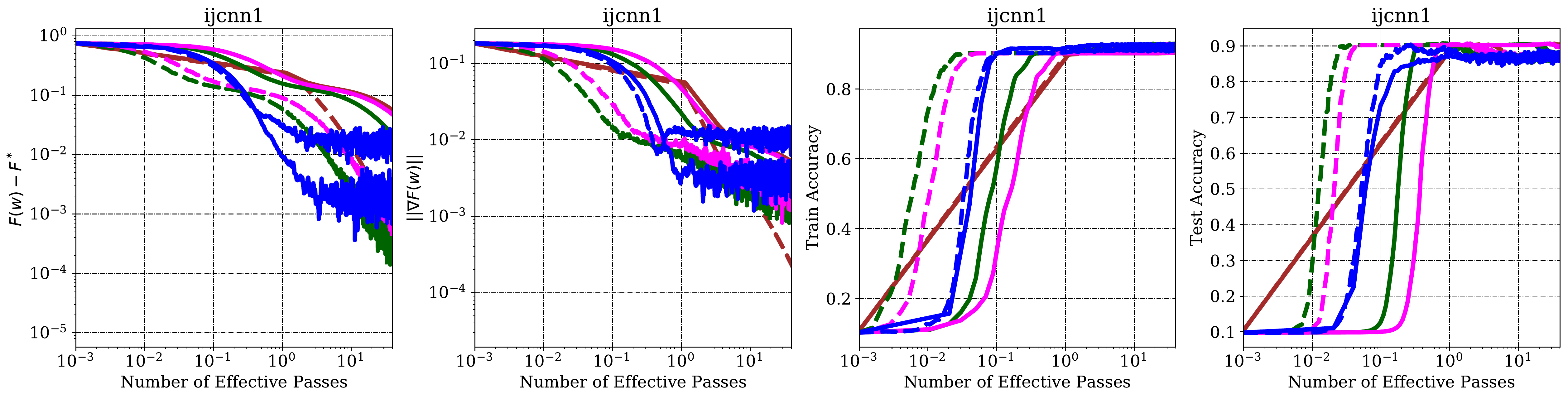}
	\caption{\texttt{ijcnn1}: Stochastic Logistic Regression ($\lambda = 10^{-6}$).}
	\label{fig:stochijcnn1_6}
\end{figure*}

\clearpage

\subsection{Additional Numerical Experiments: Stochastic Nonconvex Functions} \label{sec:moreStochNonCOnvResults}

\begin{figure*}[htb]
	\centering
{	\includegraphics[trim=10 100 10 110,clip, width=0.95\textwidth]{Figures/a1a_comp_stoch_NLLS_leg.pdf}
}
	\includegraphics[trim=10 10 10 10,clip,width=0.99\textwidth]{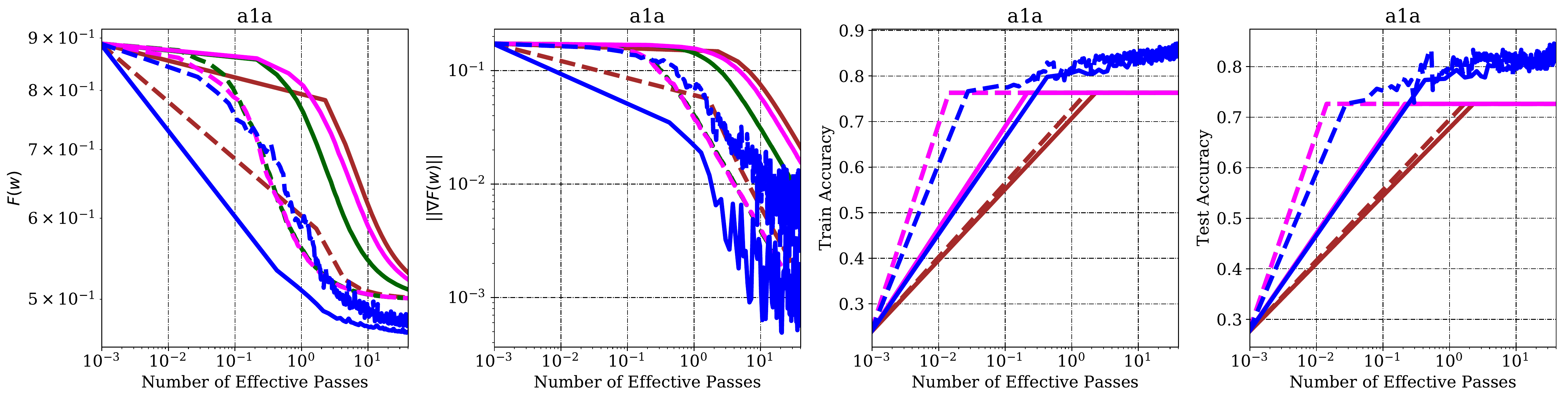}
	\caption{\texttt{a1a}: Stochastic Nonlinear Least Squares.}
	\label{fig:a1a_Stoch_NL_LS}
\end{figure*}

\begin{figure*}[!h]
	\centering
	\includegraphics[trim=10 10 10 10,clip,width=0.99\textwidth]{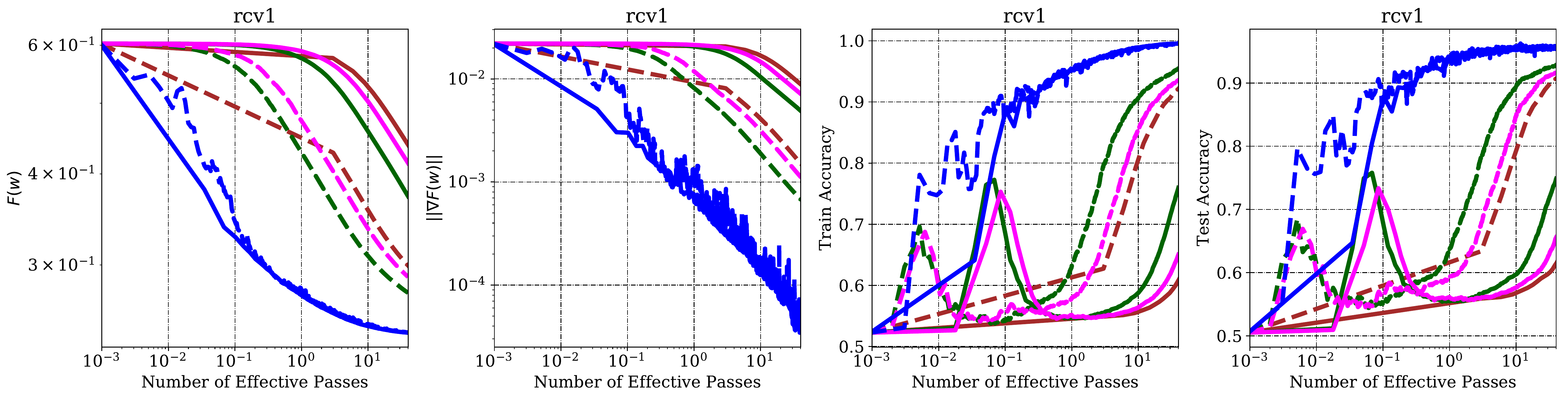}
	\caption{\texttt{rcv1}: Stochastic Nonlinear Least Squares.}
	\label{fig:rcv1_train_Stoch_NL_LS}
\end{figure*}

\begin{figure*}[htb]
	\centering
	\includegraphics[trim=10 10 10 10,clip,width=0.99\textwidth]{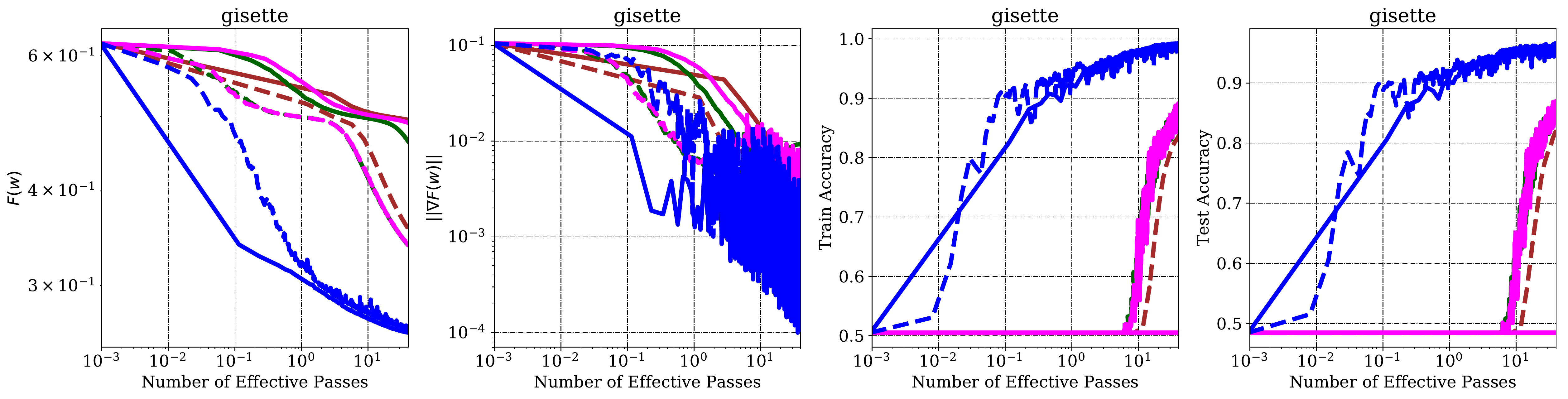}
	\caption{\texttt{gisette}: Stochastic Nonlinear Least Squares.}
	\label{fig:gisette_scale_Stoch_NL_LS}
\end{figure*}

\begin{figure*}[!h]
	\centering
	\includegraphics[trim=10 10 10 10,clip,width=0.99\textwidth]{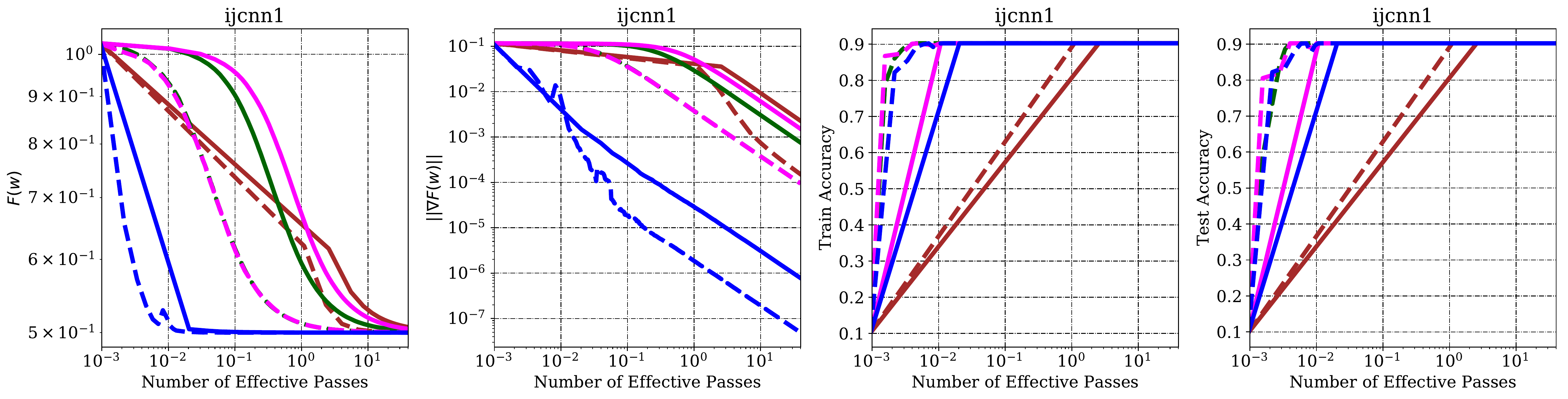}
	\caption{\texttt{ijcnn1}: Stochastic Nonlinear Least Squares.}
	\label{fig:ijcnn1_Stoch_NL_LS}
\end{figure*}

\end{document}